\documentclass[11pt,leqno,oneside]{amsart}
\usepackage[width=6.5in,height=8.8in]{geometry}
\usepackage[mathscr]{eucal}
\usepackage{amssymb, amsmath,array, amscd, bbm, comment, charter}
\usepackage{enumerate}
\usepackage{graphicx}
\usepackage{url}
\usepackage{tikz}
\usepackage[colorlinks,plainpages]{hyperref}

\newtheorem{thm}{Theorem}[section]
\newtheorem{defn}[thm]{Definition}
\newtheorem{lem}[thm]{Lemma}
\newtheorem{cor}[thm]{Corollary}
\newtheorem{pro}[thm]{Proposition}

\theoremstyle{definition}
\newtheorem{rem}[thm]{Remark}

\begin{document}

\title[Amoeba Measures of Random Plane Curves ]{Amoeba Measures of Random Plane Curves }

\author{Al\.i Ula\c{s} \"Ozg\"ur K\.i\c{s}\.isel}
\thanks{A. U. \"O K\.i\c{s}\.isel is partially supported by CIMPA Research in Pairs Fellowship, T\"{U}B\.{I}TAK 2219 grant and Camille Jordan Institute.} 
\address{Mathematics Department, Middle East Technical University, Ankara, Turkey and Camille Jordan Institute, Lyon, France}
\email{akisisel@metu.edu.tr}

\author{Jean-Yves Welschinger}
\thanks{ }
\address{CNRS, ICJ UMR5208, Ecole Centrale de Lyon, INSA Lyon, Universite Claude
Bernard Lyon 1, Université Jean Monnet,
69622 Villeurbanne, France}
\email{welschinger@math.univ-lyon1.fr}

\date{\today}


\begin{abstract}
We prove that the expected area of the amoeba of a complex plane curve of degree $d$ is less than $\displaystyle{3\ln(d)^2/2+9\ln(d)+9}$ and once rescaled by $\ln(d)^2$, is asymptotically bounded from below by $3/4$. In order to get this lower bound, given disjoint isometric embeddings of a bidisc of size $1/\sqrt{d}$ in the complex projective plane, we lower estimate the probability that one of them is a submanifold chart of a complex plane curve. It exponentially converges to one as the number of bidiscs grow to $+\infty$. 
\end{abstract}

\maketitle

\section{Introduction} 

Amoebas of algebraic varieties appear in a plethora of mathematical contexts, such as Hilbert's 16th problem, tropical geometry, complex analysis, dynamical systems and mirror symmetry \cite{GKZ, Mik1, Mik2, Mik3, Vir, FPT, Rul, ELMW, Rua, Aur}. If $V$ is an algebraic subvariety of $(\mathbb{C}^{\times})^{n}$, then its amoeba $\mathcal{A}(V)=Log(V)$ is defined to be its image in $\mathbb{R}^{n}$ under the coordinatewise logarithm map $Log: (z_{1},\ldots, z_{n})\in (\mathbb{C}^{\times})^{n}\rightarrow (\ln|z_{1}|,\ldots,\ln|z_{n}|)\in \mathbb{R}^{n}$.

Amoebas were defined by Gelfand, Kapranov and Zelevinsky in \cite{GKZ}, and named due to some of their geometric features, such as the convexity of the connected components of the complement in the case of hypersurface amoebas. In another direction, Passare and Rullg\r{a}rd \cite{PR} proved the remarkable fact that amoebas of  plane curves have Lebesgue areas bounded from above by $\pi^2 d^2/2$ in degree $d$. Moreover, this upper bound is sharp and the curves which realize the bound turn out to be real modulo some toric isometry and have real locus with rich topology and special placement with respect to the coordinate axes, they are called ``special Harnack curves'' \cite{Mik1, MR}. This finiteness result was generalized to volumes of amoebas of half-dimensional subvarieties in higher dimensions by Mikhalkin \cite{Mik2}, by making use of the rolled coamoeba map, which takes a point in $(\mathbb{C}^{\times})^{n}$ to $(S_{\pi}^{1})^{n}$ by sending each coordinate to its argument modulo $\pi$. 

In this paper, we adopt the viewpoint of random algebraic geometry and address the question of estimating the Lebesgue area of the amoeba of a typical complex plane algebraic curve, as opposed to an extremal one. Indeed, special Harnack curves are real curves having the maximal number of connected components in their degree, prescribed by the Harnack-Klein inequality. They become exponentially rare as their degree increase \cite{GW1}, and are thus far from being typical, from a probabilistic point of view. 

Some first steps regarding random amoebas were taken in \cite{BK}, where it was shown that with respect to the Kostlan distribution, the expected multiarea of the amoeba of a degree $d$ plane curve is equal to $\pi^2 d$, see subsection \ref{amoebameasures}. This implies that its expected Lebesgue area $\mathbb{E}_{d}(Vol(\mathcal{A}))$ gets bounded above by $\pi^2 d/2$, but the latter turns out to be much lower. 

\begin{thm} \label{mainbound} 
For all $d>0$, $\displaystyle{\mathbb{E}_{d}(Vol(A))\leq \frac{3}{2}\ln(d)^2+9\ln(d)+9}. $ 
Furthermore, as $d$ tends to infinity, 
\[\liminf_{d\rightarrow +\infty} \frac{\mathbb{E}_{d}(Vol(\mathcal{A}))}{\ln(d)^2} \geq \frac{3}{4}\cdot\]
\end{thm} 

This $\ln(d)^2$-asymptotic contrasts with the ``square-root law'' satisfied for instance by the number of solutions of a system of real polynomial equations, or by the volume or Betti numbers of real projective submanifolds, according to which the expected number grows as the square-root of the maximal number \cite{EK, SS, GW2, GW3, GW4, Let} with respect to the Kostlan distribution. The expected multivolume of plane amoebas obey the law though. 

In fact, the amoeba of a plane curve supports two different measures, its area measure $\lambda$ which is nothing but the Lebesgue measure of $\mathbb{R}^2$ restricted to $\mathcal{A}$ and its multiarea measure $\nu$, obtained by ``pulling back'' this Lebesgue measure to the curve and then pushing it forward  again to $\mathbb{R}^2$, both via the $Log$ map, so that the latter encodes as local weights the number of preimages of the $Log$ map restricted to the plane curve, see subsection \ref{amoebameasures}. In average, both measures have densities with respect to the Lebesgue measure, respectively denoted by $p$ and $\displaystyle{da/(2\pi)}$, where $p(t)$ denotes the probability that a plane curve intersects the compact torus $Log^{-1}(t), t\in \mathbb{R}^2$, whereas $a(t)$ denotes the area of the latter computed with respect to the Fubini-Study metric of the complex projective plane, see Proposition \ref{expectedmeasure} and subsection \ref{FSmetric}. The density $\displaystyle{da/(2\pi)}$ actually equals the average number of intersection points of a complex plane curve with the torus $Log^{-1}(t)$, by Crofton's formula, see Theorem \ref{crofton}. We deduce that $\displaystyle{p\leq \min\left(1,da/(4\pi)\right)}$, so that $\displaystyle{\mathbb{E}_{d}(Vol(\mathcal{A}))\leq \int_{\mathbb{R}^2}\min\left(1,da(t)/(4\pi)\right)dt}$, see Corollary  \ref{expbound}. The upper bound in Theorem \ref{mainbound} is obtained by estimating the latter integral and we prove in particular that as $d$ grows to $+\infty$, 
\begin{equation} \label{upasym}
\int_{\mathbb{R}^{2}}\min\left(1, da(t)/(4\pi)\right)dt = \frac{3}{2}\ln(d)^{2}+3\left(1+\ln\left(\pi\right)\right)\ln(d)+O(1),
\end{equation} 
see Theorem \ref{upperboundasymptotic} and Appendix \ref{B}. 

In order to lower estimate the density $p$ and deduce the asymptotic lower estimate in Theorem \ref{mainbound}, we use a different strategy. We define a {\it degree $d$ marked Fubini-Study chart of capacity $\kappa\in (0,1]$} centered at $z\in \mathbb{CP}^2$ to be the image of the bidisc $\mathbb{B}_{d}(\kappa)=D(0,1/\sqrt{6d})\times D(0,\kappa/\sqrt{6d})\subset \mathbb{C}^2\subset \mathbb{CP}^2$ under some isometry $\psi\in PU_{3}(\mathbb{C})$ such that $\psi(0)=z$ . It is a {\it submanifold chart} for a complex plane curve $V_{P}\subset \mathbb{CP}^2$ whenever $\psi^{-1}(V_{P})$ is the graph of a map $D(0,1/\sqrt{6d})\rightarrow D(0,\kappa/\sqrt{6d})$, see subsection \ref{allthecharts}. Given $N$ disjoint degree $d$ marked FS-charts of capacity $\kappa\in (0,1]$, we lower estimate the probability that one of them be a submanifold chart for a plane curve, namely,
\begin{thm} \label{exponential} 
For every $\kappa\in (0,1]$, there exists $d_{0}\in \mathbb{N}^{\times}$ such that for every $d\geq d_{0}$, every $N>0$ and every choice of degree $d$ marked Fubini-Study charts of capacity $\kappa$ centered at points $z_{1},\ldots,z_{N}\in \mathbb{CP}^2$ such that
\[ \min_{i\neq j} d_{FS}(z_{i},z_{j})\geq \sqrt{\frac{20\ln(d)}{d}}, \] 
the probability that one of these charts is a  submanifold chart for a degree $d$ complex plane curve is greater than $1-14\gamma(\kappa)^N$. 
\end{thm} 
In Theorem \ref{exponential}, $d_{FS}$ denotes the Fubini-Study distance in $\mathbb{CP}^2$, see subsection \ref{FSmetric}, and $\gamma$ the function 
\begin{equation} 
\gamma: \kappa\in (0,1] \mapsto \gamma(\kappa) = 1-\exp\left(-\frac{2^{11} 3^3}{\kappa^2}\right)(1-e^{-1/2}) \in \mathbb{R} 
\end{equation} 
Since $\gamma<1$, $\gamma^N$ decreases to zero as $N$ grows to $+\infty$. In order to lower estimate the density $p(t)$, $t=(t_{1},t_{2})\in \mathbb{R}^2$, and deduce Theorem \ref{mainbound} we manage to pack disjoint degree $d$ marked FS-charts of capacity $\kappa<1$ centered at points on the torus $Log^{-1}(t)$ in such a way that one of them being a submanifold chart for a complex plane curve $V_{P}$ guarantees the intersection $V_{P}\cap Log^{-1}(t)\neq \emptyset$. In order to be able to do so, we need $e^{t_{1}}$ and $e^{t_{2}}$ to be far away from the coordinate axes, so that the torus has low curvature and high diameter and we actually need $\kappa<1$ to get
\begin{thm} \label{lowerest} 
For every $\kappa\in (0,1)$, there exists $d_{0}\in \mathbb{N}^{\times}$ such that for every $d\geq d_{0}$, every $\displaystyle{1>\delta\geq 3\sqrt{\frac{5\ln(d)}{d}}}$ and every $t=(t_{1},t_{2})\in \mathbb{R}^2$ such that $\delta\leq e^{t_{1}}\leq 1$ and $\delta\leq e^{t_{2}}\leq 1$, 
\[  p(t)= \mathbb{P}(Log^{-1}(t)\cap V_{P}\neq \emptyset) \geq 1-14\gamma(\kappa)^{N}, \] 
with $\displaystyle{ N=\left\lfloor \frac{2\delta \sqrt{d}}{3\sqrt{5\ln(d)}} \right\rfloor ^{2} }$. 
\end{thm} 
We deduce that the density $p$ converges to $1$ on $\mathbb{R}^2$ in the sense of distributions, see Corollaries \ref{liminf} and \ref{distributions}. 

Let us point out that the disjoint charts in Theorem \ref{exponential} do not lead to independent events, even though they become so in the limit, at least heuristically. We need to set a probabilistic result to overcome this issue, Theorem \ref{probmain} and  Corollary \ref{cormain}, which actually represents a considerable part of the proof. 

The organization of the paper is as follows. Section \ref{preliminaries} introduces the Fubini-Study metric together with the Bergman kernel and gathers all the material we need throughout the paper, as preliminaries. Section \ref{upper} introduces the amoeba measures $\lambda$ and $\nu$, sets their first properties and contains the proofs of the asymptotic upper bounds given by \eqref{upasym} and Theorem \ref{mainbound}, see Theorem \ref{upperboundasymptotic} and Corollary \ref{upperlimit}. Section \ref{barrier} introduces marked FS-charts and sets an effective criterion for being a submanifold chart, in the spirit of the barrier method \cite{NS, GW2}, in terms of $L^2$ ball-norms of complex polynomials. Section  \ref{lower} develops the probabilistic tool we use to overcome the lack of independence which we quantify, leading to  Theorem \ref{probmain} and  Corollary \ref{cormain}. It then contains the proofs of Theorems \ref{exponential}, \ref{lowerest} and the lower estimate part of Theorem \ref{mainbound}. The paper ends with an Appendix containing a proof of Crofton's formula for the reader's convenience, see Theorem \ref{crofton}, and a proof of the first part of Theorem \ref{mainbound} in a slightly stronger form, see Theorem \ref{upperboundall}  and Remark \ref{last}. 

\section{Preliminaries} \label{preliminaries} 
Let us first recall all the material we will need throughout the paper. 

\subsection{Hermitian matrices} \label{Hermitian} 
Let $\displaystyle{C=(C_{ij})_{1\leq i,j\leq N}}$ be a Hermitian matrix of size $N$. Let $||C||_{\infty}=\sup_{1\leq i,j\leq N} |C_{ij}|$ and  $||C||$ denote its norm as an operator on $\mathbb{C}^{N}$ equipped with its standard Hermitian inner product, so that
\[  ||C||=\sup_{v\in \mathbb{C}^{N}\setminus \{\mathbf{0}\}} \frac{ ||Cv||}{||v||}. \]
It coincides with the largest eigenvalue of $C$ in absolute value, the latter being all real. 

\begin{lem} \label{firstlemma}
Let $C$ be a positive definite Hermitian matrix of size $N$. Then, for every $a\in \mathbb{C}^{N}$, $\displaystyle{\langle a, C^{-1}a\rangle \geq \frac{ ||a||^2}{||C||}}$ and  $\det(C)\geq ||C^{-1}||^{-N}$.  If moreover one of the diagonal entries of $C$ equals $1$, then $\min(||C||,||C^{-1}||)\geq 1$. 
\end{lem} 

\noindent \textit{Proof:} Let us write $C=A^2$ for some positive definite Hermitian matrix $A$ of size $N$. Then, for every $a\in \mathbb{C}^{N}$, $||a||^2=||AA^{-1}a||^2\leq  ||A||^2 ||A^{-1}a||^2=||C||\langle a, C^{-1}a\rangle$. 
Moreover, if $c_{1},\ldots,c_{N}\in \mathbb{R}^{\times}_{+}$ denote the eigenvalues of $C$,  
\[ \det(C)=\prod_{i=1}^{N} c_{i} \geq \left(\inf_{1\leq i\leq N} c_{i}\right)^{N}, \] 
whereas $\displaystyle{||C^{-1}||=\sup_{1\leq i\leq N}c_{i}^{-1} = \left(\inf_{1\leq i\leq N} c_{i}\right)^{-1}}$.  

Let us now assume that a diagonal entry $C_{ii}=1$ for some $i\in \{1,\dots,N\}$. Let $\{e_{1},\ldots,e_{N}\}$ be the standard orthonormal basis of $\mathbb{C}^{N}$ and $(\widetilde{e_{j}})_{1\leq j\leq N}$ be an orthonormal basis of eigenvectors for $C$. Then $\displaystyle{e_{i}=\sum_{j=1}^{N} a_{j}\widetilde{e_{j}}}$ for some $(a_{j})_{1\leq j\leq N}\in \mathbb{C}^{N}$ and $\displaystyle{1=\langle C(e_{i}),e_{i}\rangle= \sum_{j=1}^{N} c_{j} |a_{j}|^2}$. Since $\displaystyle{\sum_{j=1}^{N} |a_{j}|^2=1}$, we deduce that $\inf_{1\leq j\leq N} c_{j} \leq 1 \leq \sup_{1\leq j\leq N} c_{j}$, so that $\min(||C||, ||C^{-1}||)\geq 1$.  \hfill $\Box$ 

\begin{lem} \label{normest}
Let $C$ be a Hermitian matrix and $I$ the identity matrix of size $N$. Then, 
\begin{enumerate} [(a)] 
\item  $||C||_{\infty}\leq ||C|| \leq N ||C||_{\infty}$. 
\item If $||C-I||<1$, then $\displaystyle{ ||C^{-1}|| \leq \frac{1}{1-||C-I||}}$ and $\displaystyle{ ||C^{-1}-I|| \leq \frac{||C-I||}{1-||C-I||}\cdot }$ 
\end{enumerate} 
\end{lem} 

\noindent \textit{Proof:}  (a) Let $\{e_{i}\}_{1\leq i\leq N}$ denote the standard basis of $\mathbb{C}^{N}$. Then, for all $1\leq i,j\leq N$ 
\[ |C_{ij}|=|\langle C(e_{j}),e_{i}\rangle| \leq ||C|| \qquad \mathrm{by} \, \mathrm{Cauchy-Schwarz}, \] 
which implies that $||C||_{\infty}\leq ||C||$. On the other hand, for any $v=\sum_{i=1}^{N}v_{i}e_{i}\in \mathbb{C}^{N}$, $C(\sum_{i=1}^{N}v_{i}e_{i})=\sum_{i=1}^{N} v_{i}(\sum_{j=1}^{N} C_{ji}e_{j})$, so that by Cauchy-Schwarz again,
\[ ||C(v)||^2= ||C(\sum_{i=1}^{N}v_{i}e_{i})||^2\leq \sum_{i=1}^{N} \sum_{k=1}^{N} |v_{i}||v_{k}| \sum_{j=1}^{N} |C_{ji}||C_{jk}|\leq N^{2}||C||_{\infty}^{2}||v||^2. \]
The second inequality follows. 

(b) Writing $C^{-1}=I+(C^{-1}-I)$, we get $||C^{-1}||\leq 1+||C^{-1}-I||$ by triangle inequality and thus 
\begin{eqnarray*} 
1&\geq& ||C^{-1}||-||C^{-1}-I||  \\
&=& ||C^{-1}||(1-||I-C||).
\end{eqnarray*} 
The first inequality follows and the second one as well, upon multiplying both sides by $||C-I||$. \hfill $\Box$

\subsection{Conditional expectation} \label{conditionalexpectation}

For any event $A$ in any sample space with probability measure $\mu$ and for any complex random variable $Z$, let 
\[ \mathbb{E}_{A}(Z)=\frac{1}{\mu(A)} \int_{A} Z(P)d\mu(P) \]
denote the conditional expectation of $Z$ over $A$. 

\begin{lem} \label{expectedbound}
Let $Z$ be a real random variable with density, such that $\mathbb{P}(Z\geq z)=o\left(1/z\right)$ and $\mathbb{P}(Z\geq z)$ is integrable. Then, for every event $A$ and every $x< \mathbb{E}_{A}(Z)$,  
\[ \mathbb{P}(A)\leq \frac{ \int_{x}^{+\infty} \mathbb{P}(Z\geq z)dz}{\mathbb{E}_{A}(Z)-x} \cdot \] 
\end{lem} 

\noindent \textit{Proof:} Since the law of $Z$ has a density with respect to the Lebesgue measure of $\mathbb{R}$, there exists $x_{A}\in \mathbb{R}$ such that $\mathbb{P}(A)=\mathbb{P}(Z\geq x_{A})$. Then,
\begin{eqnarray*} 
\int_{A}Zd\mu &=& \int_{A\cap\{Z<x_{A}\}} Zd\mu + \int_{A\cap \{Z\geq x_{A}\}} Zd\mu  \\ 
&\leq& x_{A}\mathbb{P}(A\cap \{Z<x_{A}\}) + \int_{A\cap \{Z\geq x_{A}\}} Z d\mu \\ 
&\leq& x_{A}\left(\mathbb{P}(A)-\mathbb{P}(A\cap \{Z\geq x_{A}\})\right)+\int_{A\cap \{Z\geq x_{A}\}} Zd\mu \\ 
&\leq& x_{A}\mathbb{P}(A^{c}\cap \{Z\geq x_{A}\})+\int_{A\cap \{Z\geq x_{A}\}} Zd\mu \quad \mathrm{since} \, \mathbb{P}(A)=\mathbb{P}(Z\geq x_{A})\\ 
&\leq& \int_{Z\geq x_{A}} Zd\mu \\ 
&=& \int_{x_{A}}^{+\infty} z\mu_{Z}(z), 
\end{eqnarray*} 
where $\mu_{Z}=Z_{*}\mu$ denotes the law of $Z$ and $A^{c}$ the contrary event of $A$. Integration by parts gives 
\begin{eqnarray*} 
\int_{x_{A}}^{+\infty} z\mu_{Z}(z) &=& -z\mathbb{P}(Z\geq z)|_{x_{A}}^{+\infty}+\int_{x_{A}}^{+\infty} \mathbb{P}(Z\geq z)dz \\ 
&=& x_{A}\mathbb{P}(A)+\int_{x_{A}}^{+\infty} P(Z\geq z)dz  \quad \mathrm{by}\,\,\mathrm{hypothesis,} \,\mathrm{since}\,\mathbb{P}(Z\geq x_{A})=\mathbb{P}(A). 
\end{eqnarray*} 
We deduce 
\[ \mathbb{E}_{A}(Z)\leq x_{A}+\frac{1}{\mathbb{P}(A)} \int_{x_{A}}^{+\infty} \mathbb{P}(Z\geq z)dz. \]
The function $\displaystyle{f(x)= x+\frac{1}{\mathbb{P}(A)} \int_{x}^{+\infty} \mathbb{P}(Z\geq z)dz}$ has derivative $\displaystyle{f^{\prime}(x)=1-\frac{\mathbb{P}(Z\geq x)}{\mathbb{P}(A)}}$. Hence, $f^{\prime}$ is negative on $(-\infty,x_{A})$ and positive on $(x_{A},+\infty)$, so that $f$ reaches a global minimum at $x_{A}$.  Therefore, for all $x\in \mathbb{R}$, 
\[ \mathbb{E}_{A}(Z)\leq x+\frac{1}{\mathbb{P}(A)} \int_{x}^{+\infty} \mathbb{P}(Z\geq z)dz. \]
Hence the result. \hfill $\Box$

\begin{rem} 
If the real random variable with density $Z$ posseses a variance, then the Bienaym\'{e}-Tchebychev inequality implies that for every $z> \mathbb{E}(Z)$, 
\begin{equation} \label{Bienayme}
 \mathbb{P}(Z\geq z)\leq  \mathbb{P}\left((Z-\mathbb{E}(Z))^2\geq (z-\mathbb{E}(Z))^2\right) \leq  \frac{V(Z)}{(z-\mathbb{E}(Z))^2},  
\end{equation} 
so that $\mathbb{P}(Z\geq z)=o(1/z)$, $\mathbb{P}(Z\geq z)$ is integrable and the hypotheses of Lemma \ref{expectedbound} get satisfied.  
\end{rem}

\begin{lem} \label{variance}
Let $Z$ be a real random variable with density, which possesses a variance $V(Z)$. Then, for all events $A$ such that $\mathbb{E}(Z)< \mathbb{E}_{A}(Z)$,
\[ \mathbb{P}(A)\leq \frac{4V(Z)}{\left(\mathbb{E}_{A}(Z)-\mathbb{E}(Z)\right)^2}\cdot \] 
\end{lem} 

\noindent \textit{Proof:}   
For every $x> \mathbb{E}(Z)$,  
\begin{eqnarray*} 
\int_{x}^{+\infty} \mathbb{P}(Z\geq z)dz &\leq& \int_{x}^{+\infty} \frac{V(Z)}{(z-\mathbb{E}(Z))^2} dz \quad \mathrm{by} \, \eqref{Bienayme}\\ 
&\leq& V(Z)\int_{x-\mathbb{E}(Z)}^{+\infty} \frac{1}{z^2}dz \\ 
&\leq& \frac{V(Z)}{x-\mathbb{E}(Z)}\cdot
\end{eqnarray*} 
Lemma \ref{expectedbound} then implies that for all $x\in(\mathbb{E}(Z),\mathbb{E}_{A}(Z))$, 
\[ \mathbb{P}(A)\leq \frac{V(Z)}{(x-\mathbb{E}(Z))(\mathbb{E}_{A}(Z)-x)}\cdot \] 
Choosing $\displaystyle{ x=\frac{\mathbb{E}(Z)+\mathbb{E}_{A}(Z)}{2}}$ provides the result. \hfill $ \Box$ 

\subsection{Gaussian complex vectors} 

\begin{defn} \label{normallaw}
For every positive definite Hermitian matrix $B$ of size $N$, the centered normal law $\mathcal{N}_{\mathbb{C}^{N}}(0,B)$ 
of $\mathbb{C}^{N}$ with covariance $B$ is the Gaussian law with density 
\[ \frac{e^{-\langle \mathbf{z}, B^{-1}\mathbf{z}\rangle}}{\pi^{N}\det B}\] 
with respect to the Lebesgue measure, where $\langle,\rangle$ denotes the standard Hermitian product on $\mathbb{C}^{N}$. 
\end{defn} 
Let us write $\mathbf{Z}\sim \mathcal{N}_{\mathbb{C}^{N}}(0,B)$ to denote that the random vector $\mathbf{Z}$ follows the law $\mathcal{N}_{\mathbb{C}^{N}}(0,B)$ and first recall the moment generating function of $||\mathbf{Z}||^2$.

\begin{lem} \label{centeredgaussian}
If $\mathbf{Z}\sim \mathcal{N}_{\mathbb{C}^{N}}(0,B)$, then for all $\displaystyle{0<\lambda<\frac{1}{||B||}}$,
$\displaystyle{ \mathbb{E}(e^{\lambda ||\mathbf{Z}||^2})=\frac{1}{\det(I-\lambda B)}\cdot}$ 
\end{lem} 

\noindent \textit{Proof:} Assuming $\displaystyle{0<\lambda<\frac{1}{||B||}}$, 
\begin{eqnarray*} 
\mathbb{E}(e^{\lambda ||\mathbf{Z}||^2})&=& \frac{1}{\pi^N \det(B)} \int_{\mathbb{C}^{N}} e^{\lambda ||\mathbf{z}||^2} e^{-\langle \mathbf{z},B^{-1}\mathbf{z}\rangle}d\mathbf{z} \\
&=&   \frac{1}{\pi^N \det(B)} \int_{\mathbb{C}^{N}} e^{-\langle \mathbf{z},(B^{-1}-\lambda I)\mathbf{z}\rangle}d\mathbf{z} \\
&=&  \frac{1}{\pi^N \det(B)|\det(A)|^2}\int_{\mathbb{C}^{N}} e^{-||\mathbf{z}^{\prime}||^2}d\mathbf{z}^{\prime} \\
&=& \frac{1}{\det(I-\lambda B)}\cdot
\end{eqnarray*} 
In this derivation, we write $B^{-1}-\lambda I=A^{*}A$ and perform the change of variables $\mathbf{z}^{\prime}=A\mathbf{z}$ when passing from the second to the third line. \hfill $\Box$ 

If $\displaystyle{\mathbf{Z}=(Z_{1},\ldots,Z_{N})\sim \mathcal{N}_{\mathbb{C}^{N}}(0,B)}$, then the complex random variables $Z_{i}$ satisfy $\displaystyle{(\mathbb{E}(Z_{i}\overline{Z_{j}}))_{1\leq i,j\leq N}=B}$ and the pairing 
\[ (\mathbf{Z},\mathbf{Z}^{\prime})\in \mathbb{C}^{N}\times \mathbb{C}^{N}\mapsto \langle \mathbf{Z}, B^{-1}\mathbf{Z}^{\prime}\rangle \in \mathbb{C} \] 
defines a Hermitian product on $\mathbb{C}^{N}$. Conversely, 

\begin{defn} \label{Gaussiandef}
The Gaussian measure of an $N$-dimensional complex Hermitian space $(E,\langle,\rangle)$ is the measure having density $\displaystyle{ \frac{e^{-\langle \mathbf{Z},\mathbf{Z}\rangle}}{\pi^{N}}}$ at the point $\mathbf{Z}\in E$ with respect to its Lebesgue measure. 
\end{defn} 

Given any basis $\{P_{1},\ldots,P_{N}\}$ of $E$ equipped with the Gaussian measure $\mu$ associated to its Hermitian product, the coordinates $X_{1}(P),\ldots, X_{N}(P)$ of a vector $P\in E$ in this basis become complex random variables. Let us recall their covariance matrix. 

\begin{lem}  \label{covariance} 
Let $E$ be an $N$-dimensional complex Hermitian space equipped with a basis $\{P_{1},\ldots,P_{N}\}$ and its Gaussian measure. Let $C=(\langle P_{i},P_{j}\rangle)_{1\leq i,j\leq N}$ and for every $i\in \{1,\ldots,N\}$, $X_{i}(P)$ be the $i$-th coordinate of $P\in E$ in the basis $\{P_{1},\ldots,P_{N}\}$. Then, 
$\displaystyle{ (Cov(X_{i},X_{j}))_{1\leq i,j\leq N}= \overline{C^{-1}}. }$ 
\end{lem} 

\noindent \textit{Proof:}   The random variables $X_{i}$ are centered, hence $Cov(X_{i},X_{j})=\mathbb{E}(X_{i}\overline{X_{j}})$ for all $1\leq i,j\leq N$,  with 
\[
\mathbb{E}(X_{i}\overline{X_{j}})= \frac{1}{\pi^{N}} \int_{E} X_{i}(P)\overline{X_{j}(P)}  d\mu(P).
\]
Let  $\{ \widetilde{P}_{1},\ldots,\widetilde{P}_{N}\}$ be an orthonormal basis of $E$ and let us write $P=\sum_{i=1}^{N}b_{i}\widetilde{P}_{i}$. We get 
\begin{eqnarray*} 
\mathbb{E}(X_{i}\overline{X_{j}})&=& \frac{1}{\pi^{N}} \int_{\mathbb{C}^{N}} \left(\sum_{l=1}^{N} b_{l} X_{i}(\widetilde{P}_{l})\right) \left(\sum_{m=1}^{N} \overline{b_{m}}\overline{X_{j}(\widetilde{P}_{m})}\right) e^{-||\mathbf{b}||^2}dVol(\mathbf{b}) \\ 
&=& \sum_{l=1}^{N}X_{i}(\widetilde{P}_{l})\overline{X_{j}(\widetilde{P}_{l})},
\end{eqnarray*} 
since $\displaystyle{ \frac{1}{\pi^{N}} \int_{\mathbb{C}^{N}} b_{l} \overline{b_{m}} e^{-||\mathbf{b}||^2} dVol(\mathbf{b})=\delta_{lm}}$. Let $\mathbf{X}=(X_{1},\ldots,X_{N}):E\rightarrow \mathbb{C}^{N}$ and let $\mathbf{X}^{\star}$ denote its adjoint. We deduce that 
\begin{equation} \label{XXstar}
 (\mathbb{E}(X_{i}\overline{X_{j}}))_{1\leq i,j\leq N} =\mathbf{X}\mathbf{X}^{\star},  
\end{equation} 
while $\displaystyle{(\mathbf{X}^{-1})^{\star}\mathbf{X}^{-1}=(\overline{\langle P_{i},P_{j}\rangle})_{1\leq i,j\leq N}}$. Hence the result. \hfill $\Box$ 

We will finally need the following lower estimate of the measure of large balls in  complex Gaussian vector spaces.

\begin{lem} \label{lemma18} 
Let $E$ be an $N$-dimensional complex Hermitian vector space equipped with its associated Gaussian measure. Then, for all $R>\sqrt{N}$, 
\[ \mathbb{P}(||P||\geq R) \leq \frac{R^{2N}}{N^N} \exp(-R^2+N). \] 
\end{lem} 

\noindent \textit{Proof:} Let us fix an orthogonal decomposition $\displaystyle{E=\sum_{i=1}^{N} E_{i}}$ into complex lines and write $\displaystyle{P=\sum_{i=1}^{N} P_{i}}$ with $P_{i}\in E_{i}$. Using Bernstein's trick, we observe that for all $0<\lambda<1$,
\[\quad ||P||\geq R \Leftrightarrow  e^{\lambda \sum_{i=1}^{N} ||P_{i}||^2 }\geq e^{\lambda R^2},  \] 
so that
\begin{eqnarray*} 
\mathbb{P}(||P||\geq R) &\leq &  \frac{\mathbb{E}\left(e^{\lambda ||P||^2}\right)}{e^{\lambda R^2}} \quad \mathrm{by} \, \mathrm{Markov's} \, \mathrm{inequality}\\ 
&\leq& \frac{ \prod_{i=1}^{N} \mathbb{E}\left(e^{\lambda ||P_{i}||^2}\right)}{e^{\lambda R^2}} \\ 
&\leq& \frac{e^{-\lambda R^2}}{(1-\lambda)^{N}} \quad \mathrm{by}\, \mathrm{Lemma} \, \ref{centeredgaussian}.
\end{eqnarray*} 
The result is obtained by choosing $\displaystyle{\lambda=1-\frac{N}{R^2}}$.  \hfill $\Box$

\subsection{Fubini-Study metric} \label{FSmetric} 
Recall that the tautological line bundle $\tau=\mathcal{O}_{\mathbb{CP}^{n}}(-1)$ on $\mathbb{CP}^{n}$ inherits the Fubini-Study Hermitian metric $h_{FS}$ from the Hermitian product on the trivial complex vector bundle of rank $n+1$. Consider the affine chart 
\[ U_{0}=\{[X_{0}:X_{1}:\ldots:X_{n}]\in \mathbb{CP}^{n}| X_{0}\neq 0\}\] 
and  the local section $e\in H^{0}(U_{0},\tau)$ defined by
\[ e([X_{0}:\ldots:X_{n}])=\left(1,\frac{X_{1}}{X_{0}},\ldots,\frac{X_{n}}{X_{0}}\right) \in \tau|_{[X_{0}:\ldots:X_{n}]}. \] 
We identify $z=(z_{1},\ldots,z_{n})\in \mathbb{C}^{n}$ with $[1:z_{1}:\ldots:z_{n}]\in U_{0}$ and, for all $d\in \mathbb{Z}$,  denote by $h_{d}$ the induced hermitian metric on the line bundle $L^d=\mathcal{O}_{\mathbb{CP}^{n}}(d)=(\tau^{*})^d$. Thus, for every $z\in \mathbb{C}^{n}\cong U_{0}$, 
\[
||e(z)||^{2}_{FS}=1+||z||^2 \\ 
\]
and every degree $d$ homogeneous polynomial $P\in H^{0}(\mathbb{CP}^{n}, \mathcal{O}_{\mathbb{CP}^{n}}(d))$ inherits the pointwise Fubini-Study norm,
\begin{equation} \label{sectionnorm} 
 ||P(z)||^{2}_{FS}=h_{d}(P,P)=\frac{|P\circ e(z)|^2}{||e(z)||_{FS}^{2d}}=\frac{|P(1,z)|^2}{(1+||z||^2)^d}. 
\end{equation} 
The curvature form $\omega_{FS}$ of this Fubini-Study metric on $L=\tau^{*}$ reads
\[ -\frac{1}{2i}\partial\overline{\partial} \log \left(1+\sum_{i=1}^{n}z_{i}\overline{z_{i}}\right)=-\frac{1}{2i} \partial\left(\frac{\sum_{i=1}^{n} z_{i}d\overline{z_{i}}}{1+\sum_{i=1}^{n} z_{i}\overline{z_{i}}}\right)\cdot\] 
Consequently, 
\begin{eqnarray} 
\forall z\in \mathbb{C}^{n}\cong U_{0}, \omega_{FS|_{z}}&=&-\frac{1}{2i} \frac{\sum_{i=1}^{n} dz_{i}\wedge d\overline{z_{i}} }{1+||z||^2} - \frac{1}{2i} \frac{ (\sum_{i=1}^{n} z_{i} d\overline{z_{i}})\wedge (\sum_{i=1}^{n} \overline{z_{i}}dz_{i})}{(1+||z||^2)^2}  \quad \mathrm{and} \label{FS1} \\ 
\forall h\in T_{|z}U_{0}, ||h||^2_{FS}&=&\omega_{FS}(h,ih)=\frac{||h||^2}{1+||z||^2}-\frac{|\langle h,z\rangle |^2}{(1+||z||^2)^2}\cdot  \label{FSexplicit}
\end{eqnarray} 
The Fubini-Study volume form $dVol_{FS}$ of $\mathbb{C}^{n}\cong U_{0}$ and its Euclidean volume form $dVol$ satisfy the following relation. 

\begin{lem} \label{volcompare}
For every $z\in \mathbb{C}^{n}\cong U_{0}$, 
$\displaystyle{ dVol_{FS_{|z}}=\frac{dVol}{(1+||z||^2)^{n+1}}\cdot }$ 
In particular, $\displaystyle{Vol_{FS}(\mathbb{CP}^{n})=\frac{\pi^{n}}{n!}\cdot}$
\end{lem} 
In fact, $\mathbb{CP}^{n}$ deprived of any projective hyperplane gets symplectomorphic to the open unit ball of $\mathbb{C}^{n}$. 

\noindent \textit{Proof:} For all $z\in \mathbb{C}^{n}\cong U_{0}$, 
\begin{eqnarray*} 
dVol_{FS_{|z}}&=&  \frac{1}{n!} \omega^{n}_{FS_{|z}} \\
&=& \frac{dVol}{(1+||z||^2)^{n}} + \frac{(-1)^{n}n(\sum_{i=1}^{n} dz_{i} \wedge d\overline{z_{i}})^{n-1}\wedge (\sum_{i=1}^{n} z_{i}d\overline{z_{i}})\wedge (\sum_{i=1}^{n} \overline{z_{i}}dz_{i}) }{n!(2i)^{n}(1+||z||^2)^{n+1} } \quad \mathrm{by} \, \eqref{FS1} \\ 
&=& \frac{dVol}{(1+||z||^2)^n} - \frac{||z||^2 dVol}{(1+||z||^2)^{n+1}}  \\ 
&=& \frac{dVol}{(1+||z||^2)^{n+1}}\cdot
\end{eqnarray*} 
\hfill $\Box$

\begin{lem} \label{FS_and_Euclidean_distance} 
For every $z\in \mathbb{C}^{n}$,
$\tan(d_{FS}(0,z))=d(0,z). $ 
\end{lem} 

\noindent \textit{Proof:} By applying a unitary transformation if necessary, we may assume that $z=(x,0,\ldots,0)$, where $x$ is a non-negative real number. Then the geodesic $\gamma$ between $0$ and $z$ with respect to the Fubini-Study metric is the real line segment joining them. By $\eqref{FSexplicit}$, for any $t\in[0,x]$, 

\[||\gamma^{\prime}(t)||^{2}_{FS}=\frac{|\gamma^{\prime}(t)|^2}{1+t^{2}}-\frac{t^{2}|\gamma^{\prime}(t)|^2}{(1+t^{2})^2}=\frac{|\gamma^{\prime}(t)|^2}{(1+t^{2})^2}\cdot\]
As a result, we get 
\[ d_{FS}(0,z)=\int_{\gamma}||\gamma^{\prime}||_{FS}=\int_{0}^{x} \frac{1}{1+t^{2}}dt = \arctan(x)=\arctan(d(0,z)). \] 
\hfill $\Box$ 

\begin{lem} \label{sinvolume} 
For every $0\leq \delta<\pi/2$, $\displaystyle{Vol(B_{FS}(0,\delta))=\frac{\pi^{n}}{n!}\sin^{2n}(\delta)}$. 
\end{lem} 

\noindent \textit{Proof:} Let $\delta=\arctan(R)$, then $B_{FS}(0,\delta)=B(0,R)$ by Lemma \ref{FS_and_Euclidean_distance}. Moreover, by Lemma \ref{volcompare},
\begin{eqnarray*}
Vol(B_{FS}(0,\delta))&=& \int_{B(0,R)} \frac{dVol}{(1+||z||^2)^{n+1}}  \\ 
&=&Vol(S^{2n-1}) \int_{0}^{R} \frac{r^{2n-1}}{(1+r^2)^{n+1}} dr \quad \mathrm{using}\,\mathrm{spherical}\,\mathrm{coordinates}\\
&=&\frac{Vol(S^{2n-1})}{2} \int_{1}^{1+R^2} \left(1-\frac{1}{\rho}\right)^{n-1} \frac{d\rho}{\rho^2} \qquad \mathrm{where}\, \rho=1+r^2 \\ 
&=& \frac{Vol(S^{2n-1})}{2} \int_{0}^{\frac{R^2}{R^2+1}} u^{n-1}du \qquad \mathrm{where} \, u=1-\frac{1}{\rho} \\
&=& \frac{Vol(S^{2n-1})}{2n} \left(\frac{R^2}{R^2+1}\right)^{n} \\ 
&=& \frac{\pi^{n}}{n!} \sin^{2n}(\delta),
\end{eqnarray*} 
since by the Hopf fibration $\displaystyle{Vol(S^{2n-1})=2\pi Vol_{FS}(\mathbb{CP}^{n-1})=\frac{2\pi^{n}}{(n-1)!}}\cdot$
\hfill $\Box$ 

By letting $\delta$ grow to $\pi/2$ in Lemma \ref{sinvolume}, one recovers the total volume of $\mathbb{CP}^{n}$, which can be deduced through cohomology as well, integrating the fundamental class against $\displaystyle{ \frac{ [\omega_{FS}]^{n}}{n!}}$. 

From this point on, we will concentrate on the case $n=2$ since this paper concerns plane curves. 
\begin{lem} \label{distancecompare} 
On the bidisc $\mathbb{B}=\{(z_{1},z_{2})\in \mathbb{C}^2| |z_{1}|\leq 1, |z_{2}|\leq 1\}\subset U_{0}$, the standard metric on $\mathbb{C}^2$ and its Fubini-Study one satisfy the inequalities 
$\displaystyle{ \frac{||h||}{3} \leq ||h||_{FS}\leq ||h||}$ 
for all $z\in \mathbb{B}$, and all $h\in \mathbb{C}^2=T|_{z}\mathbb{B}$. 
\end{lem} 

\noindent \textit{Proof:} By $\eqref{FSexplicit}$, for all $z\in \mathbb{B}$ and all $h\in \mathbb{C}^2=T|_{z}\mathbb{B}$, 
\[ ||h||_{FS}^2 =\frac{||h||^2}{1+||z||^2}-\frac{|\langle h,z\rangle |^2}{(1+||z||^2)^{2}}, \] 
so that $\displaystyle{||h||_{FS}^2\leq \frac{||h||^2}{1+||z||^2}\leq ||h||^2}$  and by Cauchy-Schwarz, 
\[ ||h||^{2}_{FS}\geq \frac{||h||^2}{1+||z||^2}-\frac{||h||^2 ||z||^2}{(1+||z||^2)^2} =\frac{||h||^2}{(1+||z||^2)^2}\cdot\] 
The condition $z\in \mathbb{B}$ implies that $||z||^2\leq 2$, hence the result. \hfill $\Box$  

\begin{lem} \label{thereexistpoints} 
Let $x,y>0$ and $0<\delta\leq 2\min(x,y)$. Then, there exist $N$ points $\mathbf{z}_{1},\ldots,\mathbf{z}_{N}$ on the torus $\mathbb{T}_{(x,y)}=\{(z_{1},z_{2})\in \mathbb{C}^2| |z_{1}|=x, |z_{2}|=y\}$ such that for all $i\neq j$, $||\mathbf{z}_{i}-\mathbf{z}_{j}||\geq \delta$, with 
\[ N=\left\lfloor\frac{4x}{\delta}\right\rfloor   \left\lfloor\frac{4y}{\delta}\right\rfloor\cdot\] 
\end{lem} 

\noindent \textit{Proof:} Suppose that $z_{1}, z_{1}^{\prime}$ are points in $\mathbb{C}$ such that $|z_{1}|=|z_{1}^{\prime}|=x$ and $\theta=|arg(z_{1}^{\prime}/z_{1})|$. Let $h$ denote the foot of the altitude of the isosceles triangle $(0z_{1}z_{1}^{\prime})$ emanating from the vertex $0$. The sine rule on the triangle $(0hz_{1})$ implies $|z_{1}-z_{1}^{\prime}|=2x\sin(\theta/2)$. Hence, $|z_{1}-z_{1}^{\prime}|$ gets minorized by $\delta$ provided that $\displaystyle{\sin(\theta/2)\geq \frac{\delta}{2x}}$, which imposes the condition $\displaystyle{\pi-2\arcsin\left(\frac{\delta}{2x}\right)\geq \theta \geq 2\arcsin\left(\frac{\delta}{2x}\right)}$. The maximal number of points with distance $\geq \delta$ that one can choose on the circle of radius $x$ is thus 
\[ \left\lfloor\frac{\pi}{\arcsin\left(\delta/(2x)\right)}\right\rfloor. \] 
The function $\displaystyle{\frac{\sin(t)}{t}}$ is decreasing on the interval $[0,\pi/2]$,  so that for all $t\in [0,\pi/2]$, $\displaystyle{\sin(t)\geq \frac{2t}{\pi}}$ and thus for all $u\in[0,1]$, $\displaystyle{\arcsin(u)\leq \frac{\pi u}{2} }$. Hence, the maximal number of points of distance $\geq \delta$ that one can place on the circle of radius $x$ is minorized by $\displaystyle{ \left\lfloor \frac{4x}{\delta} \right\rfloor}$. 

Let $(z_{1,i})_{1\leq i\leq N_{1}}$ be such a family of points and similarly $(z_{2,j})_{1\leq j\leq N_{2}}$ be a family of $\displaystyle{\left \lfloor \frac{4y}{\delta}\right\rfloor}$ points on the circle of radius $y$ in $\mathbb{C}$ centered at the origin such that $|z_{2,j}-z_{2,k}|\geq \delta$ for $j\neq k$. Then the points $\mathbf{z}_{ij}=(z_{1,i},z_{2,j})$ satisfy the requirement since $||\mathbf{z}_{ij}-\mathbf{z}_{kl}||^2=|z_{1,i}-z_{1,k}|^2+|z_{2,j}-z_{2,l}|^2$ gets minorized by $\delta^2$ as soon as one of the terms on the right hand side is, which happens whenever $\mathbf{z}_{ij}\neq \mathbf{z}_{kl}$. \hfill $\Box$ 

\subsection{Bergman kernel} 

Recall that $h_{d}$ induces an $L^2$ Fubini-Study hermitian product on $H^{0}(\mathbb{CP}^{n}, L^d)$ as follows: For all $P,Q\in H^{0}(\mathbb{CP}^{n}, L^d)$,
\[ \langle P,Q\rangle_{FS}=\frac{1}{Vol_{FS}(\mathbb{CP}^{n})}\int_{\mathbb{CP}^{n}} h_{d}(P(z),Q(z))_{FS} dVol_{FS}. \] 
The multinomials $X^{\mathbf{i}}$, where $\mathbf{i}$ ranges over all multiindices $\mathbf{i}=(i_{0},i_{1},\ldots,i_{n})\in \mathbb{N}^{n+1}$ with $d=i_{0}+\ldots+i_{n}$, form an orthogonal basis of 
$H^{0}(\mathbb{CP}^{n}, L^d)$ and we recall the norms of these basis elements.

\begin{lem} \label{FSnorm}
For every $\mathbf{i}=(i_{0},i_{1},\ldots,i_{n})\in \mathbb{N}^{n+1}$,
$\displaystyle{ ||X^{\mathbf{i}}||_{FS}^2 = {d+n \choose \mathbf{i}}^{-1}, }$
where $d=i_{0}+\ldots+i_{n}$ and 
\[ {d+n \choose \mathbf{i}}= \frac{(d+n)!}{i_{0}!i_{1}!\ldots i_{n}! n! }\cdot\]
Moreover, if $\mathbf{j}=(j_{0},\ldots,j_{n})\neq \mathbf{i}$, then $\langle X^{\mathbf{i}},X^{\mathbf{j}}\rangle=0$.  
\end{lem}

\noindent \textit{Proof:} By Lemma \ref{volcompare} and \eqref{sectionnorm},
\begin{eqnarray*} 
||X^{\mathbf{i}}||_{FS}^2 &=& \frac{1}{Vol_{FS}(\mathbb{CP}^{n})} \int_{\mathbb{C}^{n}} \frac{ |z_{1}|^{2i_{1}}\ldots |z_{n}|^{2i_{n}}}{(1+|z_{1}|^2+\ldots+|z_{n}|^2)^{d+n+1}} dVol \\ 
&=&  \frac{(2\pi)^{n}}{Vol_{FS}(\mathbb{CP}^{n})} \int_{0}^{\infty}\ldots\int_{0}^{\infty} \frac{r_{1}^{2i_{1}+1}\ldots r_{n}^{2i_{n}+1}}{(1+r_{1}^2+\ldots+r_{n}^2)^{d+n+1}} dr_{1}\ldots dr_{n}, 
\end{eqnarray*} 
where the last equality is obtained by passing to polar coordinates in each factor of $\mathbb{C}^{n}$. Consider the generating function: 
\[ F_{n}(t_{0},t_{1},\ldots,t_{n})=\int_{0}^{\infty}\ldots \int_{0}^{\infty} \frac{r_{1}\ldots r_{n}}{(t_{0}+t_{1}r_{1}^2+\ldots+t_{n}r_{n}^2)^{n+1}} dr_{1}\ldots dr_{n}. \] 
By Fubini's theorem, for every $i_{0},\ldots,i_{n}\geq 0$, 
\[  \int_{0}^{\infty}\ldots\int_{0}^{\infty} \frac{r_{1}^{2i_{1}+1}\ldots r_{n}^{2i_{n}+1}}{(1+r_{1}^2+\ldots+r_{n}^2)^{d+n+1}} dr_{1}\ldots dr_{n}  =\frac{(-1)^{d}}{(n+1)\ldots (n+d)}\frac{\partial^{d}}{\partial t_{0}^{i_{0}}\ldots \partial t_{n}^{i_{n}}} F_{n}|_{(t_{0},\ldots,t_{n})=(1,\ldots,1)}, \] 
with $d=i_{0}+\ldots+i_{n}$. By integrating with respect to the last variable $r_{n}$, we get the recursion formula
\begin{eqnarray*} 
F_{n}(t_{0},\ldots,t_{n})&=& \int_{0}^{\infty}\ldots \int_{0}^{\infty} \frac{r_{1}\ldots r_{n}}{(t_{0}+t_{1}r_{1}^2+\ldots+t_{n}r_{n}^2)^{n+1}} dr_{1}\ldots dr_{n} \\ 
&=& \frac{1}{2t_{n}n}\int_{0}^{\infty}\ldots \int_{0}^{\infty}  \frac{r_{1}\ldots r_{n-1}}{(t_{0}+t_{1}r_{1}^2+\ldots+t_{n-1}r_{n-1}^2)^{n}} dr_{1}\ldots dr_{n-1} \\
&=&\frac{1}{2t_{n}n} F_{n-1}(t_{0},\ldots,t_{n-1}). 
\end{eqnarray*} 
Observing that $\displaystyle{F_{1}(t_{0},t_{1})=\frac{1}{2t_{0}t_{1}}}$, we deduce 
\[ F_{n}(t_{0},\ldots,t_{n})=\frac{1}{2^{n}n! t_{0}\ldots t_{n}}\cdot \] 
Therefore, 
\begin{eqnarray*} 
||X^{\mathbf{i}}||_{FS}^2 &=&\frac{(2\pi)^{n}}{Vol_{FS}(\mathbb{CP}^{n})} 
\frac{(-1)^{d}}{(n+1)\ldots (n+d)}\frac{\partial^{d}}{\partial t_{0}^{i_{0}}\ldots \partial t_{n}^{i_{n}}} F_{n}|_{(t_{0},\ldots,t_{n})=(1,\ldots,1)} \\ 
&=& \frac{(2\pi)^{n} n! i_{0}!\ldots i_{n}!}{\pi^{n}2^{n}(d+n)! } \quad \mathrm{by} \, \mathrm{Lemma}\, \ref{volcompare} \, \mathrm{or}\, \ref{sinvolume} \\
&=& {d+n \choose \mathbf{i}}^{-1}. 
\end{eqnarray*} 
In order to prove the orthogonality statement, we note that if $\mathbf{j}=(j_{0},\ldots,j_{n})\neq \mathbf{i}=(i_{0},\ldots,i_{n})$, there exists $k$ such that $j_{k}\neq i_{k}$ and without loss of generality, we may assume that $i_{k}>j_{k}$. 
Computing
\[ \langle X^{\mathbf{i}},X^{\mathbf{j}}\rangle_{FS} = \frac{1}{Vol_{FS}(\mathbb{CP}^{n})} \int_{\mathbb{C}^{n}} \frac{ z_{1}^{i_{1}}\ldots z_{n}^{i_{n}}\overline{z_{1}}^{j_{1}}\ldots \overline{z_{n}}^{j_{n}}}{(1+||\mathbf{z}||^2)^{d+n+1}} dVol\]
as an iterated integral, we may take the innermost integral to be
\[ \mathcal{I}= \int_{\mathbb{C}} \frac{|z_{k}|^{2j_{k}}z_{k}^{i_{k}-j_{k}}}{(1+||\mathbf{z}||^2)^{\beta}} dVol \quad \mathrm{for}\,\mathrm{some}\, \mathrm{power} \, \beta\geq 0.\] 
But then, the change of variables $z_{k}\mapsto \xi z_{k}$, where $\xi$ denotes a primitive $2(i_{k}-j_{k})$th root of unity,  takes $\mathcal{I}$ to $-\mathcal{I}$, so that $\mathcal{I}=0$ and $\langle X^{\mathbf{i}},X^{\mathbf{j}}\rangle_{FS}=0$. 
\hfill $\Box$ 

Let us now denote the space of global $L^2$-sections of $\mathcal{O}_{\mathbb{CP}^{n}}(d)$ by $W^2(\mathbb{CP}^{n},\mathcal{O}_{\mathbb{CP}^{n}}(d))$ and let $\pi: W^2(\mathbb{CP}^{n},\mathcal{O}_{\mathbb{CP}^{n}}(d))\rightarrow H^{0}(\mathbb{CP}^{n},\mathcal{O}_{\mathbb{CP}^{n}}(d))$ be the orthogonal projection. Then, for every $\sigma\in W^2(\mathbb{CP}^{n},\mathcal{O}_{\mathbb{CP}^{n}}(d))$ and every $u\in \mathbb{CP}^{n}$, 
\[ \pi(\sigma)(u)=\frac{1}{Vol(\mathbb{CP}^{n})} \int_{\mathbb{CP}^{n}} \mathcal{B}(u,v)\otimes \sigma(v)dVol_{FS}(v), \] 
where $\mathcal{B}\in \Gamma(X\times X,L^{d} \boxtimes (L^{*})^{d})$ is the Bergman kernel and $L=\mathcal{O}_{\mathbb{CP}^{n}}(1)$, see \cite{MM}. 
For any orthonormal basis $\{\sigma_{i}\}$ of $H^{0}(\mathbb{CP}^{n}, L^{d})$,
\[ \mathcal{B}(u,v)=\sum_{i=1}^{N_{d}} \sigma_{i}(u)h_{FS}^{d}(\cdot, \sigma_{i}(v)), \quad \mathrm{where} \, N_{d}={d+n \choose n} \, \mathrm{denotes} \, \mathrm{the}\, \mathrm{dimension} \,\mathrm{of}\, H^{0}(\mathbb{CP}^{n},L^d). \]
In the open chart $U_{0}$, the Bergman kernel can be expressed as follows. One can write $\sigma_{i}(u)=f_{i}(u)(e^{*})^{d}$, where $f_{i}$ is a polynomial in the affine coordinates and $e^{*}=X_{0}$ the section of $L$ defined by $X_{0}$. Then, for all $u,v\in U_{0}$, 
\begin{equation} \label{Bexplicit}
 \mathcal{B}(u,v)=\left(\sum_{i=1}^{N_{d}} f_{i}(u)\overline{f_{i}}(v)||(e^{*})^{d}(u)||||(e^{*})^{d}(v)||\right)\left(\frac{(e^{*})^{d}(u)}{||(e^{*})^{d}(u)||}h_{FS}^{d}\left(\cdot, \frac{(e^{*})^{d}(v)}{||(e^{*})^{d}(v)||}\right)\right). 
\end{equation}
We rename the first of these factors, which is now a function, as $B(u,v)$. 
\begin{lem} \label{bergman}
For every $u,v\in \mathbb{C}^{n}\cong U_{0}$, 
\[ B(u,v)= \frac{(1+\langle u,v\rangle )^{d}N_{d}}{(1+||u||^2)^{d/2}(1+||v||^2)^{d/2}}\cdot \] 
\end{lem} 
\noindent \textit{Proof:}  For every $u,v\in \mathbb{C}^{n}\cong U_{0}$ choosing the orthonormal basis of multinomials given by Lemma \ref{FSnorm}, we deduce from \eqref{Bexplicit} and \eqref{sectionnorm} that
\begin{eqnarray*} 
& & (1+||u||^2)^{d/2}(1+||v||^2)^{d/2} B(u,v) =  \sum_{0\leq i_{1}+\ldots+i_{n}\leq d} {d+n \choose \mathbf{i}} (u_{1}\overline{v_{1}})^{i_{1}}\ldots (u_{n}\overline{v_{n}})^{i_{n}} \\ 
&=& {d+n \choose n} \sum_{0\leq i_{1}+\ldots + i_{n-1}\leq d} \frac{d!(u_{1}\overline{v_{1}})^{i_{1}}\ldots (u_{n-1}\overline{v_{n-1}})^{i_{n-1}}}{i_{1}!i_{2}!\ldots i_{n-1}!(d-i_{1}-\ldots-i_{n-1})!} \sum_{i_{n}=0}^{d-i_{1}-\ldots-i_{n-1}} {d-i_{1}-\ldots-i_{n-1} \choose i_{n}} (u_{n}\overline{v_{n}})^{i_{n}} \\ 
&=&N_{d} \sum_{0\leq i_{1}+\ldots + i_{n-1}\leq d} \frac{d!(u_{1}\overline{v_{1}})^{i_{1}}\ldots (u_{n-1}\overline{v_{n-1}})^{i_{n-1}}}{i_{1}!i_{2}!\ldots i_{n-1}!(d-i_{1}-\ldots-i_{n-1})!} (1+u_{n}\overline{v_{n}})^{d-i_{1}-\ldots-i_{n-1}} \\  
&=& N_{d}(1+\langle u,v\rangle )^{d} \qquad \mathrm{by} \, \mathrm{recursion.} \\
\end{eqnarray*} 
\hfill $\Box$ 

\begin{rem} 
(1) In particular, $B$ is constant and equals $N_{d}$ on the diagonal of $\mathbb{CP}^{n}\times \mathbb{CP}^{n}$. 

(2) In contrast, $B$ decreases exponentially outside the diagonal, by the Cauchy-Schwarz inequality, since
\[ B(u,v)=\frac{\langle (1,u),(1,v)\rangle ^{d} N_{d}}{||(1,u)||^{d}||(1,v)||^{d}}\cdot\]  
\end{rem} 

\iffalse
For every $u,v\in \mathbb{C}^n\cong U_{0}$ and every $i\in \{1,\ldots,n\}$, $d\geq 2$, 
\begin{eqnarray*} 
\frac{\partial^{2}}{\partial u_{i}\partial{\overline{v_{i}}}|_{(u,v)}} (1+\langle u,v\rangle )^{d} &=& d(1+\langle u,v\rangle )^{d-1}+d(d-1)u_{i}\overline{v_{i}}(1+\langle u,v\rangle )^{d-2} \\ 
&=& d(1+\langle u,v\rangle )^{d-2}(1+\langle u,v\rangle +(d-1)u_{i}\overline{v_{i}}) \\ 
&=& \frac{\partial}{\partial u_{i}|_{(u,v)}}(du_{i}(1+\langle u,v\rangle )^{d-1}) 
\end{eqnarray*} 
\fi

\begin{lem} \label{bergmanestimate} 
For all $u,v\in \mathbb{C}^{n}\cong U_{0}$, 
$ |B(u,v)|\leq N_{d}\cos(d_{FS}(u,v))^{d}.$ 
\end{lem} 

\noindent\textit{Proof:} We may suppose without loss of generality that $v=0$, since the isometry group of $\mathbb{CP}^{n}$ acts transitively. Then, by Lemmas \ref{bergman} and \ref{FS_and_Euclidean_distance}, 
\[ |B(u,0)|=\frac{N_{d}}{(\sqrt{1+||u||^2})^{d}}=N_{d}\cos(d_{FS}(u,0))^{d}. \]
\hfill $\Box$

Let $N\geq 1$ and $\mathbf{z}_{1},\ldots,\mathbf{z}_{N}$ be $N$ distinct points of $\mathbb{C}^{n}\cong U_{0}$. For every $P\in H^{0}(\mathbb{CP}^n,L^{d})$ and every $i\in \{1,\ldots,N\}$, we set 
\[  P(\mathbf{z}_{i})=Z_{i}(P) \frac{(e^{*})^{d}(\mathbf{z}_{i})}{||(e^{*})^{d}(\mathbf{z}_{i})||} \in L^{d}_{|_{\mathbf{z}_{i}}}. \] 
It defines a complex Gaussian vector $\mathbf{Z}=(Z_{1},\ldots,Z_{N})$, since the Hermitian space $H^{0}(\mathbf{CP}^{n}, L^d)$ inherits its associated Gaussian measure $\mu$ having density $\displaystyle{\frac{1}{\pi^{N_{d}}}e^{-||P||^2}}$
with respect to the Lebesgue measure, see Definition \ref{Gaussiandef}. The Bergman kernel carries the following probabilistic meaning. 

\begin{lem} \label{covariancebergman} 
Let $N\geq 1$ and $\mathbf{z}_{1},\ldots,\mathbf{z}_{N}$ be $N$ distinct points of $\mathbb{C}^n\cong U_{0}$. Then, for every $1\leq i,j\leq N$, 
\[ Cov(Z_{i},Z_{j})=B(\mathbf{z}_{i},\mathbf{z}_{j}). \] 
\end{lem} 

\noindent \textit{Proof:} For every $1\leq i,j\leq N$, 
\begin{eqnarray*} 
Cov(Z_{i},Z_{j}) &=& \mathbb{E}(Z_{i}\overline{Z_{j}}) \\ 
&=& \frac{1}{\pi^{N_{d}}} \int_{H^{0}(\mathbb{CP}^n,L^{d})} Z_{i}(P)\overline{Z_{j}(P)} d\mu(P) \\ 
&=& \sum_{1\leq k,l\leq N_{d}} Z_{i}(\sigma_{k})\overline{Z_{j}(\sigma_{l})} \left( \frac{1}{\pi^{N_{d}}} \int_{\mathbb{C}^{N_{d}}} a_{k}\overline{a_{l}}d\mu(a)\right) \\
&=& \sum_{k=1}^{N_{d}} Z_{i}(\sigma_{k})\overline{Z_{j}(\sigma_{k})} 
\end{eqnarray*} 
where $(\sigma_{1},\ldots,\sigma_{N_{d}})$ denotes an orthonormal basis of $H^{0}(\mathbb{CP}^n,L^{d})$ and $P=\sum_{k=1}^{N_{d}} a_{k}\sigma_{k}$ with $a=(a_{1},\ldots,a_{N_{d}})\in \mathbb{C}^{N_{d}}$. Writing $\sigma_{k}=f_{k}(e^{*})^{d}$, we get $Z_{i}(\sigma_{k})=f_{k}(\mathbf{z}_{i})||(e^{*})^{d}(\mathbf{z}_{i})||$ and thus by \eqref{Bexplicit}, $Cov(Z_{i},Z_{j})=B(\mathbf{z}_{i},\mathbf{z}_{j})$. \hfill $\Box$

\section{Amoeba measures and upper bounds} \label{upper}

\subsection{Amoeba measures} \label{amoebameasures} 
For every $(x,y)\in (\mathbb{R}_{+}^{\times})^2$, let 
\[ \mathbb{T}_{(x,y)}=\{[1:z_{1}:z_{2}]\in \mathbb{CP}^2 | |z_{1}|=x, |z_{2}|=y \}, \]
and let $A(x,y)$ be its area computed with respect to the Fubini-Study metric of $\mathbb{CP}^{2}$ restricted to $\mathbb{T}_{(x,y)}$, so that
\begin{equation} \label{arearestrict} 
A(x,y)=\int_{\mathbb{T}_{(x,y)}} \frac{d Vol_{\mathbb{T}}}{(1+||z||^2)^{3/2}} = \frac{4\pi^2 x y}{(1+x^2+y^2)^{3/2}}\cdot 
\end{equation} 
The compact torus $\mathbb{T}_{(x,y)}$ is indeed Lagrangian in $\mathbb{CP}^2$, so that the density $(1+||z||^2)^{-3/2}$ of its Fubini-Study area form with respect to the standard one is the square-root of the density given by Lemma \ref{volcompare}.

\begin{rem} \label{volrem}
(1) Integration over $(\mathbb{R}_{+}^{\times})^2$ recovers the volume of $\mathbb{CP}^2$ from the co-area formula. Indeed, $\mathbb{RP}^2$ is Lagrangian as well in $\mathbb{CP}^2$, so that its Fubini-Study area form at $z\in \mathbb{RP}^2$ equals $\displaystyle{ \frac{dx dy}{\left(1+||z||^2\right)^{3/2}}}$ and
\begin{eqnarray*} 
Vol_{FS}(\mathbb{CP}^2)&=& \iint_{(\mathbb{R}_{+}^{\times})^2}\frac{A(x,y)}{(1+||z||^2)^{3/2}}dxdy=4\pi^2 \int_{0}^{+\infty}\int_{0}^{+\infty} \frac{xy}{(1+x^2+y^2)^3}dxdy \\ 
&=& \pi^2 \int_{0}^{+\infty} x \left( \frac{-1}{(1+x^2+y^2)^2} |_{0}^{+\infty} \right) dx = \pi^2 \left( \frac{-1}{2(1+x^2)}|_{0}^{+\infty} \right) =\frac{\pi^2}{2}\cdot
\end{eqnarray*} 

\noindent (2) The area functional $A(x,y)$ reaches its maximum at $(1,1)$, since 
$\displaystyle{\frac{\partial A}{\partial x}=\frac{4\pi^{2}y(1-2x^2+y^2)}{(1+x^2+y^2)^{5/2}}}$, $\displaystyle{\frac{\partial A}{\partial y}=\frac{4\pi^{2}x(1-2y^2+x^2)}{(1+x^2+y^2)^{5/2}}}$ and both vanish only at this point in $(\mathbb{R}_{+}^{\times})^2$. Moreover, the area of the Clifford torus $\mathbb{T}_{(1,1)}$ equals  
$\displaystyle{A(1,1)=\frac{4\pi^2}{3\sqrt{3}}}\cdot$ 
\end{rem}

Let $C_{b}^{0}(\mathbb{R}^2,\mathbb{R})$ be the space of bounded continuous functions on $\mathbb{R}^2$ and for every $P\in H^{0}(\mathbb{CP}^2,L^d)$, let
\begin{equation} \label{measures}
\lambda_{P}:\varphi\in C^{0}_{b}(\mathbb{R}^2,\mathbb{R})\mapsto \int_{\mathcal{A}(P)}\varphi(t)|dt| \quad \mathrm{and} \quad \nu_{P}: \varphi\in C^{0}_{b}(\mathbb{R}^2,\mathbb{R})\mapsto \int_{V_{P}}Log^{*}(\varphi)|Log^{*}dt|, 
\end{equation} 
be the amoeba measures of $V_{P}$, so that $\langle\lambda_{P},1\rangle=Vol(\mathcal{A}(P))$ and $\langle\nu_{P},1\rangle=Mvol(\mathcal{A}(P))$ compute the volume and multivolume of the amoeba of $P$ respectively, see \cite{BK}. We denote by
 \[\mathbb{E}_{d}(\nu)=\int_{H^{0}(\mathbb{CP}^2,L^d)} \nu_{P}d\mu(P) \qquad \mathrm{and} \qquad \mathbb{E}_{d}(\lambda)=\int_{H^{0}(\mathbb{CP}^2,L^d)} \lambda_{P}d\mu(P)\]  
the corresponding average measures and finally set 
\[ a:t=(t_{1},t_{2})\in \mathbb{R}^2\mapsto Area_{FS}(Log^{-1}(t))=A(e^{t_{1}},e^{t_{2}}). \] 
and $p:t\in \mathbb{R}^2\mapsto \mathbb{P}(t\in \mathcal{A}(P))$, where $\mathbb{P}(t\in \mathcal{A}(P))$ denotes the probability that $t$ belongs to the amoeba of $P$, that is the probability that $V_{P}\cap Log^{-1}(t)\neq \emptyset$.

\begin{pro} \label{expectedmeasure}
\begin{enumerate} 
\item For all $d>0,$ $\displaystyle{\mathbb{E}_{d}(\nu)=\frac{d}{2\pi} a(t)|dt|=2\pi d Log_{*}|dVol_{FS}|_{|(\mathbb{R}^{*}_{+})^2}}$.  \\ 
\item For all $d>0,$ $\displaystyle{\mathbb{E}_{d}(\lambda)=p(t) |dt|}$. 
\end{enumerate} 
\end{pro} 
\noindent\textit{Proof:} (1) Say $\varphi\in C_{b}^{0}(\mathbb{R}^2,\mathbb{R})$. Then, 
\begin{eqnarray*} 
\langle\mathbb{E}_{d}(\nu),\varphi\rangle&=& \int_{H^{0}(\mathbb{CP}^2,L^d)} \langle\nu_{P},\varphi\rangle d\mu(P) =\int_{H_{0}(\mathbb{CP}^2,L^d)}\int_{V_{P}}Log^{*}\varphi|Log^{*}dt|d\mu(P) \\
&=& \int_{H^{0}(\mathbb{CP}^2,L^d)}\left(\int_{\mathbb{R}^2}\varphi(t) \#(V_{P}\cap Log^{-1}(t))|dt|\right) d\mu(P) \quad \mathrm{by}\,\mathrm{Sard's}\,\mathrm{theorem} \\ 
&=& \int_{\mathbb{R}^2}\varphi(t)\mathbb{E}_{d}\left(\#(V_{P}\cap Log^{-1}(t))\right)|dt| \\ 
&=& \frac{d}{2\pi} \int_{\mathbb{R}^2}a(t)\varphi(t) |dt| \quad \mathrm{by}\, \mathrm{Crofton's}\,\mathrm{formula}, 
\end{eqnarray*} 
so that $\displaystyle{\mathbb{E}_{d}(\nu)=\frac{d}{2\pi}a(t)|dt|}$, see Theorem \ref{crofton}. 

In order to obtain the second equality, we observe that by Remark \ref{volrem} and \eqref{arearestrict}
\[ Log_{*}|dVol_{FS}|_{|(\mathbb{R}^{*}_{+})^2}=\frac{e^{t_{1}}e^{t_{2}}|dt_{1}dt_{2}|}{(1+e^{2t_{1}}+e^{2t_{2}})^{3/2}}=\frac{a(t)|dt|}{4\pi^2}=\frac{1}{2\pi d}\mathbb{E}_{d}(\nu). \] 

(2)  For every $\varphi\in C_{b}^{0}(\mathbb{R}^2,\mathbb{R})$,
\begin{eqnarray*} 
\langle\mathbb{E}_{d}(\lambda),\varphi\rangle&=& \int_{H^{0}(\mathbb{CP}^2,L^d)}\langle\lambda_{P},\varphi\rangle d\mu(P) = \int_{H^{0}(\mathbb{CP}^2,L^d)}\int_{\mathbb{R}^2} \mathbf{1}_{\mathcal{A}(P)}(t)\varphi(t)|dt| d\mu(P)\\ 
&=& \int_{\mathbb{R}^{2}}\varphi(t) \left(\int_{H^{0}(\mathbb{CP}^2,L^d)} \mathbf{1}_{\mathcal{A}(P)}(t)d\mu(P)\right) |dt| =\int_{\mathbb{R}^2}\varphi(t)p(t)|dt|. 
\end{eqnarray*} 
\hfill $\Box$ 

\begin{cor} \label{expbound}
For every $d>0$, 
\begin{enumerate} 
\item  For all $ t\in \mathbb{R}^2,$ $\displaystyle{0\leq p(t)\leq \min\left(1,da(t)/(4\pi)\right)}$, \\ 
\item $\mathbb{E}_{d}(Mvol(\mathcal{A}))=\pi^2 d,$ \\ 
\item $\displaystyle{\mathbb{E}_{d}(Vol(\mathcal{A}))\leq \int_{\mathbb{R}^2} \min\left(1, da(t)/(4\pi)\right)dt.}$
\end{enumerate}
\end{cor} 

\noindent\textit{Proof:} (1) For all $\varphi\in C_{b}^{0}(\mathbb{R}^2,\mathbb{R})$, $\displaystyle{\langle \lambda_{P},\varphi\rangle \leq \langle \nu_{P},\varphi\rangle/2}$, since for all regular values $t$ of $Log$ restricted to $V_{P}$, $Log^{-1}(t)$ is cobordant to $\emptyset$ and thus consists of an even number of points. Hence, $\displaystyle{\langle\mathbb{E}_{d}(\lambda),\varphi\rangle \leq \langle \mathbb{E}_{d}(\nu),\varphi\rangle/2}$. By choosing $(\varphi_{n})_{n\in \mathbb{N}}$ to approximate the Dirac delta function at $t$, we deduce from Proposition \ref{expectedmeasure} that $\displaystyle{p(t)\leq da(t)/(4\pi)}$ for any $t\in \mathbb{R}^2$, while $p(t)\leq 1$ being a probability. 
 
(2) 
\begin{eqnarray*}
\mathbb{E}_{d}(Mvol(\mathcal{A}))&=& \langle \mathbb{E}_{d}(\nu),1\rangle \\ 
&=&2\pi d\int_{\mathbb{R}^2} Log_{*}|dVol_{FS}|_{(\mathbb{R}_{+}^{\times})^2} \qquad \mathrm{by}\,\,\mathrm{Proposition} \,\,\mathrm{\ref{expectedmeasure}} \\ 
&=& 2\pi d Vol_{FS}(\mathbb{R}_{+}^{\times})^2  \\ 
&=& \frac{\pi d}{2}Vol_{FS}(\mathbb{RP}^2)=\pi^2 d, 
\end{eqnarray*}  
since $\mathbb{RP}^2$ is a 2-to-1 quotient of the unit sphere, which has volume $4\pi$. 

(3) Follows from part (1) of Corollary \ref{expbound} and part (2) of Proposition \ref{expectedmeasure}, since $\mathbb{E}_{d}(Vol(\mathcal{A}))=\langle \mathbb{E}_{d}(\lambda),1\rangle$.
\hfill  $\Box$  

\begin{rem} \label{BKremark}
Part (2) of Corollary \ref{expbound} has already been proved in \cite{BK}  using the argument map instead of the $Log$ map. 
\end{rem} 

\subsection{Level set of the area functional and asymptotic upper estimates} \label{levelset}
In the light of Corollary \ref{expbound}, we now study the level set $\{t\in \mathbb{R}^2|a(t)=4\pi/d\}$ of the area functional, or rather, for convenience, the corresponding level set $\{(x,y)\in (\mathbb{R}_{+}^{\times})^2|A(x,y)=4\pi/d\}$. The restriction of $A$ to the diagonal reads, up to a constant, as the function $g:\rho\in \mathbb{R}_{+}^{\times}\mapsto g(\rho)=f(\rho^2)$ with $\displaystyle{f:t\in \mathbb{R}_{+}\mapsto t(1+t)^{-3/2}}$. Since $\displaystyle{f^{\prime}(t)=(1-t/2)(1+t)^{-5/2}}$, the latter reaches a unique extremum point at $t=2$ which is a global maximum, and the maximum value is $\displaystyle{f(2)=2/(3\sqrt{3})}$. The function  $f$ is strictly increasing on the interval $[0,2]$ and strictly decreasing on $[2,+\infty)$ with $f(0)=\lim_{t\rightarrow +\infty}f(t)=0$. Therefore, for any positive $k$ such that $\displaystyle{k<2/(3\sqrt{3})}$, the equation $f(t)=k$ has two solutions $0<t_{1}<2<t_{2}$. The same remark holds for the function $g(\rho)=f(\rho^{2})$. Since $\displaystyle{2/(\pi d)<2/(3\sqrt{3})}$ for $d\geq 2$, we conclude that the equation 
\[ \frac{\rho^{2}}{(1+\rho^2)^{3/2}}=\frac{2}{\pi d} \] 
has precisely two positive solutions $0< \rho_{0}<\sqrt{2}<\rho_{1}$ for each $d\geq 2$ .

\begin{lem} \label{asymptotic}
$\displaystyle{\rho_{0}^2=\frac{2}{\pi d}+O(d^{-2})}$ and $\displaystyle{\frac{1}{\rho_{1}^2}=\frac{4}{\pi^2 d^2}+O(d^{-4})}$. 
\end{lem} 

\noindent \textit{Proof:} By definition, $\displaystyle{\rho_{0}^4=(4/(\pi^2 d^2))(1+3\rho_{0}^2+3\rho_{0}^4+\rho_{0}^6)}$ so that $\displaystyle{\rho_{0}^2=2/(\pi d)+O(d^{-2})}$. Likewise,  $\displaystyle{\rho_{1}^4=(4/(\pi^2 d^2))(1+3\rho_{1}^2+3\rho_{1}^4+\rho_{1}^6)}$, so that  
\[ \frac{1}{\rho_{1}^2}=\frac{4}{\pi^2 d^2}\left(1+\frac{3}{\rho_{1}^2}+\frac{3}{\rho_{1}^4}+\frac{1}{\rho_{1}^6}\right)\cdot \] 
Consequently, $\displaystyle{\rho_{1}^{-2}=4/(\pi^{2}d^2)+O\left(d^{-4}\right)}$.
 \hfill $\Box$ 

\begin{rem} 
In fact, $\displaystyle{\rho_{0}^2=\frac{2}{\pi d}+\frac{6}{\pi^2 d^2}+\frac{21}{\pi^3 d^3}+o(d^{-3})}$ and $\displaystyle{\frac{1}{\rho_{1}^2}=\frac{4}{\pi^2 d^2}+\frac{48}{\pi^4 d^4}+\frac{768}{\pi^6 d^6}+o\left(d^{-6}\right)}$.
\end{rem}

\begin{cor} \label{lnest}
\begin{eqnarray*} 
\ln(\rho_{0})&=&-\ln(d)/2-\ln\left(\pi/2\right)/2+O\left(\frac{1}{d}\right) \quad \mathrm{and} \\ 
\ln(\rho_{1})&=& \ln(d)+\ln\left(\pi/2\right)+O\left(\frac{1}{d^2}\right). \\ 
\end{eqnarray*} 
\end{cor} 

\noindent\textit{Proof:} By Lemma \ref{asymptotic}, $\displaystyle{\ln(\rho_{0}^2)=\ln\left(2/(\pi d) \right)+\ln\left(1+O\left(\frac{1}{d}\right)\right)}$ and thus
\[\ln(\rho_{0})=-\frac{\ln(d)}{2}-\frac{\ln\left(\pi/2\right)}{2}+O\left(\frac{1}{d}\right).\] 
Similarly, by Lemma \ref{asymptotic}, 
\begin{eqnarray*} 
\ln\left(\frac{1}{\rho_{1}^2}\right)&=& \ln\left(\frac{4}{\pi^2 d^2}\right)+\ln\left(1+O\left(\frac{1}{d^2}\right)\right), \quad \mathrm{so} \, \mathrm{that} \\
\ln(\rho_{1})&=&\ln(d)+\ln\left(\pi/2\right)+O\left(\frac{1}{d^2}\right). 
\end{eqnarray*} 
\hfill $\Box$ 

\begin{rem} \label{remark38}
In particular, $\ln(\rho_{0})=-\ln(\rho_{1})/2+O(1/d).$
In fact, $\ln(\rho_{0})=-\ln(d)/2-\ln\left(\pi/2\right)/2+3/(2\pi d)+o\left(d^{-1}\right)$ and $\ln(\rho_{1})= \ln(d)+\ln\left(\pi/2\right)-6/(\pi^2 d^2)+o\left(d^{-2}\right)$.
\end{rem} 

In polar coordinates $(\rho,\theta)\in \mathbb{R}_{+}\times [0,\pi/2]$, the area functional reads
\[ A:(x,y)\in \mathbb{R}_{+}^2\mapsto \frac{4\pi^2 xy}{(1+x^2+y^2)^{3/2}}=\frac{2\pi^{2}\sin(2\theta)\rho^2}{(1+\rho^2)^{3/2}}\cdot\] 
This functional is symmetric with respect to the diagonal via the symmetry $\displaystyle{\theta\mapsto \pi/2-\theta}$, and vanishes on the axes. For a fixed value of $\rho$, the function $A$ increases on $\displaystyle{\theta\in \left[0,\pi/4\right]}$ and decreases for $\displaystyle{\theta\in \left[\pi/4,\pi/2\right]}$, hence it is maximal on the diagonal. On the diagonal, as a function of $\rho$, the maximum is reached for $\rho=\sqrt{2}$, as discussed above. 

By Corollary \ref{expbound}, 
\[ \mathbb{E}_{d}(Vol(\mathcal{A}))\leq \int_{\mathbb{R}^2}\min\left(1,da(t)/(4\pi)\right)|dt|=\int_{(\mathbb{R}_{+}^{\times})^2} \min\left(1,dA(x,y)/(4\pi)\right)\frac{dxdy}{xy}\cdot \] 
Let us now prove the asymptotic expansion \eqref{upasym} of the right hand side, given in the introduction. 
\begin{thm} \label{upperboundasymptotic} 
As $d$ tends to $+\infty$, 
\[ \int_{\mathbb{R}^{2}}\min\left(1, da(t)/(4\pi)\right)dt = \frac{3}{2}\ln(d)^{2}+3\left(1+\ln\left(\pi\right)\right)\ln(d)+O(1). \] 
\end{thm} 

\noindent\textit{Proof:} 
Let us decompose the locus of integration into three parts: (i) The set where $\rho\leq \rho_{0}$ or $\rho\geq \rho_{1}$, (ii) the set where $\rho_{0}\leq \rho \leq \rho_{1}$ and either $0\leq \theta\leq \theta_{\rho}$ or $\displaystyle{\pi/2-\theta_{\rho}\leq \theta\leq \pi/2}$, where 
\begin{equation} \label{sintheta}
 \sin(2\theta_{\rho})=\frac{2(1+\rho^2)^{3/2}}{d\pi\rho^2},
\end{equation}
so that $(\rho, \theta_{\rho})$ parameterizes the locus below the diagonal $1=dA/(4\pi)$,  
and (iii) the set where $\rho_{0}\leq \rho\leq \rho_{1}$ and $\displaystyle{\theta_{\rho}\leq \theta \leq \pi/2-\theta_{\rho}}$. 

(i) 
\begin{eqnarray*}
\iint_{\{\rho\leq \rho_{0}\}\cup \{\rho\geq \rho_{1}\}} \frac{d}{4\pi} A(x,y)\frac{dxdy}{xy} 
&=& \pi d \int_{0}^{\pi/2} \int_{[0,\rho_{0}] \cup [\rho_{1},\infty)} \frac{\rho}{(1+\rho^2)^{3/2}} d\rho d\theta \\
&=& \frac{\pi^2 d}{2} \left(1-(1+\rho_{0}^2)^{-1/2}+(1+\rho_{1}^2)^{-1/2}\right) \\ 
&=& \frac{\pi^2 d}{2} \left( 1-\left(1-\frac{1}{2}\rho_{0}^2+o(\rho_{0}^2)\right)+\frac{1}{\rho_{1}}\left(1-\frac{1}{2\rho_{1}^2}+o\left(\rho_{1}^{-2}\right)\right)\right)\\ 
&=& O(1) 
\end{eqnarray*} 
by Lemma \ref{asymptotic}. 

(ii) Let $\omega(d)$ be any positive valued increasing function of $d$ going to infinity as $d\rightarrow +\infty$ such that if $\displaystyle{\widetilde{\rho}_{1}=\frac{\rho_{1}}{\omega(d)}}$ and $\widetilde{\rho}_{0}=\omega(d)\rho_{0}$, then $\rho_{0}<\widetilde{\rho}_{0}<\sqrt{2}<\widetilde{\rho}_{1}<\rho_{1}$ for all $d$. We split the integral in the following way: 
\begin{eqnarray*} 
\int_{\rho_{0}}^{\rho_{1}} \int_{[0,\theta_{\rho}]\cup [\frac{\pi}{2}-\theta_{\rho},\frac{\pi}{2}]} \frac{d}{4\pi}A(x,y) \frac{dxdy}{xy}&=& 2\pi d \int_{\rho_{0}}^{\rho_{1}} \int_{0}^{\theta_{\rho}} \frac{\rho}{(1+\rho^2)^{3/2}} d\theta d\rho \\
&=& 4\int_{\rho_{0}}^{\rho_{1}} \frac{\theta_{\rho}}{\rho \sin(2\theta_{\rho})} d\rho  \quad \mathrm{by} \, \eqref{sintheta} \\ 
&=& 4\int_{\rho_{0}}^{\widetilde{\rho}_{0}} \frac{\theta_{\rho}}{\rho \sin(2\theta_{\rho})}d\rho + 4\int_{\widetilde{\rho}_{0}}^{\widetilde{\rho}_{1}} \frac{\theta_{\rho}}{\rho \sin(2\theta_{\rho})}d\rho +4\int_{\widetilde{\rho}_{1}}^{\rho_{1}} \frac{\theta_{\rho}}{\rho \sin(2\theta_{\rho})}d\rho. 
\end{eqnarray*} 
Then,
\[  0\leq 4\int_{\rho_{0}}^{\widetilde{\rho}_{0}} \frac{\theta_{\rho}}{\rho \sin(2\theta_{\rho})}d\rho  \leq \pi\ln\left(\widetilde{\rho}_{0}/\rho_{0}\right)=\pi\ln(\omega(d)), \] 
since $\displaystyle{x\mapsto \frac{x}{\sin x}}$ is increasing and bounded from above by $\pi/2$ on $[0,\pi/2]$ and similarly, 
\[ 0\leq 4\int_{\widetilde{\rho}_{1}}^{\rho_{1}} \frac{\theta_{\rho}}{\rho \sin(2\theta_{\rho})}d\rho  \leq \pi\ln\left(\rho_{1}/\widetilde{\rho}_{1}\right)=\pi\ln(\omega(d)). \]

In order to estimate the second integral, we observe that  $\displaystyle{\sin^2(2\theta_{\rho})=\left(2/(\pi d)\right)^2\left(\rho^{-4}+3\rho^{-2}+3+\rho^2\right)}$ by \eqref{sintheta} and that $\theta_{\rho}=o(1)$ for $\rho\in [\widetilde{\rho}_{0},\widetilde{\rho}_{1}]$, by the choice of $\omega(d)$. We deduce
\begin{eqnarray*} 
 4\int_{\widetilde{\rho}_{0}}^{\widetilde{\rho}_{1}} \frac{\theta_{\rho}}{\rho \sin(2\theta_{\rho})}d\rho&=& 4\int_{\widetilde{\rho}_{0}}^{\widetilde{\rho}_{1}} \frac{\theta_{\rho}}{\rho \left(2\theta_{\rho}-4\theta_{\rho}^3/3+o(\theta_{\rho}^3)\right)}d\rho \\ 
&=&  2\int_{\widetilde{\rho}_{0}}^{\widetilde{\rho}_{1}}  (1+2\theta_{\rho}^2/3+o(\theta_{\rho}^2)) d\rho/  \rho \\ 
&=& 2\ln\left(\widetilde{\rho}_{1}/\widetilde{\rho}_{0}\right)+O\left(\int_{\widetilde{\rho}_{0}}^{\widetilde{\rho}_{1}} \theta_{\rho}^2 d\rho/\rho \right), \,\,\, \mathrm{while} \\
\int_{\widetilde{\rho}_{0}}^{\widetilde{\rho}_{1}} \theta_{\rho}^2 d\rho/\rho &=& O\left(\frac{1}{\pi^2 d^2} \int_{\widetilde{\rho}_{0}}^{\widetilde{\rho}_{1}} \left(\rho^{-5}+3\rho^{-3}+3\rho^{-1}+\rho\right)d\rho\right) \\ 
&=& O\left(\frac{1}{\pi^2 d^2}  \left( -\rho^{-4}/4-3\rho^{-2}/2+3\ln(\rho)+\rho^2/2\right)|_{\widetilde{\rho}_{0}}^{\widetilde{\rho}_{1}}\right) \\ 
&=& O(1) 
\end{eqnarray*} 
by Lemma \ref{asymptotic}. Combining the estimates above, we obtain 
\[
\int_{\rho_{0}}^{\rho_{1}} \int_{[0,\theta_{\rho}]\cup [\frac{\pi}{2}-\theta_{\rho},\frac{\pi}{2}]} \frac{d}{4\pi}A(x,y) \frac{dxdy}{xy}=2\ln\left(\widetilde{\rho}_{1}/\widetilde{\rho}_{0}\right)+O(\ln(\omega(d))) 
= 3\ln(d)+O(1)
\]
by Corollary \ref{lnest}, since this estimate holds for any function $\omega(d)$ going to $+\infty$ satisfying the additional conditions above.

(iii)   
\begin{eqnarray*} 
 \int_{\rho_{0}}^{\rho_{1}} \int_{\theta_{\rho}}^{\pi/2-\theta_{\rho}} \frac{dxdy}{xy} &=& 2\int_{\rho_{0}}^{\rho_{1}}\int_{\theta_{\rho}}^{\pi/4} \frac{d\rho d\theta}{\rho\sin\theta \cos \theta} =2\int_{\rho_{0}}^{\rho_{1}} \frac{1}{\rho} \ln(\tan \theta) |_{\theta_{\rho}}^{\pi/4}d\rho \\ 
&=& 2\int_{\rho_{0}}^{\rho_{1}} \frac{1}{\rho} \ln\left( \frac{1}{\tan \theta_{\rho}}\right) d\rho = 2\int_{\rho_{0}}^{\rho_{1}} \ln\left( \frac{2\cos^2 \theta_{\rho}}{2\cos \theta_{\rho}\sin\theta_{\rho}}\right) \frac{d\rho}{\rho} \\
&=& 2\int_{\rho_{0}}^{\rho_{1}}\left(\ln\left(\frac{d\pi \rho^2}{2(1+\rho^2)^{3/2}}\right) +\ln\left(1+\cos(2\theta_{\rho})\right)\right)\frac{d\rho}{\rho} \quad \mathrm{by} \, \eqref{sintheta} \\
&=& 2\ln\left(d \pi/2\right)\ln\left(\rho_{1}/\rho_{0}\right)+2\left(\ln(\rho_{1})^2-\ln(\rho_{0})^2\right)\\
&&-3\int_{\rho_{0}}^{\rho_{1}}\ln(1+\rho^2)d\rho/\rho+2\int_{\rho_{0}}^{\rho_{1}}\ln(1+\cos(2\theta_{\rho}))d\rho/\rho. 
\end{eqnarray*} 
We now estimate the two integrals in the last expression:  
\begin{eqnarray*} 
\int_{\rho_{0}}^{\rho_{1}} \ln(1+\rho^2)d\rho/\rho &=& \int_{\rho_{0}}^{1} \ln(1+\rho^2)d\rho/\rho+\int_{1}^{\rho_{1}}\ln(1+\rho^2)d\rho/\rho \\
&=&  \int_{\rho_{0}}^{1} \ln(1+\rho^2)d\rho/\rho+\int_{1}^{\rho_{1}}\ln(\rho^2)d\rho/\rho+\int_{1}^{\rho_{1}}\ln\left(1+\rho^{-2}\right)d\rho/\rho \\ 
&=& \ln(\rho_{1})^2+O(1)
\end{eqnarray*} 
since if $\rho\in[0,1],$ then $\ln(1+\rho^2)\leq \rho^2$. Also, using the preceding choice of $\rho_{0}<\widetilde{\rho}_{0}<\sqrt{2}<\widetilde{\rho}_{1}<\rho_{1}$, 
\begin{eqnarray*} 
\int_{\rho_{0}}^{\rho_{1}}\ln(1+\cos(2\theta_{\rho}))\frac{d\rho}{\rho} &=& \int_{\rho_{0}}^{\widetilde{\rho}_{0}}\ln(1+\cos(2\theta_{\rho}))\frac{d\rho}{\rho}
+\int_{\widetilde{\rho}_{0}}^{\widetilde{\rho}_{1}}\ln(1+\cos(2\theta_{\rho}))\frac{d\rho}{\rho}+\int_{\widetilde{\rho}_{1}}^{\rho_{1}} \ln(1+\cos(2\theta_{\rho}))\frac{d\rho}{\rho} \\ 
&=& \int_{\widetilde{\rho}_{0}}^{\widetilde{\rho}_{1}}\ln(2+O(\theta_{\rho}^2))\frac{d\rho}{\rho} + O(\ln(\omega(d))) \quad \mathrm{since} \,\,\, 0\leq \ln(1+\cos(2\theta_{\rho}))\leq \ln(2) \\ 
&=& \ln(2)\ln\left(\widetilde{\rho}_{1}/\widetilde{\rho}_{0}\right)+O\left(\int_{\widetilde{\rho}_{0}}^{\widetilde{\rho}_{1}} \frac{\theta_{\rho}^2}{\rho}d\rho\right)+O(\ln(\omega(d))) \\
&=& \ln(2)\ln\left(\rho_{1}/\rho_{0}\right)+O(1) \quad \mathrm{as}\,\, \,\mathrm{in}\, \,\,\mathrm{(ii)}. 
\end{eqnarray*} 
Finally, 
\begin{eqnarray*} 
\int_{\rho_{0}}^{\rho_{1}}\int_{\theta_{\rho}}^{\frac{\pi}{2}-\theta_{\rho}}  \frac{dxdy}{xy} &=& 2\ln\left(d\pi/2\right)\ln\left(\rho_{1}/\rho_{0}\right)+2(\ln(\rho_{1})^2-\ln(\rho_{0})^2)-3\ln(\rho_{1})^2+2\ln(2)\ln\left(\rho_{1}/\rho_{0}\right)+O(1) \\ 
&=& 3\ln(\rho_{1})\ln(d\pi)-\frac{3}{2}\ln(\rho_{1})^2+O(1) \,\, \mathrm{by} \, \mathrm{Remark} \, \ref{remark38} \\ 
&=& 3\ln(d)^2/2+3\ln\left(\pi\right)\ln(d)+O(1) \,\, \mathrm{by}\,\mathrm{Corollary}\, \ref{lnest}.  
\end{eqnarray*} 
Adding up the contributions from all three parts gives the result.  \hfill $\Box$ 

\begin{cor} \label{upperlimit}
The expected area of the amoeba of a complex plane curve satisfies 
\[  \limsup_{d\rightarrow +\infty} \frac{\mathbb{E}_{d}(Vol(\mathcal{A}))}{\ln(d)^2} \leq \frac{3}{2}\cdot\] 
\end{cor} 

\noindent\textit{Proof:} This follows from Corollary \ref{expbound} and Theorem \ref{upperboundasymptotic}.  \hfill $\Box$

\section{Fubini-Study Charts of Plane Curves} \label{barrier} 

\subsection{Fubini-Study Charts} \label{allthecharts}
Let $D(0,r)$ and $\mathring{D}(0,r)$ be the closed and open disks of radius $r>0$ in $\mathbb{C}$ respectively. For every $\kappa\in (0,1]$, let $\mathbb{B}(\kappa)=D(0,1)\times D(0,\kappa)\subset \mathbb{C}^2\cong U_{0}\subset \mathbb{CP}^2$ and $\partial^{-}\mathbb{B}(\kappa)=\partial D(0,1)\times D(0,\kappa)$, $\partial^{+}\mathbb{B}(\kappa)=D(0,1)\times \partial D(0,\kappa)$. Let $\pi_{1}: \mathbb{B}(\kappa)\rightarrow D(0,1)$ and $\pi_{2}:\mathbb{B}(\kappa)\rightarrow D(0,\kappa)$ be the projection maps onto each factor. We finally set $\mathbb{B}=\mathbb{B}(1)$.  

\begin{defn}  \label{markedchart}
A triple $(B, \partial^{-}B, \partial^{+}B)$ is called a marked chart of a four-dimensional manifold $M$ if there exists a neighborhood $V$ of $\mathbb{B}$ in $\mathbb{C}^2$, a neighborhood $U$ of $B$ in $M$ and a diffeomorphism $\varphi: V\rightarrow U$ which maps $(\mathbb{B}, \partial^{-}\mathbb{B}, \partial^{+}\mathbb{B})$ onto $(B, \partial^{-}B, \partial^{+}B)$. 
\end{defn} 

\begin{defn} 
A marked chart $(B, \partial^{-}B, \partial^{+}B)$ is a submanifold chart of a surface $S$ in a four-manifold $M$ if the diffeomorphism $\varphi: V\rightarrow U$ mapping  $(\mathbb{B}, \partial^{-}\mathbb{B}, \partial^{+}\mathbb{B})$ to $(B, \partial^{-}B, \partial^{+}B)$ can be chosen in such a way that $\varphi^{-1}(B\cap S)$ is the graph of a map $\pi_{1}(\mathbb{B})=D(0,1)\rightarrow \mathring{D}(0,1)=\pi_{2}(B\setminus \partial^{+}B)$. 
\end{defn} 

\begin{defn} \label{markedFSchart}
A marked chart of $M=\mathbb{CP}^2$ is called a degree $d$ marked Fubini-Study chart of capacity $\kappa\in (0,1]$ if the diffeomorphism $\varphi$ given by Definition \ref{markedchart} can be chosen to be $\Psi\circ h_{d}$ and defined on $(\mathbb{B}(\kappa),\partial^{-}\mathbb{B}(\kappa),\partial^{+}\mathbb{B}(\kappa))$, where $\Psi\in PU(3,\mathbb{C})$ is an isometry of $\mathbb{CP}^2$ and $h_{d}$  a homothety of $\mathbb{C}^2$ of scaling factor $\displaystyle{ \frac{1}{\sqrt{6d}}}$. The point $\varphi(0)$ is called its center. 
\end{defn} 
We thus set, for every $d>0$, 
\begin{eqnarray*}
 \mathbb{B}_{d}(\kappa)&=&D\left(0,\frac{1}{\sqrt{6d}}\right)\times  D\left(0,\frac{\kappa}{\sqrt{6d}}\right)\subset \mathbb{C}^2\cong U_{0} \subset \mathbb{CP}^2 \\ 
\partial^{+}\mathbb{B}_{d}(\kappa)&=&D\left(0,\frac{1}{\sqrt{6d}}\right)\times \partial D\left(0,\frac{\kappa}{\sqrt{6d}}\right), \\
\partial^{-} \mathbb{B}_{d}(\kappa)&=&\partial D\left(0,\frac{1}{\sqrt{6d}}\right)\times  D\left(0,\frac{\kappa}{\sqrt{6d}}\right) 
\end{eqnarray*} 
and $\mathbb{B}_{d}=\mathbb{B}_{d}(1)$. These sets are contained in $\mathbb{C}^2\subset \mathbb{CP}^2$ and a degree $d$ marked Fubini-Study chart of capacity $\kappa$ is the image of the triple $(\mathbb{B}_{d}(\kappa), \partial^{-}\mathbb{B}_{d}(\kappa), \partial^{+}\mathbb{B}_{d}(\kappa))$ under some isometry of $\mathbb{CP}^2$. We also set $\displaystyle{\rho_{d}=\arctan\left( \frac{2}{\sqrt{d}}\right)}$ and for every $z\in \mathbb{CP}^2$, denote by $B_{FS}(z,\rho_{d})$ the Fubini-Study ball with center $z$ and radius $\rho_{d}$. 

\begin{rem} \label{44}
We need, in the proof of Theorem \ref{exponential}, $\rho_{d}$ to satisfy $\displaystyle{ \lim_{d\rightarrow +\infty} (N_{d} \sin^{4} \rho_{d})>1}$ and with our choice of $\rho_{d}$, this limit equals 8.  
\end{rem} 

\begin{lem} \label{ballcontain} 
Let $(B, \partial^{-}B, \partial^{+}B)$ be a degree $d$ marked Fubini-Study chart of capacity $\kappa\in (0,1]$ with center $z\in \mathbb{CP}^2$. Then, $B\subset B_{FS}(z,\rho_{d})$. 
\end{lem} 

\noindent \textit{Proof:} By Definition \ref{markedFSchart}, it suffices to prove that $\mathbb{B}_{d}\subset B_{FS}(0,\rho_{d})$ since for every $\kappa\in (0,1]$, $\mathbb{B}_{d}(\kappa)\subset \mathbb{B}_{d}$. But for all $d>0$, 
\[ \mathbb{B}_{d}\subset B\left(0,\frac{1}{\sqrt{3d}}\right) = B_{FS}\left(0,\arctan\left(\frac{1}{\sqrt{3d}}\right)\right)\subset B_{FS}(0,\rho_{d}) \]
by Lemma \ref{FS_and_Euclidean_distance}.     (Here, $B(0,r)$ denotes an Euclidean ball with center $0$ and radius $r$.) \hfill $\Box$ 

\begin{defn} \label{ballaverage}
For every $z\in \mathbb{CP}^2$ and $P\in H^{0}(\mathbb{CP}^2, L^d)$, the Fubini-Study ball-average norm of $P$ around $z$ is 
\[ ||P||^2_{B_{FS}(z,\rho_{d})} =\frac{1}{Vol(B_{FS}(z,\rho_{d}))} \int_{B_{FS}(z,\rho_{d})} ||P(t)||_{FS}^2 dVol_{FS}(t). \] 
\end{defn}

\subsection{Linear Polynomials} \label{linearpoly}
For each $d>0$, we set $P_{d}(Z,X,Y)=\sqrt{dN_{d}}YZ^{d-1}\in H^{0}(\mathbb{CP}^2,L^d)$, where $\displaystyle{N_{d}=\frac{(d+1)(d+2)}{2}=\dim_{\mathbb{C}} H^{0}(\mathbb{CP}^2,L^d)}$. We need some estimates about these reference polynomials that we establish in this subsection.

\begin{pro} \label{PolyFSestimates}
For each $d>0$, $||P_{d}||_{FS}=1$ and 
\begin{enumerate} 
\item for all $z\in \mathbb{CP}^2$,   $\displaystyle{||P_{d}(z)||_{FS}\leq \sqrt{dN_{d}} \cos (d_{FS}(0,z))^{d-1}}$, \\
\item for all $\kappa\in (0,1]$ and all $z\in \partial^{+} \mathbb{B}_{d}(\kappa)$, $\displaystyle{||P_{d}(z)||_{FS}\geq \kappa\sqrt{\frac{N_{d}}{6}} e^{-1/6}}$, \\ 
\item for all $\Psi \in PU(3,\mathbb{C})$, $\displaystyle{ |\langle P_{d},P_{d}\circ \Psi^{-1} \rangle | \leq d N_{d} \left( \frac{1+\cos\left(d_{FS}(0,\Psi(0))\right)}{2}\right)^{d-1}} $. \\
\end{enumerate} 
\end{pro} 

\noindent \textit{Proof:}  By Lemma \ref{FSnorm} with $n=2$ and multiindex $(0,1,d-1)$, 
\[ ||YZ^{d-1}||_{FS}^2 = \frac{(d-1)! 2!}{(d+2)!} = \frac{2}{(d+2)(d+1)d}, \]
so that $||P_{d}||_{FS}=1$. 

In order to prove (1), we observe that by equation \eqref{sectionnorm} of subsection \ref{FSmetric}, for every $z=(z_{1},z_{2})\in \mathbb{C}^2$, 
\[  ||P_{d}(z)||^2_{FS} = dN_{d} \frac{|z_{2}|^2}{(1+||z||^2)^{d}} \leq \frac{dN_{d}}{(1+||z||^2)^{d-1}}\cdot \] 
Since $\tan(d_{FS}(0,z))=d(0,z)=||z||$ by Lemma \ref{FS_and_Euclidean_distance}, we get $\displaystyle{ \frac{1}{1+||z||^2}=\cos^2(d_{FS}(0,z))}$, which completes the proof of (1). 

For (2), we note that for every $z=(z_{1},z_{2})\in \partial^{+} \mathbb{B}_{d}(\kappa)$,  $\displaystyle{|z_{2}|= \frac{\kappa}{\sqrt{6d}}}$ and $\displaystyle{||z||\leq \frac{1}{\sqrt{3d}}}$, so that 
\begin{eqnarray*} 
||P_{d}(z)||^2_{FS}&=& \frac{\kappa^2 dN_{d}}{6d (1+||z||^2)^{d}}  \\ 
&=&  \frac{\kappa^2 N_{d}}{6} \exp(-d \ln(1+||z||^2))   \\
&\geq&  \frac{\kappa^2 N_{d}}{6} e^{-1/3},
\end{eqnarray*} 
which implies (2). 

Finally, for (3), 
\begin{eqnarray*} 
 |\langle P_{d},P_{d}\circ \Psi^{-1} \rangle | &\leq& \frac{1}{Vol_{FS}(\mathbb{CP}^{2})}\int_{\mathbb{CP}^2} |h_{d}(P_{d}(z),P_{d}(\Psi^{-1}(z)))|dVol_{FS}(z) \\
&\leq&  \frac{1}{Vol_{FS}(\mathbb{CP}^{2})}\int_{\mathbb{CP}^2} ||P_{d}(z)||_{FS} ||P_{d}(\Psi^{-1}(z))||_{FS} dVol_{FS}(z) \\ 
&\leq& \frac{dN_{d}}{Vol_{FS}(\mathbb{CP}^{2})} \int_{\mathbb{CP}^2} \left(\cos(d_{FS}(0,z))\cos(d_{FS}(0,\Psi^{-1}(z)))\right)^{d-1} dVol_{FS}  
\end{eqnarray*} 
by part (1). Since $\displaystyle{\cos(a)\cos(b)=\frac{1}{2}(\cos(a+b)+\cos(a-b))\leq \frac{1}{2}(1+\cos(a+b))}$ and by the triangle inequality $\pi>d_{FS}(0,z)+d_{FS}(0,\Psi^{-1}(z))=d_{FS}(0,z)+d_{FS}(\Psi(0),z)\geq d_{FS}(0,\Psi(0))$, the proof of (3) gets complete.  \hfill $\Box$ 

\begin{cor}  \label{normestimate}
For every $d>0$, 
\begin{enumerate} 
\item $\displaystyle{  \frac{1}{\pi^2 \left( 1+4/d\right)^3} \int_{B\left(0,2\right)} |t_{2}|^2 e^{-||t||^2}dVol(t) \leq \sin^{4}(\rho_{d}) ||P_{d}||^2_{B_{FS}(0,\rho_{d})} \leq 1} $, 
\item for all $z\in \mathbb{CP}^2$ such that $d_{FS}(0,z)>\rho_{d}$,  $\displaystyle{||P_{d}||_{B_{FS}(z,\rho_{d})} \leq \sqrt{dN_{d}} \cos(d_{FS}(0,z)-\rho_{d})^{d-1}} $.  
\end{enumerate} 
\end{cor}

\noindent \textit{Proof:}  To get the first inequality of (1), we observe that
\begin{eqnarray*} 
||P_{d}||^2_{B_{FS}(0,\rho_{d})} &=& \frac{1}{Vol(B_{FS}(0,\rho_{d}))} \int_{B_{FS}(0,\rho_{d})} ||P_{d}(t)||^2_{FS} dVol_{FS}(t)  \\ 
&=& \frac{dN_{d}}{Vol_{FS}(\mathbb{CP}^2)\sin^{4}(\rho_{d})} \int_{B\left(0,2/\sqrt{d}\right)} \frac{|z_{2}|^2}{(1+||z||^2)^{d+3}} dVol(z) \, \, \mathrm{by} \, \eqref {sectionnorm}   \, \mathrm{and}\, \mathrm{Lemmas} \, \ref{volcompare}-\ref{sinvolume} \\ 
 &=& \frac{N_{d}}{d^2 Vol_{FS}(\mathbb{CP}^2)\sin^{4}(\rho_{d})}  \int_{B\left(0,2\right)}  \frac{|t_{2}|^2 \exp\left( -d \ln\left(1+||t||^2/d\right)\right) }{ \left(1+||t||^2/d\right)^3} dVol(t) \\
&\geq& \frac{1}{\pi^2 \left(1+4/d\right)^3\sin^{4}(\rho_{d})} \int_{B\left(0, 2\right)} |t_{2}|^2 \exp(-||t||^2) dVol(t). 
\end{eqnarray*}  
Here, the substitution $t=\sqrt{d}z$ was used. This establishes the first inequality in (1). The second inequality in (1) follows from the equality $||P_{d}||^2_{FS}=1$ given by Proposition \ref{PolyFSestimates}, the inequality $\displaystyle{||P_{d}||^2_{B_{FS}(0,\rho_{d})}\leq \frac{1}{Vol(B_{FS}(0,\rho_{d}))} \int_{\mathbb{CP}^2} ||P_{d}(t)||^2_{FS} dVol_{FS}(t)}$ and Lemma \ref{sinvolume}.

To get (2), we observe that
\begin{eqnarray*} 
||P_{d}||^2_{B_{FS}(z,\rho_{d})}\leq \frac{d N_{d}}{Vol(B_{FS}(z,\rho_{d}))} \int_{B_{FS}(z,\rho_{d})} \cos(d_{FS}(0,t))^{2(d-1)} dVol_{FS}(t)
\end{eqnarray*}  
by the first part of Proposition \ref{PolyFSestimates}, while by the triangle inequality, $\pi/2\geq d_{FS}(0,t)\geq d_{FS}(0,z)-\rho_{d}$ for any $t\in B_{FS}(z,\rho_{d})$. \hfill $\Box$ 

\subsection{An Effective Criterion} 
The zero scheme of the reference polynomial $P_{d}(Z,X,Y)=\sqrt{dN_{d}}YZ^{d-1}$ in $\mathbb{CP}^2$ is the union of the axes $Y=0$ and a $(d-1)$-fold thick line supported on the ``line at infinity'' $Z=0$. In an affine neighborhood of the origin $[1:0:0]$, this zero locus agrees with the $x$-axis, hence $(\mathbb{B}_{d}(\kappa), \partial^{-} \mathbb{B}_{d}(\kappa), \partial^{+} \mathbb{B}_{d}(\kappa))$ is a  submanifold chart for it whatever $\kappa\in(0,1]$ is and sustains being a  submanifold chart for any small enough perturbation of it. The goal of this subsection is to quantify the latter, to deduce an effective criterion that will be used in the proof of Theorem \ref{exponential}, see Corollary \ref{corexplicitsupupperbound}. 

\begin{pro} \label{supupperbound}
For every $d>0$, every $\kappa\in (0,1]$ and every $Q\in H^{0}(\mathbb{CP}^2,L^d)$, 
\[  \sup_{z\in \partial^{+} \mathbb{B}_{d}(\kappa)} \left( \frac{||Q(z)||^2_{FS} }{||P_{d}(z)||^2_{FS}} \right) \leq  \frac{3^2 2^4 e \left(1+4/d\right)^3 }{\inf_{z\in \partial^{+} \mathbb{B}_{d}(\kappa)}  ||P_{d}(z)||^2_{FS} }  ||Q||^2_{B_{FS}(0,\rho_{d})}. \] 
\end{pro} 

\noindent \textit{Proof:}  For all $z\in \partial^{+}\mathbb{B}_{d}(\kappa)$ and all $Q\in H^{0}(\mathbb{CP}^2,L^d)$, by the mean value inequality,
\begin{eqnarray*} 
   \frac{||Q(z)||^2_{FS} }{||P_{d}(z)||^2_{FS}} &=& \frac{|Q(z)|^2}{|P_{d}(z)|^2} \leq  \frac{1}{|P_{d}(z)|^2 Vol(B(z,1/\sqrt{3d}))} \int_{B\left(z,1/\sqrt{3d}\right)} |Q(t)|^2 dVol(t) \\
&\leq&  \frac{1}{||P_{d}(z)||^2_{FS} Vol(B(z,1/\sqrt{3d}))} \int_{B\left(z,1/\sqrt{3d}\right)} ||Q(t)||^2_{FS} \left( \frac{1+||t||^2}{1+||z||^2}\right)^{d} dVol(t).  \\
\end{eqnarray*}
Moreover, for all $t\in B(z,1/\sqrt{3d})$ and all $z\in \partial^{+}\mathbb{B}_{d}(\kappa)$, $||z||\leq 1/\sqrt{3d}$, $||t||-||z||\leq 1/\sqrt{3d}$ and 
\begin{eqnarray*} 
 \left( \frac{1+||t||^2}{1+||z||^2}\right)^{d} &=& \left( 1+ \frac{(||t||-||z||)(||t||-||z||+2||z||)}{1+||z||^2}\right)^{d} \\
&\leq& \left( 1+\frac{1}{d}\right)^{d}  \\ 
&\leq& e. 
\end{eqnarray*} 
For any $z\in \partial^{+}\mathbb{B}_{d}(\kappa)$, the ball $B(z,1/\sqrt{d})$ is contained in $B(0,2/\sqrt{d})$, while $Vol(B(0,2/\sqrt{d}))=3^2 2^4 Vol(B(z,1/\sqrt{3d}))$. We deduce from Lemma \ref{volcompare} that 
\[ \frac{||Q(z)||^2_{FS}}{||P_{d}(z)||^2_{FS}} \leq \frac{3^2 2^4 e}{||P_{d}(z)||^2_{FS} Vol(B(0,2/\sqrt{d}))} \int_{B\left(0, 2/\sqrt{d}\right)} ||Q(t)||^2_{FS}(1+||t||^2)^3 dVol_{FS}(t). \] 
Since $\displaystyle{ Vol(B(0,2/\sqrt{d}))= \int_{B\left(0,2/\sqrt{d}\right)} (1+||t||^2)^3 dVol_{FS}(t) \geq Vol(B_{FS}(0, \rho_{d}))}$, we get 
\[  \sup_{z\in \partial^{+} \mathbb{B}_{d}(\kappa)} \left( \frac{||Q(z)||^2_{FS} }{||P_{d}(z)||^2_{FS}} \right) \leq  \frac{3^2 2^4 e\left(1+4/d\right)^3 }{\inf_{z\in \partial^{+} \mathbb{B}_{d}(\kappa)}  ||P_{d}(z)||^2_{FS} }  ||Q||^2_{B_{FS}(0,\rho_{d})}.  \]  
\hfill $\Box$ 

\begin{cor} \label{explicitsupupperbound}
For every $d>0$, every $\kappa\in (0,1]$ and every $Q\in H^{0}(\mathbb{CP}^2,L^{d})$, 
\[  N_{d} \sup_{z\in \partial^{+} \mathbb{B}_{d}(\kappa)} \left( \frac{||Q(z)||^2_{FS} }{||P_{d}(z)||^2_{FS}} \right) \leq \frac{c_{1}(d)}{\kappa^2} ||Q||^2_{B_{FS}(0,\rho_{d})} \] 
with $\displaystyle{c_{1}(d)=3^3 2^5  e^{4/3} \left( 1+\frac{4}{d} \right)^3 \xrightarrow[d\rightarrow +\infty]{} 3^3 2^5 e^{4/3} }$. 
\end{cor} 

\noindent \textit{Proof:}   The result follows from Proposition \ref{PolyFSestimates} part (2) and Proposition \ref{supupperbound}.  \hfill $\Box$ 

\begin{pro} \label{homotopy}
For every $d>0$, every $\kappa\in (0,1]$ and every $(a,Q)\in \mathbb{C}\times H^{0}(\mathbb{CP}^2,L^d)$ such that $\displaystyle{ |a|^2 > \sup_{z\in \partial^{+} \mathbb{B}_{d}(\kappa)} \left( \frac{||Q(z)||^2_{FS}}{||P_{d}(z)||^2_{FS}} \right)}$, the triple $(\mathbb{B}_{d}(\kappa), \partial^{-} \mathbb{B}_{d}(\kappa), \partial^{+} \mathbb{B}_{d}(\kappa))$ is a submanifold chart for $V_{aP_{d}+Q}$. 
\end{pro} 

\noindent \textit{Proof:}  We already observed that $(\mathbb{B}_{d}(\kappa), \partial^{-} \mathbb{B}_{d}(\kappa), \partial^{+} \mathbb{B}_{d}(\kappa))$ is a submanifold chart for $V_{P_{d}}=V_{aP_{d}}$. Now, for all $a\in \mathbb{C}$ satisfying the hypothesis, all $t\in[0,1]$ and all $z\in \partial^{+}\mathbb{B}_{d}(\kappa)$, the evaluation $(aP_{d}+tQ)(z)$ does not vanish since
\begin{eqnarray*} 
|| ( aP_{d}+tQ)(z)||_{FS} &\geq&  |a| ||P_{d}(z)||_{FS}-||Q(z)||_{FS}  \\ 
&=& ||P_{d}(z)||_{FS} \left( |a|- \frac{||Q(z)||_{FS}}{||P_{d}(z)||_{FS}}\ \right)  \\ 
&>& 0.  
\end{eqnarray*} 
The result then follows by homotopy, since the number of roots of the polynomial $aP_{d}+tQ$ in each disc $\{z_{1}\}\times D(0,\kappa/\sqrt{6d})\subset \mathbb{B}_{d}(\kappa)$, $z_{1}\in D(0,1/\sqrt{6d})$, does not depend on $t\in [0,1]$ by Cauchy formula and thus equals one.
\hfill $\Box$

\begin{cor} \label{corexplicitsupupperbound}
For every $\kappa\in (0,1]$ and every $(a,Q)\in \mathbb{C}\times H^{0}(\mathbb{CP}^2,L^d)$ such that $\displaystyle{ |a|^2 > \frac{c_{1}(d)}{\kappa^{2}N_{d}} ||Q||^2_{B_{FS}(0,\rho_{d})} }$, the triple $(\mathbb{B}_{d}(\kappa), \partial^{-} \mathbb{B}_{d}(\kappa), \partial^{+} \mathbb{B}_{d}(\kappa))$ is a submanifold chart for $V_{aP_{d}+Q}$. 
\end{cor} 

\noindent \textit{Proof:} The result follows from Proposition \ref{homotopy} and Corollary \ref{explicitsupupperbound}.  \hfill $\Box$

\subsection{The Ball-Average Norm} \label{ballnorm}
We now investigate, as a random variable, the ball-average norm given by Definition  \ref{ballaverage}    which appeared in Corollary \ref{corexplicitsupupperbound}. Let thus $Y_{d}$ denote the random variable $P\in H^{0}(\mathbb{CP}^2,L^d)\mapsto ||P||_{B_{FS}(0,\rho_{d})}\in \mathbb{R}_{+}$. Given an orthogonal direct sum decomposition $H^{0}(\mathbb{CP}^2,L^{d})=H\oplus H^{\perp}$, we denote by $\pi_{H}$ (resp. $\pi_{H^{\perp}}$) the orthogonal projection to $H$ (resp. $H^{\perp}$) and by $s_{H}$ (resp. $s_{H^{\perp}}$) the orthogonal reflection with axis $H$ (resp. $H^{\perp}$), so that $Fix(s_{H})=H=Im(\pi_{H})$. We also set $Y_{d,H}=Y_{d}\circ \pi_{H}$ and $Y_{d,H^{\perp}}=Y_{d}\circ \pi_{H^{\perp}}$. 

\begin{lem} \label{directsum} 
For every $d>0$, $\mathbb{E}(Y_{d}^2)=N_{d}$ and for every orthogonal decomposition $H^{0}(\mathbb{CP}^2,L^d)=H\oplus H^{\perp}$, $\mathbb{E}(Y_{d}^2)=\mathbb{E}(Y_{d,H}^2)+\mathbb{E}(Y_{d,H^{\perp}}^2)$. 
\end{lem} 

\noindent \textit{Proof:}  For all $d>0$,
\begin{eqnarray*} 
\mathbb{E}(Y_{d}^2)&=& \int_{H^{0}(\mathbb{CP}^2,L^d)} \left(\frac{1}{Vol(B_{FS}(0,\rho_{d}))} \int_{B_{FS}(0,\rho_{d})} ||P(t)||^{2}_{FS} dVol_{FS}(t)\right) d\mu(P) \\ 
&=& \frac{1}{Vol(B_{FS}(0,\rho_{d}))} \int_{B_{FS}(0,\rho_{d})} \left( \int_{H^{0}(\mathbb{CP}^2,L^d)} ||P(t)||^2_{FS} d\mu(P) \right) dVol_{FS}(t) 
\end{eqnarray*} 
The average  $\displaystyle{ \int_{H^{0}(\mathbb{CP}^2,L^d)} ||P(t)||^2_{FS} d\mu(P)}$ does not depend on $t\in \mathbb{CP}^2$ since the isometries of $\mathbb{CP}^2$ act transitively. We deduce that
\begin{eqnarray*} 
\mathbb{E}(Y_{d}^2)&=& \int_{H^{0}(\mathbb{CP}^2,L^d)} ||P(t)||^2_{FS} d\mu(P)   \qquad \mathrm{for} \, \mathrm{any} \, t\in \mathbb{CP}^2   \\ 
&=& \frac{1}{Vol_{FS}(\mathbb{CP}^2)} \int_{\mathbb{CP}^2} \int_{H^{0}(\mathbb{CP}^2,L^{d})} ||P(t)||^2_{FS} d\mu(P)  dVol_{FS}(t) \\ 
&=& \mathbb{E}(||P||^2_{FS}) \\
&=& Var(\mu) \\
&=&N_{d}. 
\end{eqnarray*} 
Hence the first claim. 

Next, we observe that for any $P\in H$ and $Q\in H^{\perp}$,
\begin{eqnarray*} 
\frac{1}{2}(Y_{d}^2(P+Q)+Y_{d}^2(P-Q)) &=&  \frac{1}{2} \left(||P+Q||^2_{B_{FS}(0,\rho_{d})}+||P-Q||^2_{B_{FS}(0,\rho_{d})}\right) \\
&=& ||P||^2_{B_{FS}(0,\rho_{d})} +||Q||^2_{B_{FS}(0,\rho_{d})}. 
\end{eqnarray*}
by the parallelogram rule for Hermitian norms. This shows that 
\begin{equation}\label{symmetry}
\frac{1}{2}(Y_{d}^2+Y_{d}^2\circ s_{H}) =Y^2_{d,H}+Y^2_{d,H^{\perp}}
\end{equation}
 and consequently, 
\[ \frac{1}{2}(\mathbb{E}(Y_{d}^2)+\mathbb{E}(Y_{d}^2\circ s_{H}))=\mathbb{E}(Y_{d,H}^2)+\mathbb{E}(Y_{d,H^{\perp}}^2).  \] 
The result now follows from the fact that $s_{H}$ preserves the measure $\mu$ on $H^{0}(\mathbb{CP}^2,L^d)$. \hfill $\Box$ 

\begin{lem} \label{expectedestimate}
Let $P_{1},\ldots, P_{N}\in H^{0}(\mathbb{CP}^2,L^d)$ be linearly independent, $H=Span(P_{1},\ldots,P_{N})$ and $C=\left(\langle P_{i},P_{j}\rangle \right)_{1\leq i,j\leq N}$. Then, 
\[ \sup_{1\leq i\leq N} \left( \frac{||P_{i}||^2_{B_{FS}(0,\rho_{d})}}{||P_{i}||^2_{FS}}\right) \leq \mathbb{E}(Y^2_{d,H})\leq ||C^{-1}||\sum_{i=1}^{N} ||P_{i}||^{2}_{B_{FS}(0,\rho_{d})}. \]
\end{lem} 

\noindent \textit{Proof:}
\[ \mathbb{E}(Y_{d,H}^2)= \frac{1}{Vol(B_{FS}(0,\rho_{d}))} \int_{B_{FS}(0,\rho_{d})} \left( \int_{H} ||P(t)||_{FS}^2 d\mu_{H}(P) \right) dVol_{FS}(t), \] 
where $\mu_{H}$ denotes the Gaussian measure of $H$ given by Definition \ref{Gaussiandef}. 
Say $t\in B_{FS}(0,\rho_{d})$ and $\sigma_{t,H}$ is a vector of norm $1$ in $H$ which is orthogonal to the hyperplane of polynomials vanishing at $t$. Then,
\[ \int_{H}||P(t)||^2_{FS} d\mu_{H}(P) = ||\sigma_{t,H}(t)||^2_{FS} \geq \frac{||P_{i}(t)||^2_{FS}}{||P_{i}||^2_{FS}} \] 
for every $i\in \{1,\ldots,N\}$. This implies the first inequality. 

In order to prove the second inequality, we write $P=\sum_{i=1}^{N} X_{i}(P) P_{i}$ and deduce that for every $t\in B_{FS}(0,\rho_{d})$, 
\begin{eqnarray*} 
 \int_{H} ||P(t)||_{FS}^2 d\mu_{H}(P) & = &\sum_{1\leq i,j \leq N} Cov(X_{i},X_{j}) h_{FS}(P_{i}(t),P_{j}(t))  \\ 
&\leq& ||C^{-1}|| \left( \sum_{j=1}^{N} ||P_{j}(t)||^2_{FS} \right)  
\end{eqnarray*} 
by Lemma \ref{covariance} and Cauchy-Schwarz inequality. Averaging over the ball $B_{FS}(0,\rho_{d})$ provides the second inequality. \hfill $\Box$ 

\begin{lem} \label{lemmaxA}
For every  $d>0$, every orthogonal decomposition $H^{0}(\mathbb{CP}^2,L^d)=H\oplus H^{\perp}$ and every event $A\subset H^{0}(\mathbb{CP}^2,L^d)$ invariant under $s_{H}$, there exists $z_{A}\in B_{FS}(0,\rho_{d})$ such that $\mathbb{E}_{A}(Y_{d,H}^2)\leq \mathbb{E}_{A}(||Z||^2_{FS})$, with $Z: P\in H^{0}(\mathbb{CP}^2,L^{d})\mapsto P(z_{A})\in L^{d}|_{z_{A}}$. 
\end{lem} 
Recall that $\mathbb{E}_{A}$ denotes the conditional expectation over $A$, see subsection \ref{conditionalexpectation}. 

\noindent \textit{Proof:} By definition  and  Fubini's theorem,
\begin{eqnarray*} 
\mathbb{E}_{A}(Y_{d,H}^2) &=& \frac{1}{\mu(A)} \int_{A} \left( \frac{1}{Vol(B_{FS}(0,\rho_{d}))} \int_{B_{FS}(0,\rho_{d})} ||\pi_{H}(P)(t)||^2_{FS} dVol_{FS}(t) \right) d\mu(P) \\ 
&=&   \frac{1}{Vol(B_{FS}(0,\rho_{d}))}  \int_{B_{FS}(0,\rho_{d})}  \left(  \frac{1}{\mu(A)}  \int_{A} ||\pi_{H}(P)(t)||^2_{FS}d\mu(P)  \right)  dVol_{FS}(t).
\end{eqnarray*} 
Therefore, there exists $z_{A}\in B_{FS}(0,\rho_{d})$ such that 
\[ \mathbb{E}_{A}(Y_{d,H}^2) \leq \frac{1}{\mu(A)}  \int_{A} ||\pi_{H}(P)(z_{A})||^2_{FS}d\mu(P) = \mathbb{E}_{A}(||Z\circ \pi_{H}||^2_{FS}). \] 
Now, by the parallelogram rule for Hermitian norms
\[ \frac{1}{2}(||Z||^2_{FS}+||Z\circ s_{H}||^2_{FS}) =||Z\circ \pi_{H}||^2_{FS} +||Z\circ \pi_{H^{\perp}}||^2_{FS}, \quad \mathrm{compare} \,\eqref{symmetry}. \] 
Taking conditional expectations over $A$, the two summands on the left hand side contribute equally since $s_{H}$ preserves $\mu$ and $s_{H}(A)=A$, so that we get $\mathbb{E}_{A}(||Z||^{2}_{FS})$ on the left. We deduce that $ \mathbb{E}_{A}(||Z\circ \pi_{H}||^2_{FS}) \leq 
\mathbb{E}_{A}(||Z||^{2}_{FS})$ and the result. \hfill  $\Box$

\section{Asymptotic lower estimate}\label{lower}

\subsection{The probabilistic tool} 
A substantial part of the proof of Theorem \ref{exponential} relies on a probabilistic phenomenon which we now establish. Namely, the proof provides $N$ pairs of real random variables $(X_{i},Y_{i})$ together with the event $\exists i\in \{1,\ldots,N\}$ such that $X_{i}\geq x$ and $Y_{i}\leq y$ for some $x,y\in \mathbb{R}$. The aim of this section is to bound from below the probability of this event by $1-c\gamma^{N}$ with $c>0$ and $0\leq \gamma <1$. If these pairs of random variables were independent and identically distributed, such a lower estimate would be given by the following proposition, which we recall for convenience. 

\begin{pro} \label{elementary} 
Suppose that $(X_{i},Y_{i})_ {1\leq i\leq N}$ are $N$ independent pairs of real random variables which follow the law of $(X,Y)$, with $X\perp Y$. Then, for all $x,y\in \mathbb{R}$, 
\[ \mathbb{P}\left(\exists i\in\{1,\ldots,N\}| X_{i}\geq x\,\,\mathrm{and}\,\,Y_{i}\leq y\right)=1-\left(1-\mathbb{P}(X\geq x)\mathbb{P}(Y\leq y)\right)^N. \] 
\end{pro} 

\noindent \textit{First proof:} 
\begin{eqnarray*} 
1-\mathbb{P}\left(\exists i\in\{1,\ldots,N\}| X_{i}\geq x \,\,\mathrm{and}\,\, Y_{i}\leq y \right) &=& \mathbb{P}\left(\forall i\in \{1,\ldots, N\}, X_{i}<x \,\,\mathrm{or}\,\, Y_{i}>y\right) \\ 
&=& \mathbb{P}\left(\bigcap_{i=1}^{N}(\{X_{i}<x\}\cup \{Y_{i}>y\})\right) \\   
&=& \prod_{i=1}^{N}\mathbb{P}\left(\{X_{i}<x\}\cup \{Y_{i}>y\}\right)  \\
&=& \prod_{i=1}^{N}\left(1-\mathbb{P}\left(\{X_{i}\geq x\} \cap \{Y_{i}\leq y\}\right)\right) \\ 
&=& \prod_{i=1}^{N} \left(1-\mathbb{P}(X_{i}\geq x)\mathbb{P}(Y_{i}\leq y)\right) \\
&=& \left(1-\mathbb{P}(X\geq x)\mathbb{P}(Y\leq y)\right)^N. 
\end{eqnarray*} 
\hfill $\Box$ 

\noindent\textit{Second proof:} In general, assuming only that $X_{i} \perp Y_{j}$ for all $i,j$, we get
\begin{eqnarray*} 
\mathbb{P}\left(\bigcap_{i=1}^{N}\left(\{X_{i}<x\}\cup\{Y_{i}>y\}\right)\right) &=& \mathbb{P}\left(\bigsqcup_{k=0}^{N} \bigsqcup_{I_{k}\subset \{1,\ldots,N\}, |I_{k}|=k} \left(\{X_{i}<x \, \mathrm{iff} i\in I_{k}\}\cap \{\forall i\notin I_{k}, Y_{i}>y\}\right)\right) \\
&=& \sum_{I\subset \{1,\ldots,N\}} \mathbb{P}\left(X_{i}<x \,\, \mathrm{iff}\,\, i\in I\right)\mathbb{P}\left(\forall i \notin I, Y_{i}>y\right). 
\end{eqnarray*} 
Now, assuming moreover that the $X_{i}$'s are i.i.d. and the $Y_{i}$'s are i.i.d., we deduce 
\begin{eqnarray*} 
\mathbb{P}(X_{i}<x \,\,\mathrm{iff}\,\, i\in I)&=&\mathbb{P}(X<x)^{\#I}\mathbb{P}(X\geq x)^{N-\#I} \qquad \mathrm{and} \\
\mathbb{P}(\forall i\notin I, Y_{i}>y)&=&\mathbb{P}(Y>y)^{N-\#I}.
\end{eqnarray*} 
Therefore, 
\begin{eqnarray*} 
\mathbb{P}\left(\bigcap_{i=1}^{N}(\{X_{i}<x\}\cup\{Y_{i}>y\})\right) &=& \sum_{k=0}^{N} {n\choose k} \mathbb{P}(X<x)^{k}\left(\mathbb{P}(X\geq x)\mathbb{P}(Y>y)\right)^{N-k} \\ 
&=& \left(\mathbb{P}(X<x)+\mathbb{P}(X\geq x)\mathbb{P}(Y>y)\right)^{N} \\ 
&=& \left(1-\mathbb{P}(X\geq x)\mathbb{P}(Y\leq y)\right)^{N}. 
\end{eqnarray*} 
\hfill $\Box$

We now aim to prove a similar estimate under weaker independence assumptions, but for random variables which can be compared to Gaussian ones with controlled covariance, see Corollary \ref{cormain}.

\begin{pro} \label{bernstein}
(\`{a} la Bernstein) Say $\mathbf{Z}\sim \mathcal{N}_{\mathbb{C}^{N}}(0,B)$. Then, 
for all $y>||B||$, 
\[ \mathbb{P}(||\mathbf{Z}||^2\geq Ny)\leq \frac{y^{N}}{\det(B)} \exp\left(N\left(1-\frac{y}{||B||}\right)\right). \] 
\end{pro} 
Recall that the normal law $\mathcal{N}_{\mathbb{C}^{N}}(0,B)$ has been introduced in Definition \ref{normallaw}. 

\noindent \textit{Proof:} Following Bernstein's trick, for all $y>||B||$ and all $\displaystyle{0<\lambda<\frac{1}{||B||}}$, 
\begin{eqnarray*} 
\mathbb{P}(||\mathbf{Z}||^2\geq Ny)&=& \mathbb{P}\left( e^{\lambda ||\mathbf{Z}||^2}\geq e^{\lambda Ny}\right) \\ 
&\leq& \mathbb{E}\left(e^{\lambda ||\mathbf{Z}||^2}\right) e^{-\lambda Ny}  \qquad \mathrm{by}\,\,\mathrm{Markov's} \,\, \mathrm{inequality} \\ 
&\leq& \frac{e^{-\lambda Ny}}{\det(I-\lambda B)} \qquad \mathrm{by}\,\, \mathrm{Lemma}\,\, \ref{centeredgaussian}. 
\end{eqnarray*}
Since $B$ is a positive definite Hermitian matrix, it is diagonalizable with real eigenvalues $b_{1},\ldots,b_{N}$ and orthogonal eigenvectors, so that $\displaystyle{||B||=\sup_{1\leq i\leq N} b_{i}}$. Choosing $\displaystyle{\lambda=\frac{1}{||B||}-\frac{1}{y}}$, we get
\[ \det(Id-\lambda B)=\prod_{i=1}^{N} (1-\lambda b_{i}) =\prod_{i=1}^{N} \left(1-\frac{b_{i}}{||B||}+\frac{b_{i}}{y}\right)\geq \frac{\prod_{i=1}^{N} b_{i}}{y^N}=\frac{\det B}{y^N}.\]
Therefore, 
\[ \mathbb{P}(||\mathbf{Z}||^2\geq Ny)\leq \frac{y^{N}}{\det(B)} \exp\left(N\left(1-\frac{y}{||B||}\right)\right). \]  
\hfill $\Box$

\begin{pro} \label{expboundprop} 
Say $\mathbf{Z}=(Z_{1},Z_{2},\ldots,Z_{N})\sim \mathcal{N}_{\mathbb{C}^{N}}(0,B)$ satisfies $V(Z_{i})=1$ for some $i\in \{1,\ldots,N\}$. Let $A$ be an event such that $\displaystyle{a:=\mathbb{E}_{A}\left(||\mathbf{Z}||^2/N \right)>2||B||||B^{-1}||=:2b}$. Then, 
\[ \mathbb{P}(A)\leq \left( \frac{e\cdot a}{a-2b}\right) \exp\left(1-\frac{a}{b}+\ln a\right)^{N}. \] 
\end{pro} 

\noindent\textit{Proof:} By Lemma \ref{firstlemma}, $\min(||B||,||B^{-1}||)\geq 1$ and Lemma \ref{expectedbound} applied to the random variable $\displaystyle{||\mathbf{Z}||^2/N}$ implies that for all $x< a $, 
\[ \mathbb{P}(A)\leq \frac{ \int_{x}^{+\infty} \mathbb{P}(||\mathbf{Z}||^2\geq Nz)dz}{a-x}. \] 
Then, Proposition \ref{bernstein} implies that for all $\displaystyle{x\in \left(||B||,a\right)}$, 
\[ \mathbb{P}(A)\leq \frac{1}{\det(B)\left(a-x\right)} \int_{x}^{+\infty}\exp\left(N\left(1-\frac{z}{||B||}+\ln z\right)\right)dz. \] 
Noting that for all $0<\alpha<1$ and all $\displaystyle{t\geq 1/\alpha}$, $\displaystyle{\ln(t)\leq \alpha t+\ln\left(1/\alpha\right)-1}$, we get  by choosing $\displaystyle{\alpha=\det(B)^{1/N}/x}$ that for all $z\geq x$, 
\[ \ln\left(\frac{z}{\det(B)^{1/N}}\right)\leq \frac{z}{x}+\ln\left(\frac{x}{\det(B)^{1/N}}\right)-1. \] 
Therefore,
\begin{eqnarray*} 
\mathbb{P}(A)&\leq& \frac{x^N  \exp\left(Nz\left(1/x-1/||B||\right)\right)|_{x}^{+\infty}}{\det(B)\left(a-x\right)N\left(1/x-1/||B||\right)}  \\ 
&\leq& \frac{ x^{N+1}||B|| \exp\left(N\left(1-x/||B||\right)\right)}{N\det(B)\left(a-x\right)(x-||B||)}\cdot 
\end{eqnarray*} 
Choosing $\displaystyle{x=a-||B||/N}$, which lies in the interval $\displaystyle{ \left(||B||,a\right)}$ by hypothesis, we deduce
\begin{eqnarray} \label{stronger}
\mathbb{P}(A) &\leq& \frac{ea^{N+1}\exp(N(1-a^{\prime}))}{\det(B)\left(a-\left(1+1/N \right)||B||\right)} \nonumber\\
&\leq& \left(\frac{e}{a^{\prime}-1-1/N}\right)\exp\left(1-a^{\prime}+\left(1+1/N\right)\ln\left(a^{\prime}\left(\frac{||B||^N}{\det(B)}\right)^{1/(N+1)}\right)\right)^N, 
\end{eqnarray} 
where $\displaystyle{a^{\prime}=\frac{a}{||B||}}$. Since $\displaystyle{\det(B)\geq ||B^{-1}||^{-N}}$ by Lemma \ref{firstlemma}, this implies 
\[ \mathbb{P}(A)\leq \left( \frac{e}{a^{\prime}-1-1/N}\right) \exp\left(1+\ln b-a^{\prime}+\left(1+1/N\right)\ln a^{\prime}\right)^{N}. \] 
But $\displaystyle{a^{\prime}-\left(1+1/N\right)\ln a^{\prime}\geq \frac{a}{b}-\left(1+1/N\right)\ln \frac{a}{b}}$, since the function $t\in \mathbb{R}_{+}^{\times}\mapsto t-\left(1+1/N\right)\ln t$ is increasing on $(1+1/N,+\infty)$ and $\displaystyle{1+1/N\leq \frac{a}{b}\leq a^{\prime}}$. Hence, 
\[ \mathbb{P}(A)\leq \left( \frac{e\cdot a}{a-2b}\right) \exp\left(1-\frac{a}{b}+\ln a\right)^N. \] 
\hfill $\Box$

\begin{rem} 
We have actually proved the stronger estimate \eqref{stronger} under the weaker assumption $\mathbb{E}_{A}(||Z||^2/N)\geq (1+1/N)||B||$. 
\end{rem} 

\begin{pro} \label{Bergmanbound} 
Say $\mathbf{X}=(X_{1},\ldots,X_{N})\sim \mathcal{N}_{\mathbb{C}^{N}}(0,B)$. Then, for all $x\in \mathbb{R}^{+}$ and all $I\subset \{1,\ldots,N\}$,
\[ \mathbb{P}\left(|X_{i}|^2<x^2 V(X_{i}) \,\, \mathrm{iff}\,\, i\in I\right) \leq \left(||C|| ||C^{-1}||\right)^{N} \mu(x)^{\# I} \left(1-\mu(x)\right)^{N-\# I}, \] 
where $\displaystyle{C=\left(\frac{Cov(X_{i},X_{j})}{\sqrt{V(X_{i})V(X_{j})}}\right)_{1\leq i,j\leq N} }$ and $\displaystyle{\mu(x)=1-e^{-x^2/||C||}}$.
\end{pro} 

\noindent \textit{Proof:} We observe that 
\begin{eqnarray*} 
\mathbb{P}\left(|X_{i}|^2<x^2 V(X_{i}) \,\, \mathrm{iff}\,\, i\in I\right) &=& \mathbb{P}\left( \frac{ |X_{i}|}{\sqrt{V(X_{i})}} <x  \,\, \mathrm{iff}\,\, i\in I \right) \\ 
&=& \int_{\prod_{i\in I} D(0,x) \times \prod_{i\notin I} D(0,x)^{c}}  \frac{e^{-\langle a,C^{-1}a\rangle}}{\pi^{N}\det C} da_{1}\ldots da_{N}, 
\end{eqnarray*} 
where $D(0,x)$ denotes the disk of radius $x$, centered at the origin of $\mathbb{C}$. 

Moreover,  by Lemma \ref{firstlemma}, $\displaystyle{\det A\geq ||A^{-1}||^{-N}}$ and for all $a\in \mathbb{C}^{N}$, $\displaystyle{\langle a,C^{-1}a\rangle \geq \frac{||a||^2}{||C||}}$, so that 
\begin{eqnarray*} 
\mathbb{P}\left(|X_{i}|^2<x^2 V(X_{i}) \,\, \mathrm{iff}\,\, i\in I\right) &\leq& ||C^{-1}||^N  \int_{\prod_{i\in I} D(0,x) \times \prod_{i\notin I} D(0,x)^{c}} \frac{e^{-||a||^2/||C||}}{\pi^{N}} da_{1}\ldots da_{N} \\ 
&\leq& \left(||C|| ||C^{-1}||\right)^{N} \mu(x)^{\# I} \left(1-\mu(x)\right)^{N-\# I}  
\end{eqnarray*} 
with $\displaystyle{\mu(x)=\mu\left(D\left(0,\frac{x}{\sqrt{||C||}}\right)\right)=1-e^{-x^2/||C||}}$ since $\mu \sim \mathcal{N}_{\mathbb{C}}(0,1)$.
\hfill $\Box$

We are now ready to prove the main result of this subsection. Let us denote by $\mathcal{B}(\mathbb{R}^{n})$ the $\sigma$-algebra of Lebesgue measurable sets of $\mathbb{R}^{N}$. 

\begin{thm} \label{probmain} 
Let $c_{N}\geq 1$, $\alpha>0$ and $\mathbf{X}=(X_{1},\ldots,X_{N})\sim \mathcal{N}_{\mathbb{C}^{N}}(0,B)$, $\mathbf{Y}=(Y_{1},\ldots,Y_{N})\in \mathbb{R}^{N}$ be two independent random vectors such that for every event $A\in \mathbf{Y}^{-1}(\mathcal{B}(\mathbb{R}^{N}))$, there exists $\mathbf{Z}=(Z_{1},\ldots,Z_{N})\sim\mathcal{N}_{\mathbb{C}^{N}}(0,B_{A})$  satisfying
\begin{itemize}
\item  $\forall i\in \{1,\ldots,N\}$, $\mathbb{E}(|Y_{i}|^2)\geq \alpha^2 \mathbb{E}(|Z_{i}|^2)$ and $\mathbb{E}_{A}(|Y_{i}|^2)\leq \mathbb{E}_{A}(|Z_{i}|^2)$, 
\item  $\forall I\subset \{1,\ldots,N\}$, $||C_{I}||||C_{I}^{-1}||\leq c_{N}$, with
\[ C_{I}=\left( \frac{Cov(Z_{i},Z_{j})}{\sqrt{V(Z_{i})V(Z_{j})}}\right)_{i,j\in I^{c}}.\]
\end{itemize}
Then, for all $x\in \mathbb{R}^{+}$ and all $\displaystyle{ y\geq \sqrt{2c_{N}}/\alpha}$, 
\begin{eqnarray*} 
1&-&\mathbb{P}\left(\exists i\in \{1,\ldots,N\}| |X_{i}|^2\geq x^2 V(X_{i}) \,\, \mathrm{and} |Y_{i}|^2\leq y^2 \mathbb{E}(|Y_{i}|^2)\right)  \\
&\leq& \left(\frac{e\alpha^2 y^2}{\alpha^2 y^2 -2c_{N}}\right) \left( b\left( \mu(x)+(1-\mu(x))\exp\left(1-\frac{\alpha^2 y^2}{c_{N}}+2\ln(\alpha y)\right)\right)\right)^{N},
\end{eqnarray*}  
where $b=||C|| ||C^{-1}||$ with  $\displaystyle{C=\left(\frac{Cov(X_{i},X_{j})}{\sqrt{V(X_{i})V(X_{j})}}\right)_{1\leq i,j \leq N}}$ and $\displaystyle{\mu(x)=1-e^{-x^2/||C||}}$. 
\end{thm} 

\noindent \textit{Proof:} 
We proceed as in the second proof of Proposition \ref{elementary} to deduce that
\begin{eqnarray*} 
1&-&\mathbb{P}\left(\exists i\in \{1,\ldots,N\}| |X_{i}|^2   \geq    x^2 V(X_{i}) \,\, \mathrm{and}\,\, |Y_{i}|^2\leq y^2 \mathbb{E}(|Y_{i}|^2)\right) \\
&=& \sum_{I\subset \{1,\ldots,N\}} \mathbb{P}\left(|X_{i}|^2<x^2V(X_{i})\,\,\mathrm{iff}\,\,i\in I\right) \mathbb{P}\left(\forall i\notin I, |Y_{i}|^2>y^2\mathbb{E}(|Y_{i}|^2)\right) \quad \mathrm{since} \, \mathbf{X}\perp \mathbf{Y}  \\ 
&=& b^{N} \sum_{I\subset \{1,\ldots,N\}} \mu(x)^{\#I}(1-\mu(x))^{N-\#I} \mathbb{P}(\forall i\notin I, |Y_{i}|^2>y^2\mathbb{E}(|Y_{i}|^2)) \quad \mathrm{by}\, \mathrm{Proposition} \, \ref{Bergmanbound}.
\end{eqnarray*} 
For every $I\subset \{1,\ldots,N\}$, let us denote by $A_{I}$ the event 
$\displaystyle{ \{\forall i\notin I, |Y_{i}|^2>y^2 \mathbb{E}(|Y_{i}|^2) \} }$. By hypothesis, there exists a complex centered Gaussian random vector $\mathbf{Z}_{I}=(Z_{1},\ldots,Z_{N})$ such that $\mathbb{E}(|Y_{i}|^2)\geq \alpha^2 \mathbb{E}(|Z_{i}|^2)$ and $\mathbb{E}_{A_{I}}(|Y_{i}|^2)\leq \mathbb{E}_{A_{I}}(|Z_{i}|^2)$ for all $i\in \{1,\ldots,N\}$.  
Thus, for all $i\notin I$, $\mathbb{E}_{A_{I}}(|Y_{i}|^2)>y^2 \mathbb{E}(|Y_{i}|^2)\geq \alpha^2 y^2 \mathbb{E}(|Z_{i}|^2)$ and therefore $\mathbb{E}_{A_{I}}(|Z_{i}|^2)>\alpha^2 y^2 \mathbb{E}(|Z_{i}|^2)=\alpha^2 y^2 V(Z_{i})$. We set 
\[ \widetilde{\mathbf{Z}}_{I}=\left( \frac{Z_{j_{1}}}{\sqrt{V(Z_{j_{1}})}},\ldots, \frac{Z_{j_{N-\# I}}}{\sqrt{V(Z_{j_{N-\# I}})}}\right),\]
where $\{j_{1},\ldots,j_{N-\# I}\}=I^{c}$, and deduce by summation that
$\displaystyle{\mathbb{E}_{A_{I}}\left( \frac{||\widetilde{\mathbf{Z}}_{I}||^2}{N-\# I} \right)> \alpha^2 y^2.}$ 
As soon as $I^{c}\neq \emptyset$, Proposition \ref{expboundprop} can be applied to $\widetilde{\mathbf{Z}}_{I}$ since by hypothesis $\alpha^2 y^2 \geq 2||C_{I}||||C_{I}^{-1}||:=2c_{I}$. It  provides
\[  \mathbb{P}(A_{I})\leq \left( \frac{e \alpha^2 y^2}{\alpha^2 y^2 -2c_{I}}\right)\exp\left(1- \frac{\alpha^2 y^2}{c_{I}}+2\ln(\alpha y)\right)^{N-\#I}, \]
the function $t\in \mathbb{R}_{+}^{\times}\mapsto t-c_{I}\ln(t)\in \mathbb{R}$ being increasing on $(c_{I},+\infty)$. The latter estimate remains valid when $I^{c}=\emptyset$, since $\displaystyle{ \frac{e \alpha^2 y^2}{\alpha^2 y^2 -2c_{I}} \geq 1}$.

Finally, setting $\displaystyle{c=\sup_{I\subset \{1,\ldots,N\}} ||C_{I}|| ||C_{I}^{-1}|| }$, we obtain
\begin{eqnarray*} 
& & 1-\mathbb{P}\left(\exists i\in \{1,\ldots,N\}| |X_{i}|^2\geq x^2 V(X_{i}) \,\, \mathrm{and} |Y_{i}|^2\leq y^2 \mathbb{E}(|Y_{i}|^2)\right)  \\
&\leq& \left( \frac{e\alpha^2 y^2}{\alpha^2 y^2-2c}\right) b^{N}\sum_{k=0}^{N} {N\choose k} \mu(x)^k 
\left((1-\mu(x))\exp\left(1-\frac{\alpha^2 y^2}{c}+2\ln(\alpha y)\right)\right)^{N-k} \\ 
&\leq& \left( \frac{e\alpha^2 y^2}{\alpha^2 y^2-2c}\right)  \left( b\left(\mu(x)+(1-\mu(x))\exp\left(1-\frac{\alpha^2 y^2}{c}+2\ln(\alpha y)\right)\right)\right)^N. 
\end{eqnarray*} 
Hence the result.  \hfill $ \Box$   

\begin{cor} \label{cormain} 
Under the hypotheses of Theorem \ref{probmain} and provided that $c_{N}\leq 2e/5$, for all $x\in \mathbb{R}_{+}$, 
\[ \mathbb{P}\left(\exists i\in \{1,\ldots,N\}| |X_{i}|^2\geq x^2 V(X_{i}) \,\, \mathrm{and} |Y_{i}|^2\leq e\frac{ \mathbb{E}(|Y_{i}|^2)}{\alpha^2} \right) \geq 1-5e\left(b\left( \mu(x)+(1-\mu(x))e^{-1/2}\right)\right)^{N}. \] 
\end{cor} 

\noindent\textit{Proof:}  We apply Theorem \ref{probmain} with $\displaystyle{y^2=\frac{5c_{N}}{2\alpha^2}}$, and deduce that for all $x\in \mathbb{R}_{+}$, 
\begin{eqnarray*} 
& &  \mathbb{P}\left(\exists i\in \{1,\ldots,N\}, |X_{i}|^2\geq x^2 V(X_{i}) \,\, \mathrm{and} |Y_{i}|^2\leq e\frac{ \mathbb{E}(|Y_{i}|^2)}{\alpha^2} \right) \\ 
&\geq&  \mathbb{P}\left(\exists i\in \{1,\ldots,N\}, |X_{i}|^2\geq x^2 V(X_{i}) \,\, \mathrm{and} |Y_{i}|^2\leq \frac{5c_{N}}{2}\frac{ \mathbb{E}(|Y_{i}|^2)}{\alpha^2} \right) \\ 
&\geq& 1-5e\left(b\left(\mu(x)+(1-\mu(x))\exp\left(-\frac{3}{2}+\ln\left(\frac{5c_{N}}{2}\right) \right) \right)\right)^{N} \\
&\geq&  1-5e\left(b\left( \mu(x)+(1-\mu(x))e^{-1/2}\right)\right)^{N}.
\end{eqnarray*} 
\hfill $\Box$

\subsection{Proof of Theorem \ref{exponential} } \label{proof} 
Let $d>0$ and $(B_{1},\partial^{-}B_{1},\partial^{+}B_{1}),\ldots, (B_{N},\partial^{-}B_{N},\partial^{+}B_{N})\subset \mathbb{CP}^2$ be degree $d$ marked Fubini-Study charts of capacity $\kappa$ centered at $z_{1},\ldots,z_{N}\in \mathbb{CP}^2$. By Definition \ref{markedFSchart}, there exist $\varphi_{1},\ldots, \varphi_{N}\in PU_{3}(\mathbb{C})$ such that for every $i\in \{1,\ldots,N\}$, $\varphi_{i}(\mathbb{B}_{d}(\kappa))=B_{i}$, $\varphi_{i}(\partial^{\pm} \mathbb{B}_{d}(\kappa))=\partial^{\pm} B_{i}$ and $\varphi_{i}(0)=z_{i}$. We set $P_{i}=P_{d}\circ \varphi_{i}^{-1}\in H^{0}(\mathbb{CP}^2,L^{d})$, where $P_{d}(Z,X,Y)=\sqrt{dN_{d}}YZ^{d-1}$ has been introduced in subsection \ref{linearpoly}, $H=Span(P_{1},\ldots, P_{N})$ and $H^{\perp}$ its orthogonal complement in $H^{0}(\mathbb{CP}^2,L^d)$. Let $C=\left( \langle P_{i},P_{j}\rangle \right)_{1\leq i,j\leq N}$ and $\displaystyle{\delta=\min_{i\neq j} d_{FS}(z_{i},z_{j})\geq \sqrt{\frac{20\ln(d)}{d}}}$ by hypothesis. By Lemma \ref{normest} and Proposition \ref{PolyFSestimates}, $\displaystyle{ c:=||C-I||\leq N ||C-I||_{\infty} \leq N d N_{d} \left( \frac{1+\cos(\delta)}{2}\right)^{d-1}}$ and $\displaystyle{ ||C^{-1}||\leq \frac{1}{1-c} }.$ We know by Lemmas  \ref{volcompare} and \ref{sinvolume} that $\displaystyle{ N\leq \frac{Vol(\mathbb{CP}^2)}{Vol(B_{FS}(0,\delta/2))}=\frac{1}{\sin^{4}(\delta/2)}}$, so that 
$\displaystyle{ N=O\left( \frac{d^2}{\ln(d)^2}\right) }$, 
\[ \left(\frac{1+\cos(\delta)}{2}\right)^{d-1} =\exp\left((d-1)\ln\left(1+\frac{\cos(\delta)-1}{2}\right)\right) = O\left(\frac{1}{d^{5}}\right) \] 
and $c=o(1)$. We deduce in particular that $P_{1},\ldots, P_{N}$ are linearly independent for $d$ large enough. For every $P\in H^{0}(\mathbb{CP}^2,L^d)$, we then denote by $X_{1}(P),\ldots,X_{N}(P)\in \mathbb{C}$ the coefficients of $\pi_{H}(P)$ in the basis $\{P_{1},\ldots, P_{N}\}$ of $H$ and for every $i\in \{1,\ldots, N\}$, we set $Y_{i}(P)=||P||_{B_{FS}(z_{i},\rho_{d})}$, $Y_{i,H}=Y_{i}\circ \pi_{H}$ and $Y_{i,H^{\perp}}=Y_{i}\circ \pi_{H^{\perp}}$, see subsection \ref{ballnorm}. 

We apply Corollary \ref{cormain} to $\mathbf{X}=(X_{1},\ldots, X_{N})$ and $\mathbf{Y}_{H^{\perp}}=(Y_{1,H^{\perp}},\ldots, Y_{N,H^{\perp}})$. The hypotheses of Theorem \ref{probmain} are satisfied since $\mathbf{X}$ is a centered complex Gaussian random variable independent of $\mathbf{Y}_{H^{\perp}}$ and every event $A$ in $\mathbf{Y}_{H^{\perp}}^{-1}(\mathcal{B}(\mathbb{R}^{N}))$ is invariant under $s_{H^{\perp}}$, so that by Lemma \ref{lemmaxA}, for every $i\in \{1,\ldots,N\}$, there exists $z_{i,A}\in B_{FS}(z_{i},\rho_{d})$ such that $\mathbb{E}_{A}(Y_{i,H^{\perp}}^2)\leq \mathbb{E}_{A}(||\widetilde{Z_{i}}||^2_{FS})$ with $\widetilde{Z_{i}}:P\in H^{0}(\mathbb{CP}^2,L^d)\mapsto P(z_{i,A})\in L^{d}|_{z_{i,A}}$. We set, for every $P\in H^{0}(\mathbb{CP}^{2},L^d)$, $\displaystyle{\widetilde{Z_{i}}(P)=Z_{i}(P)\frac{(e^{*})^{d}(z_{i,A})}{||(e^{*})^{d}(z_{i,A})||_{FS}}}$, where $e^{*}$ is the section of $L$ defined by $Z$, so that $Z_{i}$ is a complex centered  random variable which satisfies $V(Z_{i})=\mathbb{E}(||\widetilde{Z_{i}}||^2_{FS})=B(z_ {i,A},z_{i,A})=N_{d}$ by Lemmas \ref{bergman} and \ref{covariancebergman}. 
Moreover, for every $i\in \{1,\ldots,N\}$, 
\[ \mathbb{E}(Y_{i,H^{\perp}}^2)=N_{d}-\mathbb{E}(Y_{i,H}^2)\geq N_{d}-||C^{-1}|| \sum_{j=1}^{N} ||P_{j}||^2_{B_{FS}(z_{i},\rho_{d})} \]
by Lemmas \ref{directsum} and \ref{expectedestimate}, whereas by Corollary \ref{normestimate}, $\displaystyle{ ||P_{i}||^2_{B_{FS}(z_{i},\rho_{d})} \leq \frac{1}{\sin^{4}(\rho_{d})} }$ and for every $j\neq i$, $\displaystyle{ ||P_{j}||^2_{B_{FS}(z_{i},\rho_{d})} \leq d N_{d} \left( \cos(\delta-\rho_{d})\right)^{2(d-1)} }$. We observe that 
\[ \left(\cos(\delta-\rho_{d})\right)^{2(d-1)}=\exp\left(2(d-1)\ln(\cos(\delta-\rho_{d}))\right)=O\left(\exp(-d(\delta-\rho_{d})^2)\right)=O\left(\frac{\exp(4\sqrt{20\ln(d)})}{d^{20}}\right) \] 
and deduce that 
\[ \liminf_{d\rightarrow +\infty} \frac{1}{N_{d}} \mathbb{E}(Y_{i,H^{\perp}}^2) \geq 1- \frac{1}{\lim_{d\rightarrow+\infty} \left(N_{d}\sin^{4}(\rho_{d})\right)} =\frac{7}{8},\] 
see Remark \ref{44}. In particular, for $d$ large enough, $\mathbb{E}(Y_{i,H^{\perp}}^2)\geq \alpha^2 \mathbb{E}(||Z_{i}||^2_{FS})$, with $\alpha^2=3/4$. 
 Let $c_{N}=\max_{I\subset \{1,\ldots,N\}} ||C_{I}|| ||C_{I}^{-1}||$, with 
\[  C_{I}= \left( \frac{Cov(Z_{i},Z_{j})}{\sqrt{V(Z_{i})V(Z_{j})}}\right)_{i,j\in I} = \left(\frac{B(z_{i,A},z_{j,A})}{N_{d}}\right)_{i,j\in I} \] 
by Lemma \ref{covariancebergman}, so that
\begin{eqnarray*} 
c_{N} &\leq& \max_{I\subset \{1,\ldots,N\}} \left( \frac{1+||C_{I}-Id||}{1-||C_{I}-Id||}\right) \,\,\mathrm{by}\,\mathrm{Lemma}\, \ref{normest}  \\ 
&\leq& \max_{I\subset \{1,\ldots,N\}} \left( \frac{1+N||C_{I}-Id||_{\infty}}{1-N||C_{I}-Id||_{\infty}}\right) \,\,\mathrm{by}\,\mathrm{Lemma}\, \ref{normest}   \\ 
&\leq& \frac{1+N\cos^d(\delta-2\rho_{d})}{1-N\cos^d(\delta-2\rho_{d})}  \,\,\mathrm{by}\,\mathrm{Lemma}\, \ref{bergmanestimate}. 
\end{eqnarray*} 
Since 
\[ \cos^d(\delta-2\rho_{d})=\exp\left(d\ln(\cos(\delta-2\rho_{d}))\right) =O\left(\exp\left(-\frac{d}{2}(\delta-2\rho_{d})^2\right)\right) =   O\left( \frac{\exp(4\sqrt{20\ln(d)})}{d^{10}}\right),  \] 
we deduce that $c_{N}\xrightarrow[d\rightarrow +\infty]{} 1$.

In particular,  $\displaystyle{c_{N}\leq \frac{2e}{5}}$ for $d$ large enough and Corollary \ref{cormain} implies that for every $x\in \mathbb{R}_{+}$, 
\begin{equation} \label{fundamental} 
\mathbb{P}\left(\exists i\in \{1,\ldots,N\} | |X_{i}|^2\geq x^2 V(X_{i}) \, \mathrm{and} \, Y_{i,H^{\perp}}^2\leq \frac{e N_{d}}{\alpha^2} \right) \geq 1-5e\left(b\left(\mu(x)+(1-\mu(x))e^{-1/2}\right)\right)^N 
\end{equation}
with $\displaystyle{\mu(x)=1-e^{-x^2/||\widetilde{C}||}, \, \widetilde{C}=\left(\frac{Cov(X_{i},X_{j})}{\sqrt{V(X_{i})V(X_{j})}} \right)_{1\leq i,j\leq N} }$ and $b=||\widetilde{C}|| ||\widetilde{C}^{-1}||$, using that $\mathbb{E}(Y_{i,H^{\perp}}^2)\leq \mathbb{E}(Y_{i}^2)=N_{d}$ by Lemma  \ref{directsum}   . 

All the quantities $||\widetilde{C}||$,  $b$ and $(V(X_{i}))_{1\leq i\leq N}$ can be estimated by using Lemmas \ref{normest} and \ref{covariance}. Indeed, let $V=Diag_{1\leq i \leq N} (V(X_{i}))$ be the diagonal matrix which coincides with the diagonal part of $(Cov(X_{i},X_{j}))_{1\leq i,j\leq N}=\overline{C^{-1}}$ by Lemma \ref{covariance}. Then,

\begin{enumerate}
\item 
\begin{eqnarray*}
\frac{1}{||V^{-1}||}=\min_{1\leq i\leq N} V(X_{i})  &\geq& 1-||V-I||= 1-||V-I||_{\infty} \,\,\,\, (V\, \mathrm{is}\,\mathrm{diagonal}) \\ 
&\geq& 1-||C^{-1}-I||   \,\,\,\,(\mathrm{Lemma}\, \ref{normest})    \\
&\geq& 1-\frac{||C-I||}{1-||C-I||}= \frac{1-2c}{1-c}  \,\,\,\,(\mathrm{Lemma}\, \ref{normest}).
\end{eqnarray*} 

\item
\[ ||V||\leq 1+||V-I||\leq 1+ \frac{c}{1-c} =\frac{1}{1-c} \quad \mathrm{by} \, (1). \] 

\item
\begin{eqnarray*} 
b&=&||\sqrt{V}^{-1}\overline{C^{-1}}\sqrt{V}^{-1}|| || (\sqrt{V}^{-1}\overline{C^{-1}}\sqrt{V}^{-1})^{-1}|| 
 \leq ||V|| ||V^{-1}|| ||C|| ||C^{-1}|| \\ 
&\leq& \frac{1+c}{(1-c)(1-2c)} \,\, \mathrm{by} \, (1),\, (2)\, \mathrm{and}\, \mathrm{Lemma} \,\, \ref{normest}.
\end{eqnarray*} 

\item By (2) and Lemma  \ref{normest},
\[  \frac{1-c}{1+c} \leq \frac{1}{||C||||V||}\leq \frac{||C^{-1}||}{||V||}\leq ||\widetilde{C}||\leq ||C^{-1}|| ||V^{-1}|| \leq \frac{1}{1-2c} \cdot \]

\end{enumerate}

Now, by Corollary \ref{corexplicitsupupperbound}, if $\displaystyle{|X_{i}(P)|^2> \frac{c_{1}(d)}{\kappa^2 N_{d}} ||\pi_{H^{\perp}}(P)+\sum_{j\neq i}X_{j}(P)P_{j}||^2_{B_{FS}(z_{i},\rho_{d})}}$ for some $i\in \{1,\ldots,N\}$, then $(B_{i},\partial^{-}B_{i},\partial^{+}B_{i})$ is a submanifold chart for $V_{P}$. By the triangle inequality, 
\[ ||\pi_{H^{\perp}}(P)+\sum_{j\neq i}X_{j}(P)P_{j}||_{B_{FS}(z_{i},\rho_{d})} \leq Y_{i,H^{\perp}}(P)+\left(\sum_{j\neq i}|X_{j}(P)|\right)\max_{j\neq i} ||P_{j}||_{B_{FS}(z_{i},\rho_{d})} \] 
and by Cauchy-Schwarz inequality, 
\[ \left( \sum_{j\neq i} |X_{j}(P)|\right)^2  \leq N ||\mathbf{X}(P)||^2 =N\langle \pi_{H}(P),\mathbf{X}^{\star} \mathbf{X}(P)\rangle \leq N||\mathbf{X}^{\star}\mathbf{X}|| ||\pi_{H}(P)||^2, \] 
while $\displaystyle{ ||\mathbf{X}^{\star} \mathbf{X}||=||\mathbf{X} \mathbf{X}^{\star}||=||Cov(X_{i},X_{j})||=||\overline{C^{-1}}||\leq \frac{1}{1-c} }$ by  Lemmas \ref{normest} and \ref{covariance}. 

Under the condition $\displaystyle{ ||\pi_{H}(P)||^2\leq \frac{e N_{d}(1-c)}{\alpha^2 N \max_{j\neq i} ||P_{j}||^2_{B_{FS}(z_{i},\rho_{d})}} =R_{d}^2(N) }$, so that $\displaystyle{d^{16}=o(R_{d}^2(N))}$, we get 
\begin{eqnarray*} 
  ||\pi_{H^{\perp}}(P)+\sum_{j\neq i} X_{j}(P)P_{j}||  &\leq&  Y_{i,H^{\perp}}(P) +\sqrt{\frac{N}{1-c}} \max_{j\neq i} ||P_{j}||_{B_{FS}(z_{i},\rho_{d})} R_{d}(N) \\
&\leq&   Y_{i,H^{\perp}}(P)  +\frac{\sqrt{eN_{d}}}{\alpha} \cdot 
\end{eqnarray*} 
Choosing  $\displaystyle{x=\frac{2}{\alpha \kappa}\sqrt{ \frac{e \cdot  c_{1}(d)}{\min_{1\leq i\leq N} V(X_{i})} } }$ in \eqref{fundamental}, 
we thus deduce from Corollary  \ref{corexplicitsupupperbound} that 
\begin{eqnarray*} 
& & \mathbb{P}\left(\exists i\in \{1,\ldots,N\} \, \mathrm{such}\,\mathrm{that}\, (B_{i},\partial^{-}B_{i},\partial^{+}B_{i}) \, \mathrm{is}\,\mathrm{a}\,  \mathrm{submanifold}\,\mathrm{chart}\,\mathrm{for} \, V_{P} \right)   \\
&\geq& 1 - 5e\left(b\left(\mu(x)+(1-\mu(x))e^{-1/2}\right)\right)^{N}-\mathbb{P}\left(||\pi_{H}(P)||\geq R_{d}(N)\right)    \\
&\geq&  1-5e\left( \frac{(1+c)}{(1-c)(1-2c)} \left(1-\exp\left(-4e \frac{c_{1}(d)}{\alpha^2 \kappa^2}\left(\frac{1+c}{1-2c}\right)\right)(1-e^{-1/2})\right)\right)^{N} -\frac{R_{d}(N)^{2N}}{N^N} \exp(-R_{d}(N)^2+N) \\ 
\end{eqnarray*} 
by (1), (3) and Lemma \ref{lemma18}, since $R_{d}(N)>\sqrt{N}$ for $d$ large enough. Continuing the chain of inequalities, the quantity above dominates 
$\displaystyle{1-5e\gamma_{1}^{N}-\gamma_{2}^{N}=1-\gamma_{1}^{N} \left(5e+\left(\frac{\gamma_{2}}{\gamma_{1}}\right)^{N}\right)}$ 
with 
\[\gamma_{1}= \frac{1+c}{(1-c)(1-2c)} \left(1-\exp\left(-4e\frac{c_{1}(d)}{\alpha^2\kappa^2}\left(\frac{1+c}{1-2c}\right)\right)(1-e^{-1/2})\right) \] 
and 
\[ \gamma_{2}= \exp\left(1-\frac{R_{d}(N)^2}{N}+2\ln(R_{d}(N))-\ln N\right). \]

Now, $\gamma_{2}$ converges to $0$ as $d$ grows to $+\infty$, uniformly in $N$, and $\gamma_{1}$ converges to $1-\exp(-2^7\cdot 3^3\cdot e^{7/3}/(\alpha^2\kappa^2))(1-e^{-1/2})$ by Corollary \ref{explicitsupupperbound}, so that for instance there exists $d_{0}>0$ which does not depend on $N$ such that for every $d\geq d_{0}$, 
\[\mathbb{P}\left(\exists i\in \{1,\ldots,N\} \, \mathrm{such}\, \mathrm{that}\, (B_{i},\partial^{-}B_{i},\partial^{+}B_{i}) \, \mathrm{is} \, \mathrm{a}\, \mathrm{submanifold} \, \mathrm{chart} \, \mathrm{for} \, V_{P} \right)  
 \geq       1-14 \gamma(\kappa)^{N}                  
\] 
with $\gamma(\kappa)=1-\exp(-2^{11} 3^3/\kappa^2)(1-e^{-1/2})$. Hence the result. \hfill  $\Box$

\subsection{The lower estimate} 
We are now ready to prove Theorem \ref{lowerest} and  the lower estimate of Theorem \ref{mainbound}. 

\noindent \textit{Proof of Theorem \ref{lowerest}:}  Let $\kappa\in(0,1)$ and $d_{0}\in \mathbb{N}^{\times}$ be given by Theorem \ref{exponential} for this choice of $\kappa$. Let $d\geq d_{0}$, $\displaystyle{1>\delta\geq 3\sqrt{ \frac{5\ln(d)}{d}}}$ and $(x,y)\in (\mathbb{R}_{+}^{\times})^2$ be such that $\delta\leq x\leq 1$ and $\delta\leq y\leq 1$.            By Lemma  \ref{thereexistpoints}, there exists $z_{1},\ldots,z_{N}\in \mathbb{T}_{(x,y)}$ such that for all $i\neq j\in \{1,\ldots,N\}$, $\displaystyle{||z_{i}-z_{j}||\geq 3\sqrt{\frac{20\ln(d)}{d}}}$ with 
\[  N\geq  \left\lfloor \frac{2\delta \sqrt{d}}{3\sqrt{5\ln(d)}} \right\rfloor^2.\] 
By Lemma \ref{distancecompare}, it implies that for all $i\neq j$, $\displaystyle{d_{FS}(z_{i},z_{j})\geq \sqrt{\frac{20\ln(d)}{d}} }$. In particular, the Fubini-Study balls $(B_{FS}(z_{i},\rho_{d}))_{1\leq i\leq N}$ are pairwise disjoint, since $\displaystyle{\rho_{d}=\arctan(2/\sqrt{d}) \leq 2/\sqrt{d}}$. Let $\Delta=\{(z,\overline{z})|z\in \mathbb{C}\}$ be the graph of the complex conjugation in $\mathbb{C}^2$, it is a Lagrangian plane which contains the origin. For every $i\in \{1,\ldots,N\}$, let $\Psi_{i}\in PU_{3}(\mathbb{C})$ be such that $\Psi_{i}(0)=z_{i}$ and $d|_{0} \Psi_{i}(\Delta)=T_{z_{i}}\mathbb{T}_{(x,y)}$. We observe that $\Psi^{-1}(B_{FS}(z_{i},\rho_{d}))=B_{FS}(0,\rho_{d})$ and since $\rho_{d}\ll\min(x,y)$, $\Psi_{i}^{-1}(\mathbb{T}_{(x,y)})$ is a smooth perturbation of $\Delta$ in the neighborhood of this ball, that is 
\[ \Psi_{i}^{-1}(\mathbb{T}_{(x,y)})=\left\{ (z,\overline{z}+o(1/\sqrt{d} ) | |z|\leq 2/\sqrt{d} \right\}. \] 
The latter does not depend on $i\in \{1,\ldots,N\}$ provided $(\Psi_{i})_{1\leq i\leq N}$ is chosen such that for all $i\neq j$, $\Psi_{i}\circ \Psi_{j}^{-1}$ is a toric isometry, preserving $\mathbb{T}_{(x,y)}$. The chart $\mathbb{B}_{d}=D(0,1/\sqrt{6d})^2$ is contained in $B_{FS}(0,\rho_{d})$ by Lemma \ref{ballcontain} and contains $\mathbb{B}_{d}(\kappa)$. Moreover, $\Delta\cap \partial \mathbb{B}_{d}$ is a trivial knot which has linking number $\pm 1$ with the graph of any map $\partial D(0,1/\sqrt{6d})\rightarrow 
\mathring{D}(0,1/\sqrt{6d})$. We deduce that when $d$ is large enough $\Psi_{i}^{-1}(\mathbb{T}_{(x,y)})\cap \partial \mathbb{B}_{d}$ has linking number $\pm 1$ with the graph of any map $\partial D(0,1/\sqrt{6d})\rightarrow D(0,\kappa/\sqrt{6d})$, for the latter is homotopic to the graph of a constant map, the former homotopic to $\Delta\cap \partial \mathbb{B}_{d}$ and such homotopies can be chosen to keep the two knots disjoint.   

Now, for every $i\in \{1,\ldots,N\}$, let $(B_{i},\partial^{-}B_{i},\partial^{+}B_{i})=\Psi_{i}(\mathbb{B}_{d},\partial^{-}\mathbb{B}_{d},\partial^{+}\mathbb{B}_{d})$. By Theorem \ref{exponential}, the probability that there exists $i\in \{1,\ldots,N\}$ such that $(V_{P}\cap B_{i})$ is the image under $\Psi_{i}$ of the graph of a map $D(0,1/\sqrt{6d})\rightarrow D(0,\kappa/\sqrt{6d})$ is greater than $1-14\gamma(\kappa)^{N}$. We deduce that $V_{P}$ and $\mathbb{T}_{(x,y)}$ intersect each other in such a chart $B_{i}$, since the linking number of $\mathbb{T}_{(x,y)}\cap \partial B_{i}$ with $V_{P}\cap \partial B_{i}$ is $\pm 1$. Hence the result. \hfill $\Box$ 

\begin{rem} 
In Theorem \ref{lowerest}, the torus $\mathbb{T}_{(x,y)}$ could be replaced by any compact surface $S_{d}$ of $\mathbb{CP}^2$ with principal curvatures neglectable with respect to $\sqrt{d}$, so that one can pick disjoint marked degree $d$ FS-charts with centers in $S_{d}$ making them submanifold charts for $S_{d}$. 
\end{rem}

\begin{cor} \label{liminf} 
\[ \lim_{d\rightarrow +\infty} \inf_{(t_{1},t_{2})\in \mathbb{R}_{-}^2, -\frac{1}{2}\ln(d)+\ln(\ln(d))\leq t_{i} } (p(t))=1. \] 
\end{cor} 
Recall that $p(t)=\mathbb{P}(t\in \mathcal{A}(P))$, see subsection \ref{amoebameasures}. 

\noindent \textit{Proof:} We apply Theorem \ref{lowerest} with $\displaystyle{\kappa=\frac{1}{2}}$ and $\displaystyle{\delta=\frac{\ln(d)}{\sqrt{d}}}$ for $d$ large enough, so that $\ln(d)=O(N)$ and $N$ converges to $+\infty$ as $d$ grows to $+\infty$. This holds for every $(x,y)\in (\mathbb{R}_{+}^{\times})^2$ such that $\delta\leq x\leq 1$ and $\delta\leq y\leq 1$, that is $\displaystyle{ -\frac{1}{2}\ln(d)+\ln(\ln(d))\leq t_{1}\leq 0}$ and $\displaystyle{ -\frac{1}{2}\ln(d)+\ln(\ln(d))\leq t_{2}\leq 0}$, with $(t_{1},t_{2})=Log(x,y)$.  The result follows from the fact that $\displaystyle{\gamma\left(1/2\right)^{N} \xrightarrow[d\rightarrow +\infty]{} 0}$, while $p\leq 1$ by definition. \hfill $\Box$ 

\begin{figure}
\centering
\begin{tikzpicture}  
\draw[very thick]  (-2.8,0) -- (-2.8,-2.8);
\draw[very thick]  (-2.8,-2.8) -- (0,-2.8); 
\draw[very thick]  (0,-2.8) -- (2.8,0); 
\draw[very thick]   (2.8,0) -- (2.8,2.8); 
\draw[very thick]  (0,0) -- (2.8,2.8); 
\draw[very thick]  (0,2.8) -- (2.8,2.8); 
\draw[very thick]  (-2.8,0) -- (0,2.8);
\draw[very thick]  (-2.8,0) -- (0,0);
\draw[very thick]  (0,-2.8) -- (0,0);
\draw[gray, very thin, ->] (-3.5,0) -- (3.5,0); 
\draw[gray, very thin, ->] (0,-3.5) -- (0,3.5);
\node at (1.7,-2.9){$-\frac{1}{2}\ln(d)+\ln(\ln(d))$};
\node at (-4.3,0.3){$-\frac{1}{2}\ln(d)+\ln(\ln(d))$};
\node at (-1.4,-1.4){\Large $\mathcal{K}_{d}$};
\node at (1.4,0){\Large $\psi(\mathcal{K}_{d})$};
\node at (0,1.4){\Large $\psi^2(\mathcal{K}_{d})$};
\end{tikzpicture} 
\caption{The domain $\mathcal{H}_{d}$ in $\mathbb{R}^2$.} \label{figure1}
\end{figure}
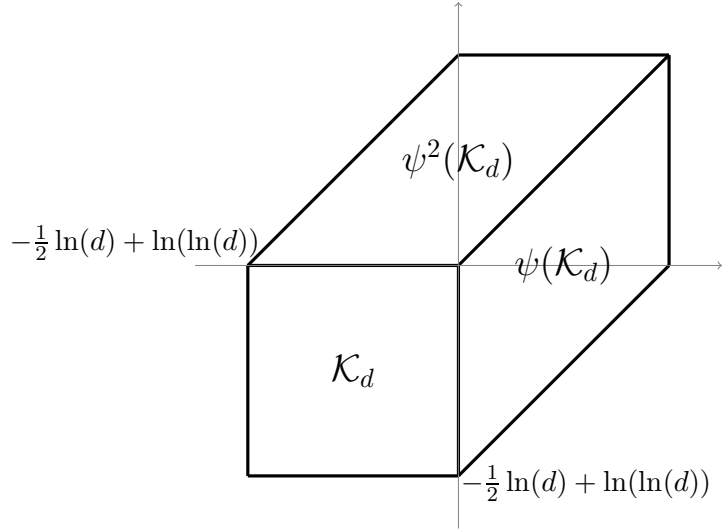

Let $\Psi:[z_{0}:z_{1}:z_{2}]\in \mathbb{CP}^2 \mapsto [z_{2}:z_{0}:z_{1}] \in \mathbb{CP}^2$, so that in the affine chart $U_{0}=\{[z_{0}:z_{1}:z_{2}]\in \mathbb{CP}^2|z_{0}=1\}$, it reads $\displaystyle{(z_{1},z_{2})\in \mathbb{C}\times \mathbb{C}^{\times} \mapsto \left(\frac{1}{z_{2}},\frac{z_{1}}{z_{2}}\right) \in \mathbb{C}^2}$. This toric isometry induces under the $Log$ map the special linear isomorphism 
$\displaystyle{ \psi: (t_{1},t_{2})\in \mathbb{R}^2 \mapsto (-t_{2}, t_{1}-t_{2})\in \mathbb{R}^2}$  
whose square equals  
$\displaystyle{ \psi^2 : (t_{1},t_{2})\in \mathbb{R}^2 \mapsto (t_{2}-t_{1},-t_{1})\in \mathbb{R}^2}$, so that $Log \circ \Psi = \psi \circ Log$.  
Let $\mathcal{H}_{d}$ be the union of the images of the square domain $\mathcal{K}_{d}=\displaystyle{ \{(t_{1},t_{2})\in \mathbb{R}^2| -\frac{1}{2}\ln(d)+\ln(\ln(d))\leq t_{i}\leq 0\} }$,  under the action of the group $\langle \psi \rangle =\{Id, \psi,\psi^2\}$, see Figure \ref{figure1}. The function $p$ is invariant under $\psi$ since for every $t\in \mathbb{R}^2$, 
\[ p\circ \psi (t)= \mathbb{P}(\psi(t)\in \mathcal{A}(P)) =\mathbb{P}(t\in \mathcal{A}(P\circ \Psi))=\mathbb{P}(t\in \mathcal{A}(Q))=p(t). \]  
We deduce from Corollary \ref{liminf} that it simply converges to $1$ over $\mathbb{R}^2$ as $d$ grows to $+\infty$, or also that it converges to 1 in the sense of distributions. 

\begin{cor} \label{distributions}
The expected amoeba measure $\mathbb{E}_{d}(\lambda)$ converges to the Lebesgue measure of $\mathbb{R}^2$ as $d$ grows to $+\infty$, in the sense of distributions. 
\end{cor} 
See subsection \ref{amoebameasures} for the definition of $\mathbb{E}_{d}(\lambda)$. 

\noindent \textit{Proof:}
Let $\varphi$ be a continuous function with compact support in $\mathbb{R}^2$. The latter gets contained in $\mathcal{H}_{d}$ for $d$ large enough, so that by Corollary \ref{liminf} and Proposition \ref{expectedmeasure}, 
\[ \langle \mathbb{E}_{d}(\lambda),\varphi\rangle = \int_{\mathbb{R}^2} \varphi(t)p(t)|dt|  \xrightarrow[d\rightarrow +\infty]{} \int_{\mathbb{R}^2} \varphi(t) |dt|. \] 
\hfill $\Box$ 

We finally deduce the asymptotic lower estimate given in Theorem \ref{mainbound}. 

\noindent \textit{Proof of the lower estimate of Theorem \ref{mainbound}:} For every $d$ large enough, 
\begin{eqnarray*} 
\frac{1}{\ln(d)^2} \mathbb{E}_{d}(Vol(\mathcal{A})) &=& \frac{1}{\ln(d)^2} \int_{\mathbb{R}^2} \langle \mathbb{E}_{d}(\lambda),1\rangle \\ 
&=& \frac{1}{\ln(d)^2} \int_{\mathbb{R}^2} p(t) |dt| \qquad \mathrm{by} \,\, \mathrm{Proposition}\,\,\mathrm{\ref{expectedmeasure}}\\
&\geq& \frac{1}{\ln(d)^2} \int_{\mathcal{H}_{d}} p(t) |dt| \qquad \mathrm{see}\,\,\mathrm{Figure} \,\,\ref{figure1}\\ 
&\geq& \frac{3}{\ln(d)^2} \iint_{[-\frac{1}{2}\ln(d)+\ln(\ln(d)),0]^2} p(t)|dt|, \quad \mathrm{since}\, \psi\in SL_{2}(\mathbb{R}). \\ 
\end{eqnarray*} 
The result follows, since the right hand side converges to $\displaystyle{\frac{3}{4}}$ by Corollary \ref{liminf}. \hfill $\Box$ 

\appendix

\section{Crofton's formula} \label{A}
Let $S$ be a smooth surface in $\mathbb{CP}^2$. The zero locus $V_{Q}$ of a section $Q\in H^{0}(\mathbb{CP}^2,L^{d})$ almost surely intersects $S$ transversely. We denote by $\delta_{V_{Q}\cap S}$ the empirical measure on this intersection locus, so that for every $\varphi\in C^{0}_{c}(S,\mathbb{C})$, ${\langle \delta_{V_{Q}\cap S},\varphi\rangle=\sum_{x\in V_{Q}\cap S} \varphi(x)}$. We then denote by $\mathbb{E}_{d}(\delta_{V_{Q}\cap S})$ the average measure, so that 
\[ \langle{E}_{d}(\delta_{V_{Q}\cap S}),\varphi\rangle =\int_{H^{0}(\mathbb{CP}^2,L^{d})}\left(\sum_{x\in V_{Q}\cap S} \varphi(x)\right) d\mu(Q), \] 
and finally we denote by $|dVol_{S}|$ the Lebesgue measure on $S$ induced by the Fubini-Study metric of $\mathbb{CP}^2$ restricted to $S$. 

\begin{thm}\label{crofton}
  Say $S$ is a smooth Lagrangian surface of $\mathbb{CP}^{2}$. Then, for all $d>0$, $\displaystyle{\mathbb{E}_{d}(\delta_{V_{Q}\cap S})=\frac{d}{2\pi}|dVol_{S}|}$.  
\end{thm} 
In particular, if $S$ has finite volume, the expected number of intersections of $V_{Q}$ and $S$ is given by 
\[ \mathbb{E}_{d}(\#(V_{Q}\cap S))=\frac{d}{2\pi}Area(S). \] 

\noindent \textit{Proof:}  Let $\mathcal{H}$ denote the correspondence subvariety $\mathcal{H}=\{ (x,Q)\in S\times H^{0}(\mathbb{CP}^2,L^d) |Q(x)=0 \}$, and $\pi_{1}, \pi_{2}$ the projections from $\mathcal{H}$ to the two factors $S$ and $H^{0}(\mathbb{CP}^2,L^d)$ respectively. The former is the restriction to $S$ of the fibration of hyperplanes in $H^{0}(\mathbb{CP}^{2}, L^{d})$ over $\mathbb{CP}^2$ with fibers $\mathcal{H}_{x}=\{Q\in H^{0}(\mathbb{CP}^2, L^{d})| Q(x)=0 \}$, $x\in \mathbb{CP}^{2}$. Then, for every $\varphi\in C^{0}_{c}(S,\mathbb{C})$, 
\begin{eqnarray*}
 \langle\mathbb{E}_{d}(\delta(V_{Q}\cap S)),\varphi\rangle &=& \int_{H^{0}(\mathbb{CP}^{2},L^d)} \left(\sum_{x\in(V_{Q}\cap S)}\varphi(x)\right)d\mu(Q) \\
&=&\int_{\mathcal{H}}(\pi_{1}^{*}\varphi) \pi_{2}^{*}\mu =\frac{1}{\pi^{N_{d}}} \int_{\mathcal{H}}e^{-||Q||^2_{FS}}(\pi_{1}^{*}\varphi) \pi_{2}^{*}|dVol(Q)|. 
\end{eqnarray*}
Writing $\pi_{2}^{*}|dVol(Q)|= (\pi_{1}^{*}|dVol_{S}(x)|)\wedge |dVol_{\mathcal{H}_{x}}(Q)| f(x,Q)$ for a suitable density function $f:\mathcal{H}\rightarrow \mathbb{R}_{+}^{\times}$, we deduce from Fubini's theorem that 
\begin{equation} \label{numberofintersections}
 \langle\mathbb{E}_{d}(\delta(V_{Q}\cap S)),\varphi\rangle= \frac{1}{\pi} \int_{S} \varphi(x)\left(\int_{\mathcal{H}_{x}}f(x,Q)d\mu(Q) \right) |dVol_{S}(x)|. 
\end{equation} 
Let us now calculate the density $f(x,Q)$. We denote by $\sigma_{x}$ a Bergman polynomial of norm 1 orthogonal to $\mathcal{H}_{x}$, so that $H^{0}(\mathbb{CP}^{n}, L^{d})=\mathcal{H}_{x}\oplus^{\perp} \mathbb{C} \sigma_{x}$ and $dVol=dVol|_{\mathcal{H}_{x}}\wedge dVol_{\mathbb{C}}$.  
Say $\gamma:t\in D(0,1)\mapsto (x(t),Q(t))\in \mathcal{H}$ is such that $\gamma(0)=(x,Q)$ and $d|_{(x,Q)} \pi_{2}(\gamma^{\prime}(0))\in \mathbb{C}\sigma_{x}$, where $D(0,1)$ is the unit disk in $\mathbb{C}$. For every $t\in D(0,1)$ we have $Q(t)(x(t))=0$, so that, by taking derivatives, we deduce 
\[ x^{\prime}(0)=-(dQ|_{T_{x}S})^{-1}\frac{\partial Q}{\partial t}|_{t=0}(x). \]
As a result, the restriction of $d|_{(x,Q)} \pi_{1}$ to $\mathcal{H}_{x}^{\perp}$ can be written as 
$\displaystyle{ t\in\mathbb{C}\mapsto -(d|_{x} Q|_{T_{x}S})^{-1}(t\sigma_{x}(x)) }$ 
and $\pi_{1}^{*}|dVol_{S}|=|dVol_{S}(d\pi_{1}, d\pi_{1})|=|Jac(d\pi_{1}|_{\mathcal{H}_{x}^{\perp}})||dVol_{\mathcal{H}_{x}^{\perp}}|$, showing that 
\begin{equation} \label{coareadensity}
f(x,Q)=\frac{1}{|det(dQ|_{T_{x}S})^{-1} \sigma_{x}(x)|}\cdot
\end{equation}
We now proceed to evaluate this density function in the case of an isotropically embedded surface. Let $\mathcal{H}_{2x}=\{Q\in H^{0}(\mathbb{CP}^{2}, L^{d})| Q(x)=0=d|_{x}Q\}$ and $\mathcal{H}_{2x}^{\perp}$ be the orthogonal complement of $\mathcal{H}_{2x}$ in $\mathcal{H}_{x}$. 

\begin{lem} \label{dilation} 
The evaluation $Q\in \mathcal{H}_{2x}^{\perp}\mapsto d|_{x}Q\in T_{x}^{*}\mathbb{CP}^{2}$ dilates the norm by a factor of $\displaystyle{\sqrt{ {d+2 \choose d-1,1,0}}}$. 
\end{lem} 

\noindent\textit{Proof:} It suffices to prove the result for $x=[1:0:0]$ since the isometry group of $\mathbb{CP}^2$ acts transitively on $\mathbb{CP}^2$.   The dimension of $\mathcal{H}_{2x}^{\perp}$ is two, since $d|_{x}Q=0$ imposes two independent conditions, and we choose the basis $\{X_{1},X_{2}\}$ for $\mathcal{H}_{2x}^{\perp}$ given by subsection \ref{FSmetric}, so that
\[ ||X_{i}||^2_{FS} = {d+2 \choose d-1, 1,0}^{-1}\]
by Lemma  \ref{FSnorm}. The result follows. \hfill $\Box$

\iffalse
\begin{eqnarray*} 
||x_{i}||_{FS}^{2} &=& \frac{1}{Vol_{FS}(\mathbb{CP}^{2})} \int_{\mathbb{C}^{2}} \frac{|x_{i}|^{2}}{(1+||z||^{2})^{d+3}}dz \\
&=& \frac{2\pi}{Vol_{FS}(\mathbb{CP}^{2})} \int_{\mathbb{C}} \left( \int_{0}^{\infty} \frac{ |x_{i}|^{2} rdr}{1+|x_{i}|^2+r^2)^{d+3}}\right) dx_{i} \\ 
&=& \frac{\pi}{Vol_{FS}(\mathbb{CP}^{2})(d+2)} \int_{\mathbb{C}} \frac{|x_{i}|^2 dx_{i}}{(1+|x_{i}|^2)^{d+2} } \\ 
&=& \frac{2\pi^{2}}{Vol_{FS}(\mathbb{CP}^{2})(d+2)} \int_{0}^{+\infty} \frac{r^3 dr}{(1+r^2)^{d+2}} 
\end{eqnarray*} 
Now, 
\begin{eqnarray*} 
\int_{0}^{+\infty} \frac{r^{3}dr}{(1+r^2)^{d+2}} &=& \int_{0}^{+\infty} \frac{r(r^2+1-1)dr}{(1+r^2)^{d+2}}\\
 &=& \int_{0}^{\infty} \frac{rdr}{(1+r^2)^{d+1}} -   \int_{0}^{\infty} \frac{rdr}{(1+r^2)^{d+2}} \\
&=& \frac{1}{2d}-\frac{1}{2d+2} = \frac{1}{d(2d+2)} 
\end{eqnarray*} 
Thus, 
\[ ||x_{i}||^2_{FS}=\frac{2\pi^{2}}{Vol_{FS}(\mathbb{CP}^{2})2(d+2)(d+1)d} = {d+2 \choose d-1, 1}^{-1}.\]
\fi

Let us now set $\mathcal{H}_{2x}^{S}=\{Q\in H^{0}(\mathbb{CP}^{n},L^d)| Q(x)=0 \, \mathrm{and}\, d|_{x}Q|_{T_{x}S}=0\}$ and write $\mathcal{H}_{x}=\mathcal{H}_{2x}^{S}\oplus^{\perp} (\mathcal{H}_{2x}^{S})^{\perp}$. 
The following corollary calculates $\int_{\mathcal{H}_{x}}f(x,Q)d\mu(Q)$ in the case of a Lagrangian surface $S$ and finishes the proof of Crofton's formula.
\begin{cor} \label{intlagrangian}
Let $S$ be a Lagrangian surface of $\mathbb{CP}^2$. Then, for every $x\in S$, $\int_{\mathcal{H}_{x}}f(x,Q)d\mu(Q)=d/2$.
\end{cor} 
\noindent\textit{Proof:} Since $S$ is isotropic, the condition $d|_{x}Q|_{T_{x}S}=0$ is a complex codimension two condition on the set of complex polynomials. Let $\{\sigma_{1},\sigma_{2}\}$ be an orthonormal basis of $(\mathcal{H}_{2x}^{S})^{\perp}$. If $Q=Q_{2}+a\sigma_{1}+b\sigma_{2}$ with $Q_{2}\in \mathcal{H}_{2x}^{S}$ and $a=a_{1}+ia_{2}, b=b_{1}+ib_{2}\in \mathbb{C}$, we obtain from equation \eqref{coareadensity} and Lemmas \ref{bergman} and \ref{dilation}
\[ f(x,Q)= \frac{1}{N_{d}}{d+2 \choose d-1,1,0}\begin{vmatrix} a_{1} & b_{1} \\ a_{2} & b_{2} \end{vmatrix}=d|Im(\overline{a}b)|. \] 
Consequently, 
\begin{eqnarray*} 
\int_{\mathcal{H}_{x}}f(x,Q)d\mu(Q)&=& d\int_{\mathbb{C}^2}|Im(\overline{a}b)|d\mu(a,b)\\
&=&\frac{d}{\pi^2}\int_{0}^{2\pi}\int_{0}^{2\pi}\int_{0}^{+\infty}\int_{0}^{+\infty} e^{-r_{1}^2-r_{2}^2}r_{1}^2 r_{2}^2|\sin(\theta_{2}-\theta_{1})|dr_{1}dr_{2}d\theta_{1}d\theta_{2} \\ 
&=& \frac{d}{4}\left(\frac{1}{\sqrt{\pi}}\int_{-\infty}^{+\infty}r^2 e^{-r^2}dr\right)^2 \frac{1}{\pi}\int_{0}^{2\pi}\int_{0}^{2\pi} |\sin\theta|d\theta_{1}d\theta \qquad \mathrm{(substitute}\, \theta=\theta_{2}-\theta_{1}\mathrm{)} \\ 
&=& \frac{d}{4} (-\cos\theta|_{0}^{\pi}) \\
&=& \frac{d}{2} \cdot 
\end{eqnarray*} 
\hfill $\Box$

\noindent\textit{Last part of the proof of Crofton's formula:} From Corollary \ref{intlagrangian} and equation \eqref{numberofintersections} we deduce that
\[ \mathbb{E}_{d}(\#(V_{Q}\cap S))= \frac{1}{\pi} \int_{S}\left(\int_{\mathcal{H}_{x}}f(x,Q)d\mu(Q) \right) |dVol_{S}(x)|=\frac{1}{\pi}\int_{S} \frac{d}{2}  |dVol_{S}(x)|=\frac{d}{2\pi}Area(S). \] 
\hfill $\Box$

\begin{rem} 
(1)  If $S$ is a complex curve instead of a Lagrangian surface, then $\displaystyle{\mathbb{E}_{d}(\delta_{V_{Q}\cap S})=\frac{d}{\pi}|dVol_{S}|}$, as follows from B\'{e}zout's Theorem, since the density function $\int_{\mathcal{H}_{x}}f(x,Q)d\mu(Q)$ is then also constant. In general, the latter depends on $x\in S$.

(2) Such a Crofton's formula has already been used by M. Kac or A. Edelman and E. Kostlan to compute the expected number of real roots of a univariate real polynomial  \cite{Kac, EK}.
\end{rem}

\section{Upper estimates in all degrees} 
\label{B}

Corollary \ref{upperlimit} provides the asymptotic upper bound on $\mathbb{E}_{d}(Vol(\mathcal{A}))$ as $d\rightarrow +\infty$. The aim of this Appendix \ref{B} is to prove the upper estimate given by Theorem \ref{mainbound} for every degree, see Theorem \ref{upperboundall}. For this purpose, we return to subsection \ref{levelset}. 

\begin{pro} \label{expupper}
Let  $d\geq 2$ and let $0<\rho_{0}<\sqrt{2}<\rho_{1}$ be the two positive solutions of the equation $\displaystyle{\frac{\rho^{2}}{(1+\rho^2)^{3/2}}=\frac{2}{\pi d}}$. Then, 
\begin{eqnarray*}
\mathbb{E}_{d}(Vol(\mathcal{A})) &\leq& 2\ln(d)\ln\left(\rho_{1}/\rho_{0}\right)-\ln(\rho_{1})^2-2\ln(\rho_{0})^2+\left(2\ln(\pi)+\pi \right)\ln\left(\rho_{1}/\rho_{0}\right) \\
&+& 3\ln(1+\rho_{0}^2)\ln(\rho_{0})+\frac{\pi^2 d}{2}\left(1-(1+\rho_{0}^2)^{-1/2}+(1+\rho_{1}^2)^{-1/2}\right). \\
\end{eqnarray*} 
\end{pro} 
\noindent \textit{Proof:} We first employ Corollary \ref{expbound} to get, for every $d>0$,  
\[ \mathbb{E}_{d}(Vol(\mathcal{A}))\leq \int_{\mathbb{R}^2}\min\left(1,da(t)/(4\pi)\right)|dt|=\int_{(\mathbb{R}_{+}^{\times})^2} \min\left(1,dA(x,y)/(4\pi)\right)\frac{dxdy}{xy}\cdot  \]
We then proceed as in the proof of Theorem \ref{upperboundasymptotic}      and split the latter integral into three regions that we estimate separately. 

(i) If $0\leq \rho \leq \rho_{0}$ or $\rho_{1}\leq \rho$, then $\displaystyle{\min\left(1,dA(x,y)/(4\pi)\right)=dA(x,y)/(4\pi)}$ and we saw in the proof of Theorem \ref{upperboundasymptotic} that
\[
\iint_{ \{\rho\leq \rho_{0}\}\cup \{\rho\geq \rho_{2}\}} \frac{d}{4\pi}A(x,y)\frac{dxdy}{xy} 
= \frac{\pi^2 d}{2} \left(1-(1+\rho_{0}^2)^{-1/2}+(1+\rho_{1}^2)^{-1/2}\right). 
\]

(ii) If $\rho_{0}\leq \rho\leq \rho_{1}$ and either $0\leq \theta \leq \theta_{\rho}$ or $\displaystyle{\pi/2-\theta_{\rho} \leq \theta \leq \pi/2}$, then $\displaystyle{\min\left(1,dA(x,y)/(4\pi)\right)=dA(x,y)/(4\pi)}$. Let us call this region $U$ temporarily. We have, as in the proof of Theorem \ref{upperboundasymptotic}, 
\begin{eqnarray*} 
\iint_{U} \frac{d}{4\pi}A(x,y) \frac{dxdy}{xy} 
&=& 4\int_{\rho_{0}}^{\rho_{1}} \frac{\theta_{\rho}}{\rho \sin(2\theta_{\rho})} d\rho  \\ 
&\leq&\pi \ln\left(\rho_{1}/\rho_{0} \right) \qquad \mathrm{using}\,\mathrm{the}\, \mathrm{inequality} \, \frac{2\theta_{\rho}}{\sin(2\theta_{\rho})}\leq \frac{\pi}{2}. 
\end{eqnarray*} 

(iii) If $\rho_{0}\leq \rho \leq \rho_{1}$ and $\displaystyle{\theta_{\rho}\leq \theta \leq \pi/2-\theta_{\rho}}$, then $\displaystyle{\min\left(1,dA(x,y)/(4\pi)\right)=1}$ and as in the proof of Theorem \ref{upperboundasymptotic}, 
\begin{eqnarray*} 
& &\int_{\rho_{0}}^{\rho_{1}} \int_{\theta_{\rho}}^{\pi/2-\theta_{\rho}} \frac{dxdy}{xy} = 2\ln\left( \frac{d\pi}{2}\right) \ln\left(\frac{\rho_{1}}{\rho_{0}}\right)+2(\ln(\rho_{1})^2-\ln(\rho_{0})^2) -3\int_{\rho_{0}}^{\rho_{1}} \ln(1+\rho^2)\frac{d\rho}{\rho}\\
& &+2 \int_{\rho_{0}}^{\rho_{1}} \ln(1+\cos(2\theta_{\rho}))\frac{d\rho}{\rho}  \\ 
&\leq& 2\ln(d\pi)\ln\left(\rho_{1}/\rho_{0}\right)+ 2(\ln(\rho_{1})^2-\ln(\rho_{0})^2)-3\int_{\rho_{0}}^{\rho_{1}} \ln(1+\rho^2) d\rho/\rho \quad \mathrm{since} \, \cos(2\theta_{\rho})\leq 1 \\ 
&\leq& 2\ln(d\pi) \ln\left(\rho_{1}/\rho_{0}\right)+2\ln(\rho_{1})^2-2\ln(\rho_{0})^2-3\int_{\rho_{0}}^{1}  \ln(1+\rho^2) d\rho/\rho - 3\int_{1}^{\rho_{1}} \ln(1+\rho^2) d\rho/\rho \\ 
&\leq&  2\ln(d\pi) \ln\left(\rho_{1}/\rho_{0}\right) +2\ln(\rho_{1})^2-2\ln(\rho_{0})^2+3\ln(1+\rho_{0}^2)\ln(\rho_{0})-3\ln(\rho_{1})^2 \\ 
&\leq &  2\ln(d\pi) \ln\left(\rho_{1}/\rho_{0}\right) - \ln(\rho_{1})^2-2\ln(\rho_{0})^2+3\ln(1+\rho_{0}^2)\ln(\rho_{0}).
\end{eqnarray*} 
Adding up the contributions from (i), (ii) and (iii) gives us the result.
\hfill $\Box$

We now estimate $\rho_{0}$ and $\rho_{1}$ in the following two lemmas. 
\begin{lem} \label{B2} 
For all $d\geq 2$, 
\[  \frac{2}{\pi d-3}\leq \rho_{0}^2 \leq \frac{4}{\left(\pi d-3\right)+\sqrt{\left(\pi d-3\right)^2-6}} \quad \mathrm{and} \quad \rho_{1}\leq \frac{\pi d}{2}\left(1-\frac{6}{(\pi d)^2-6}\right). \] 
\end{lem} 

\noindent \textit{Proof:}  From the relation $\displaystyle{\rho_{0}^2=(2/(\pi d))(1+\rho_{0}^2)^{3/2}}$ defining $\rho_{0}$, we first deduce that $\displaystyle{\rho_{0}^2\geq 2/\pi d}$. This in turn implies that
\[ \rho_{0}^{2}\geq \frac{2}{\pi d}\left(1+\frac{2}{\pi d}\right)^{3/2}\geq \frac{2}{\pi d}\left(1+\frac{3}{\pi d}\right), \] 
since $(1+t)^{3/2}\geq 1+3t/2$ for all $t\geq 0$. Iterating, we obtain
\begin{eqnarray*} 
\rho_{0}^{2}&\geq& \frac{2}{\pi d}\left(1+\frac{2}{\pi d}\left(1+\frac{3}{\pi d}\right) \right)^{3/2} \geq \frac{2}{\pi d}\left(1+\frac{3}{\pi d}\left(1+\frac{3}{\pi d}\right)\right), \quad \mathrm{so} \,\mathrm{that} \\
\rho_{0}^{2}&\geq& \frac{2}{\pi d}\left(1+\frac{2}{\pi d}\left(1+\frac{3}{\pi d}+\left(\frac{3}{\pi d}\right)^2\right)\right)^{3/2}\geq \frac{2}{\pi d}\left( 1+\frac{3}{\pi d}+\left(\frac{3}{\pi d}\right)^{2}+\left(\frac{3}{\pi d}\right)^{3} \right),\\
\end{eqnarray*} 
so that by recursion 
\[ \rho_{0}^{2}\geq \frac{2}{\pi d}\sum_{k=0}^{+\infty} \left(\frac{3}{\pi d}\right)^k = \frac{2}{\pi d}\left(\frac{1}{1-3/(\pi d)}\right)=\frac{2}{\pi d -3} \cdot\] 
Likewise, by definition, $\rho_{1}^{2}=(2/(\pi d))(1+\rho_{1}^2)^{3/2}$, so that $\rho_{1}^{-1}=(2/(\pi d))\left(1+\rho_{1}^{-2}\right)^{3/2}$. This implies 

\begin{eqnarray*}
\frac{1}{\rho_{1}}&=&\frac{2}{\pi d}\left(1+\left(\frac{2}{\pi d}\right)^2\left(1+\frac{1}{\rho_{1}^2}\right)^3\right)^{3/2} \\
&=& \frac{2}{\pi d}\left(1+\left( \frac{2}{\pi d}\right)^2\left(1+\left(\frac{2}{\pi d}\right)^2\left(1+\left(\frac{2}{\pi d}\right)^2  \left(1+\frac{1}{\rho_{1}^2}\right)^3                       \right)^3\right)^3\right)^{3/2} \\ 
&\geq& \frac{2}{\pi d}\left(1+\left( \frac{2}{\pi d}\right)^2\left(1+3\left(\frac{2}{\pi d}\right)^2\left(1+3\left(\frac{2}{\pi d}\right)^2    \left(1+\frac{1}{\rho_{1}^2}\right)^3               \right)\right)\right)^{3/2} \\
&\geq& \frac{2}{\pi d} \left( 1+\frac{3}{2}\left( \frac{2}{\pi d}\right)^2 \sum_{k=0}^{\infty} \left( 3 \left(\frac{2}{\pi d}\right)^2 \right)^{k} \right) \quad \mathrm{by}\, \mathrm{iteration.}
\end{eqnarray*} 
We deduce 
\[ \frac{1}{\rho_{1}}\geq \frac{1}{\pi d}\left(2+3\left(\frac{2}{\pi d}\right)^2\left( \frac{1}{1-3\left(2/(\pi d)\right)^2}\right)\right)\geq \frac{2}{\pi d}\left(\frac{(\pi d)^2-6}{(\pi d)^2-12}\right)\] 
and the upper bound on $\rho_{1}$. 

In order to obtain the upper bound on $\rho_{0}^2$, use the relation $\displaystyle{\rho_{0}^2=\frac{2}{\pi d}(1+\rho_{0}^2)^{3/2}}$ together with the inequality $(1+t)^{3/2}\leq 1+3t/2+3t^2/8$ valid for $t\geq 0$, to deduce 
\[ \rho_{0}^{2}\leq \frac{2}{\pi d}\left(1+\frac{3}{2}\rho_{0}^2+\frac{3}{8}\rho_{0}^4\right),\] 
that is 
\[ 1+\frac{3-\pi d}{2}\rho_{0}^2+\frac{3}{8}\rho_{0}^{4} \geq 0. \] 
The discriminant of the quadratic  polynomial obtained by substituting $u=\rho_{0}^2$ above equals $\displaystyle{\Delta=\left(\frac{\pi d-3}{2}\right)^2-\frac{3}{2}}$, which is positive for $d\geq 2$. Since $\rho_{0}^2\leq 2$, it must be less than or equal to the smaller root of the quadratic, hence 
\[ \rho_{0}^2\leq \frac{4}{3}\left(\frac{\pi d-3}{2}-\sqrt{\Delta}\right)=\frac{4}{\left(\pi d-3\right)+\sqrt{\left(\pi d-3\right)^2-6} }\cdot  \] 
\hfill $\Box$ 

An alternative to Lemma \ref{B2} is the following. 
\begin{lem} \label{iterate}
If $d\geq 2$, $0\leq \beta\leq \rho_{0}^2\leq \gamma$ and $2\leq \alpha\leq \rho_{1}^2$, then 
\begin{eqnarray*} 
\frac{2(3+\sqrt{9+(1+\beta^3)((\pi d)^2-12)})}{(\pi d)^2-12} &\leq& \rho_{0}^2  \leq \frac{2(3+\sqrt{9+(1+\gamma^3)((\pi d)^2-12)})}{(\pi d)^2-12} \quad \mathrm{and} \\
\frac{1}{2}\left(\left(\pi d/2\right)^2-3+\sqrt{\left(\left(\pi d/2\right)^2-3\right)^2-12-4/\alpha}\right) &\leq& \rho_{1}^2 \leq \frac{1}{2}\left(\left(\pi d/2\right)^2-3+\sqrt{\left(\left(\pi d/2\right)^2-3\right)^2-12}\right)
\end{eqnarray*} 
\end{lem} 

\noindent\textit{Proof:} Given $0\leq \beta\leq \rho_{0}^2 \leq \gamma$ and using the relation $\displaystyle{\rho_{0}^4=\left(2/(\pi d)\right)^2(1+\rho_{0}^2)^3}$, we deduce
\[ \left(\frac{2}{\pi d}\right)^2(1+\beta^3+3\rho_{0}^2+3\rho_{0}^4) \leq \rho_{0}^4 \leq \left(\frac{2}{\pi d}\right)^2(1+\gamma^3+3\rho_{0}^2+3\rho_{0}^4), \]
that is $\displaystyle{ \left((\pi d/2)^2-3\right)\rho_{0}^4-3\rho_{0}^2-(1+\beta^3)\geq 0}$ and  $\displaystyle{ \left((\pi d/2\right)^2-3)\rho_{0}^4-3\rho_{0}^2-(1+\gamma^3)\leq 0}$. The discriminants of these two trinomials in $\rho_{0}^2$ respectively are $\displaystyle{\Delta_{\beta}=9+(1+\beta^3)((\pi d)^2-12)}$ and \newline $\displaystyle{\Delta_{\gamma}=9+(1+\gamma^3)((\pi d)^2-12)}$. We deduce that $\rho_{0}^2$ has to be greater than the positive root of the trinomial on the left and less than the positive root of the trinomial on the right. This means that
\[ \frac{2(3+\sqrt{\Delta_{\beta}})}{(\pi d)^2-12} \leq \rho_{0}^2 \leq \frac{2(3+\sqrt{\Delta_{\gamma}})}{(\pi d)^2-12},\] 
which proves the first part of the statement. Likewise, 
\[ \rho_{1}^4=\left(2/(\pi d)\right)^2 (1+3\rho_{1}^2+3\rho_{1}^4+\rho_{1}^6) \geq \left(2/(\pi d)\right)^2(3\rho_{1}^2+3\rho_{1}^4+\rho_{1}^6), \]  
so that
\[ \rho_{1}^4+\left(3-(\pi d/2)^2\right)\rho_{1}^2+3\leq 0. \] 
The discriminant of the trinomial $\displaystyle{X^2+\left(3-(\pi d/2)^2\right)X+3}$ equals $\displaystyle{\Delta^{\prime}=\left((\pi d/2)^2-3\right)^2-12}$ and its roots are $\displaystyle{ y_{\pm}=\frac{1}{2}\left(\left(\pi d/2\right)^2-3\pm \sqrt{\Delta^{\prime}}\right)}$. We deduce that 
\[ \rho_{1}^2\leq y_{+}=\frac{1}{2}\left(\left(\pi d/2\right)^2-3+\sqrt{\left(\left(\pi d/2\right)^2-3\right)^2-12}\right). \] 
Finally,  if $\alpha\geq 2$ is a less than $\rho_{1}^2$, then $1\leq \rho_{1}^2/\alpha$ and 
\[ \rho_{1}^4=\left(2/(\pi d)\right)^2(1+3\rho_{1}^2+3\rho_{1}^4+\rho_{1}^6)\leq \left(2/(\pi d)\right)^2\left(\left(3+1/\alpha\right)\rho_{1}^2+3\rho_{1}^4+\rho_{1}^6\right), \] 
so that 
\[  \rho_{1}^4+\left(3-(\pi d/2)^2\right)\rho_{1}^2+\left(3+1/\alpha\right)\geq 0. \] 
The discriminant of the trinomial $\displaystyle{X^2+\left(3-(\pi d/2)^2\right)X+\left(3+1/\alpha\right)}$ is $\displaystyle{\Delta^{\prime\prime}=\left((\pi d/2)^2-3\right)^2-12-4/\alpha}$ and its roots are 
$\displaystyle{ z_{\pm}=\frac{1}{2}\left(\left(\pi d/2\right)^2-3\pm \sqrt{\Delta^{\prime\prime}}\right).}$
We deduce that for every $2\leq \alpha\leq \rho_{1}^2$,
\[ \rho_{1}^2\geq \frac{1}{2}\left(\left(\pi d/2\right)^2-3+\sqrt{\left(\left(\pi d/2\right)^2-3\right)^2-12-4/\alpha}\right).  \]
\hfill $\Box$

Starting from $\beta_{0}=0$, $\gamma_{0}=\alpha_{0}=2$, one can now use Lemma \ref{iterate} iteratively to improve the lower and upper bounds at each step. For small values of $d$, we obtain after few iterations the following numerical estimates: 
\begin{eqnarray*} 
d=2 \qquad 0.714 &\leq& \rho_{0}^2 \leq 0.715 \,\, \mathrm{and} \, \,   6.374\leq \rho_{1}^2\leq 6.401 \\ 
d=3 \qquad 0.322 &\leq& \rho_{0}^2\leq 0.323 \, \, \mathrm{and}\,\,  19.046\leq \rho_{1}^2 \leq 19.050 \\ 
d=4 \qquad 0.212 &\leq& \rho_{0}^2\leq 0.213 \, \, \mathrm{and} \, \, 36.395\leq \rho_{1}^2\leq 36.396 \\ 
d=5 \qquad 0.158&\leq&\rho_{0}^2\leq 0.159 \, \, \mathrm{and}\,\,  58.633 \leq \rho_{1}^2 \leq 58.634  \\
d=6 \,\, \quad 0.1269&\leq&\rho_{0}^2\leq 0.1270 \, \mathrm{and} \, 85.7913\leq \rho_{1}^2 \leq 85.7915 \\ 
\end{eqnarray*} 

We deduce from Proposition \ref{expupper} the upper bound $\mathbb{E}_{5}(Vol(\mathcal{A}))\leq 24.298$ for the expected area of the amoebas of plane complex curves in degree five, while the upper estimate given by half the expected multiarea of the amoeba is $\displaystyle{ \pi^2 d/2=24.675}$, see Corollary \ref{expbound} and Remark \ref{BKremark}.  In degrees less than five, the latter provides better estimates. As for degrees greater than five, the above values lead to $\mathbb{E}_{6}(Vol(\mathcal{A}))\leq 26.813$ and in general, 

\begin{cor} \label{datleast5}
For every $d\geq 6$, 
\[ 2/(\pi d)\leq \rho_{0}^2\leq 2/(\pi d)\cdot 1.2138 \quad \mathrm{and}\quad \left(\pi d/2\right)^2\cdot 0.9658 \leq \rho_{1}^2 \leq \left(\pi d/2\right)^2. \] 
\end{cor} 

\noindent \textit{Proof:} For all $d\geq 2$, by Lemma \ref{B2}, $\displaystyle{\rho_{0}^2\geq \frac{2}{\pi d-3} \geq \frac{2}{\pi d}}$ and by Lemma \ref{iterate} with $\gamma=0.1270$ for $d=6$ 
\begin{eqnarray*} 
\rho_{0}^2 &\leq& \frac{2}{\pi d}\left(\frac{ 3/(\pi d)+\sqrt{9/(\pi d)^2+(1+\gamma^3)}}{1-12/(\pi d)^2}\right) \\ 
&\leq& \frac{2}{\pi d}\cdot 1.2138 
\end{eqnarray*} 
Likewise, by Lemma \ref{B2}, $\displaystyle{\rho_{1}^2\leq \left(\pi d/2\right)^2}$ and by Lemma \ref{iterate} with $\alpha=85.7913$ for $d=6$,
\begin{eqnarray*} 
\rho_{1}^2&\geq& \frac{1}{2}\left(\frac{\pi d}{2}\right)^2 \left( 1-\frac{3}{\left(\pi d/2\right)^2}+\sqrt{\left(1-\frac{3}{\left(\pi d/2\right)^2}\right)^2-\frac{12+4/\alpha}{\left(\pi d/2\right)^4}}\right) \\ 
&\geq& \left(\frac{\pi d}{2}\right)^2 \cdot 0.9658
\end{eqnarray*} 
\hfill $\Box$

\begin{thm} \label{upperboundall} 
For all $d\geq 6$, $\displaystyle{\mathbb{E}_{d}(Vol(A))\leq \frac{3}{2}\ln(d)^2+8.3752\ln(d)+8.5607}$.
\end{thm} 

\noindent \textit{Proof:} By  Corollary \ref{datleast5}, $\displaystyle{ \ln\left(\rho_{1}/\rho_{0}\right)\leq \frac{3}{2}\ln\left(\pi d/2\right)}$. Using this estimate together with the other estimates in Corollary \ref{datleast5}, we deduce from Proposition \ref{expupper} that 
\begin{eqnarray*} 
& &\mathbb{E}_{d}(Vol(\mathcal{A})) \leq 3\ln(d)\ln\left(\pi d/2\right)-\left( \frac{1}{2}\ln\left(\left(\pi d/2\right)^2\cdot 0.9658\right)\right)^2-2\left(\frac{1}{2}\ln\left(\frac{2}{\pi d}\cdot 1.2138 \right)\right)^2 \\
& & +\frac{3}{2}\left(2\ln(\pi)+\pi\right)\ln\left(\pi d/2\right)+\frac{\pi^2 d}{2}\left(1-(1+2/(\pi d)\cdot 1.2138)^{-1/2}+(1+\left(\pi d/2\right)^2\cdot 0.9658)^{-1/2} \right) \\
&\leq & 3\ln(d)^2+3\ln\left(\pi/2\right)\ln(d)-\left(\ln(d)+\ln\left(\pi/2\right)+\frac{1}{2}\ln(0.9658)\right)^2\\
& &-\frac{1}{2}\left(\ln(d)+\ln\left(\pi/2\right)-\ln(1.2138)\right)^{2}+8.1466\left(\ln\left(\pi/2\right)+\ln(d)\right)+\frac{\pi^2 d}{2}\left(1-\left(1-\frac{1.2138}{\pi d}\right)+\frac{2.0352}{\pi d}\right) \\ 
&\leq& \frac{3}{2}\ln(d)^2+\left(\ln(1.2138)-\ln(0.9658)+8.1466\right)\ln(d)-\left(\ln(\pi/2)+\frac{1}{2}\ln(0.9658)\right)^2\\
& &-\frac{1}{2}\left(\ln\left(\pi/2\right)-\ln(1.2138)\right)^{2}+8.1466\ln\left(\pi/2\right)+1.6245\pi \\
&\leq& \frac{3}{2}\ln(d)^2+8.3752\ln(d)+8.5607.
\end{eqnarray*} 
\hfill $\Box$ 
\begin{rem} \label{last}
For $d=6$, Theorem \ref{upperboundall} provides the upper bound $\mathbb{E}_{6}(Vol(\mathcal{A}))\leq 28.3827$, while $\pi^2 d/2=29.6089$. The upper estimate given by Theorem \ref{mainbound} follows from Theorem \ref{upperboundall} and the above given values for $d\leq 5$. 
\end{rem}

\end{document}